\documentclass[12pt]{article}
\usepackage[latin1]{inputenc}
\usepackage{amsmath}
\usepackage{amsfonts}
\usepackage{amssymb}
\usepackage{amsthm}

\usepackage{fixltx2e} % subscripts

\usepackage{titlesec}
\titleformat{\subsection}[hang]{\normalfont\bfseries}{\thesubsection}{1em}{}
\titlespacing\section{0pt}{3.5ex plus 0.5ex minus .2ex}{0.3ex plus .2ex}
\titlespacing\subsection{0pt}{2.5ex plus 0.5ex minus .2ex}{0.3ex plus .2ex}
\titlespacing\subsubsection{0pt}{2.5ex plus 0.5ex minus .2ex}{0.3ex plus .2ex}
\usepackage{mathrsfs} % for \mathscr (script letters)

\usepackage[alphabetic,lite]{amsrefs} % for different reference style

\usepackage{fullpage}
\usepackage{setspace}
\usepackage{hyperref}
\usepackage{color}
\usepackage{enumitem}
\usepackage{ulem} %um text durchzustreichen
\usepackage{comment}

\usepackage{marginnote}
\usepackage[all,cmtip]{xy} % crazy diagrams

 \usepackage{fancyhdr}
\newtheorem{Thm}{Theorem}[section]
\newtheorem{Lemmasub}{Lemma}[Thm]

\newtheorem{Lemma}[Thm]{Lemma}

\newtheorem{Cor}[Thm]{Corollary}

\newtheorem{Prop}[Thm]{Proposition}

\theoremstyle{definition}
\newtheorem{Def}[Thm]{Definition}
\newtheorem{Rem}[Thm]{Remark}
\newtheorem{Aside}[Thm]{Aside}
\newtheorem{Assumption}[Thm]{Assumption}
\newcommand{\Proof}{\textbf{Proof.\\}}

\newcommand{\<}{\left\langle}
\renewcommand{\>}{\right\rangle}
\newcommand{\bC}{\mathbb{C}}

\newcommand{\bF}{\mathbb{F}}
\newcommand{\bG}{\mathbb{G}}

\newcommand{\bP}{\mathbb{P}}
\newcommand{\bQ}{\mathbb{Q}}
\newcommand{\bR}{\mathbb{R}}
\newcommand{\bS}{\mathbb{S}}

\newcommand{\bZ}{\mathbb{Z}}

\newcommand{\cO}{\mathcal{O}}
\newcommand{\cP}{\mathcal{P}}

\newcommand{\cR}{\mathcal{R}}
\newcommand{\cS}{\mathcal{S}}

\newcommand{\fg}{\mathfrak{g}}
\newcommand{\fh}{\mathfrak{h}}

\newcommand{\fj}{\mathfrak{j}}

\newcommand{\fr}{\mathfrak{r}}
\newcommand{\fs}{\mathfrak{s}}
\newcommand{\ft}{\mathfrak{t}}

\newcommand{\fI}{\mathfrak{I}}

\newcommand{\sA}{\mathscr{A}}
\newcommand{\sB}{\mathscr{B}}

\newcommand{\sS}{\mathscr{S}}

\DeclareMathAlphabet{\mathpzc}{OT1}{pzc}{m}{it}

\newcommand{\ra}{\rightarrow}

\newcommand{\wt}{\widetilde}

\newcommand{\eps}{\epsilon}

\newcommand{\ov}{\overline}

\providecommand{\abs}[1]{\left\lvert#1\right\rvert}

\DeclareMathOperator{\Spec}{Spec}		% Spectrum of a ring
\DeclareMathOperator{\Hom}{Hom}			% Set of arrows between two object
\DeclareMathOperator{\Id}{Id}		  	% Identity
\DeclareMathOperator{\GL}{GL}		  	% GL
\DeclareMathOperator{\SL}{SL}		  	% SL
\DeclareMathOperator{\Dyn}{Dyn}			% Dynkin diagram
\DeclareMathOperator{\Gal}{Gal}			% Galois group
\DeclareMathOperator{\Sp}{Sp}		  	% Sp (for rigid analytic geometry)
\DeclareMathOperator{\Ad}{Ad}	  		% Adjoint
\DeclareMathOperator{\Lie}{Lie}  	  % Lie algebra
\DeclareMathOperator{\ind}{ind}  	  % compact induction
\DeclareMathOperator{\val}{val}  	  % valuation
\DeclareMathOperator{\Cent}{Cent}  	  % centralizer
\newcommand{\ba}{\begin{aligned}}
\newcommand{\ea}{\end{aligned}}
\newcommand{\beqn}{\begin{eqnarray}}
\newcommand{\eeqn}{\end{eqnarray}}
\newcommand{\beqns}{\begin{eqnarray*}}
\newcommand{\eeqns}{\end{eqnarray*}}
\newcommand{\benum}{\begin{enumerate}}
\newcommand{\eenum}{\end{enumerate}}

\newcommand{\RP}{\mathbf{G}} %reductive quotient of parahoric subgroup

\def\almoststronglystable/{almost strongly stable}
\def\Almoststronglystable/{Almost strongly stable}
\def\almoststable/{almost stable}
\def\extendeddatum/{extended datum}
\def\datum/{datum}
\def\extendeddata/{extended data}
\def\data/{data}
\def\truncatedextendeddatum/{truncated extended datum}
\def\truncateddatum/{truncated datum}
\def\generic/{generic}
\def\stype/{$\fs$-type}
\def\type/{type}

\newcommand{\ff}{\mathfrak f}  	 
\newcommand{\Uff}{U_{\ff}}  
\newcommand{\Uiff}{U_{i,\ff}}

\newcommand{\repzero}{\rho_0}  	
\newcommand{\Vrepzero}{V_{\rho_0}}  	
\newcommand{\repzeroK}{\rho}  	
\newcommand{\repzeroKYu}{\wt \rho}  	
\newcommand{\VrepzeroK}{V_{\rho}}  
\newcommand{\VrepzeroKYu}{V_{\wt \rho}}  	
\newcommand{\piK}{\pi_K}  	
\newcommand{\VpiK}{V_{\pi_K}}  
\newcommand{\piKYu}{\pi_{\wt K}}  	
\newcommand{\VpiKYu}{V_{\pi_{\wt K}}}  	
\newcommand{\rhoYu}{\rho_{\text{Yu}}}  	 
\newcommand{\KYu}{\wt K}

\usepackage{titling}

\AtEndDocument{\bigskip{\footnotesize%
		\textsc{Trinity College, Cambridge, CB2 1TQ, UK} \par  
		\textit{E-mail address}: \texttt{fintzen@math.duke.edu}
	}}

\begin{document}
\author{Jessica Fintzen} 
\title{Types for tame $p$-adic groups } 
\date{\vspace{-1.5cm}}

\maketitle  
\begin{abstract}
	Let $k$ be a non-archimedean local field with residual characteristic $p$. Let $G$ be a connected reductive group over $k$ that splits over a tamely ramified field extension of $k$. Suppose $p$ does not divide the order of the Weyl group of $G$. Then we show that every smooth irreducible complex representation of $G(k)$ contains an $\fs$-type of the form constructed by Kim--Yu and that every irreducible supercuspidal representation arises from Yu's construction.
	This improves an earlier result of Kim, which held only in characteristic zero and with a very large and ineffective bound on $p$. By contrast, our bound on $p$ is explicit and tight, and our result holds in positive characteristic as well. Moreover, our approach is more explicit in extracting an input for Yu's construction from a given representation. \\[-1.2cm]
\end{abstract}

	{
		\renewcommand{\thefootnote}{}  % to delete the footnote number 	
	%	\footnotetext{Date: Nov 1, 2020}
		\footnotetext{MSC2010: 22E50} 
		\footnotetext{Keywords: representations of reductive groups over non-archimedean local fields, types, supercuspidal representations, $p$-adic groups} 
		\footnotetext{The author was partially supported by a postdoctoral fellowship of the German Academic Exchange Service (DAAD), an AMS--Simons travel grant and NSF Grants DMS-1638352 and DMS-1802234.}
	}

\tableofcontents

\section{Introduction}
The aim of the theory of types is to classify, up to some natural equivalence, the smooth irreducible complex representations of a $p$-adic group in terms of representations of compact open subgroups. For $\GL_n$ it is known that every irreducible representation contains an $\fs$-type. This theorem lies at the heart of many results in the representation theory of $\GL_n$ and plays a key role in the construction of an explicit local Langlands correspondence for $\GL_n$ as well as in the study of its fine structure. One of the main results of this paper is the existence of $\fs$-types for general $p$-adic groups and the related exhaustion of supercuspidal representations under minimal tameness assumptions. These tameness assumptions arise from the nature of the available constructions of supercuspidal representations for general $p$-adic groups.

To explain our results in more detail, let $k$ denote a non-archimedean local field with residual characteristic $p$ and let $G$ be a connected reductive group over $k$. Before introducing the notion of a \type/, let us first discuss the case of supercuspidal representations, the building blocks of all other representations. Since the constructions below of supercuspidal representations for general reductive groups $G$ assume that $G$ splits over a tamely ramified extension of $k$, we will impose this condition from now on.
Under this assumption, Yu (\cite{Yu}) gave a construction of supercuspidal representations as representations induced from compact mod center, open subgroups of $G(k)$ generalizing an earlier construction of Adler (\cite{Adler}). Yu's construction is the most general construction of supercuspidal representations for general reductive groups known at present and it has been widely used to study representations of $p$-adic groups, e.g. to obtain results about distinction, to calculate character formulas, to suggest an explicit local Langlands correspondence and to investigate the theta correspondence. However, all these results only apply to representations obtained from Yu's construction.
In this article, we prove that  \textit{all} supercuspidal representations of $G(k)$ are obtained from Yu's construction if $p$ does not divide the order of the Weyl group $W$ of $G$. This result was previously shown by Kim (\cite{Kim}) under the assumption that $k$ has characteristic zero and that $p$ is ``very large''. Note that Kim's hypotheses on $p$ depend on the field $k$ and are much stronger than our requirement that $p \nmid \abs{W}$, see \cite[\S~3.4]{Kim}. The few primes that divide the order of the Weyl group of $G$ are listed in Table \ref{table-weyl-group}, and we expect that this assumption is optimal in general when also considering types as below for the following reason. Yu's construction is limited to tori that split over a tamely ramified field extension of $k$. If $p$ does not divide the order of the Weyl group of $G$ and $G$ splits over a tamely ramified extension (our assumptions), then all tori split over a tame extension. However, if one of these assumptions is violated, then, in general, the group $G$ contains tori that do not split over a tame extension (for some non-split inner forms of split groups of type $A_n, n \geq 2, D_l, l \geq 4$ prime, or $E_6$ the condition on the prime number is slightly weaker, see \cite[Theorem~2.4 and Corollary~2.6]{Fi-tame-tori} for the details). We expect that we can use these tori to produce supercuspidal representations (of Levi subgroups) that were not constructed by Yu. Examples of such representations are provided by the construction of Reeder and Yu (\cite{ReederYu}), whose ingredients exist also when $p \mid \abs{W}$ (whenever they exist for some large prime $p$), see \cite{FR} and \cite{Fi}.

In order to study arbitrary smooth irreducible representations, we recall the theory of types introduced by Bushnell and Kutzko (\cite{BK-types}): By Bernstein (\cite{Bernstein}) the category $\cR(G)$ of smooth complex representations of $G(k)$ decomposes into a product of subcategories $\cR^{\fs}(G)$ indexed by the set of inertial equivalence classes $\fI$ of pairs $(L, \sigma)$ consisting of a Levi subgroup $L$ of (a parabolic subgroup of) $G$ together with a smooth irreducible supercuspidal representation $\sigma$ of $L(k)$: $$\cR(G)=\prod_{\fs \in \fI}\cR^{\fs}(G). $$ 

Let $\fs \in \fI$. Following Bushnell--Kutzko (\cite{BK-types}), we call a pair $(K, \rho)$ consisting of a compact open subgroup $K$ of $G(k)$ and an irreducible smooth representation $\rho$ of $K$ an \textit{$\fs$-type} if for every irreducible smooth representation $\pi$ of $G(k)$ the following holds: 
\begin{center} $\pi$ lies in $\cR^{\fs}(G)$ if and only if $\pi|_K$ contains $\rho$. \end{center}
In this case the category $\cR^\fs(G)$ is isomorphic to the category of (unital left) modules of the Hecke algebra of compactly supported $\rho$-spherical functions on $G(k)$. 
 Thus, if we know that there exists an \stype/ for a given $\fs \in \fI$, then we can study the corresponding representations $\cR^{\fs}(G)$ using the corresponding Hecke algebra.
We say that a smooth irreducible representation $(\pi, V_\pi)$ of $G(k)$ \textit{contains a \type/} if there exists an \stype/ $(K, \rho)$ for the class $\fs \in \fI$ that satisfies $(\pi, V_\pi) \in \cR^{\fs}(G)$, i.e. $\pi|_{K}$ contains $\rho$.

Using the theory of $G$-covers introduced by Bushnell and Kutzko in \cite{BK-types}, Kim and Yu (\cite{KimYu}) showed that Yu's construction of supercuspidal representations can also be used to obtain types by omitting some of the conditions that Yu imposed on his input data. In this paper we prove that every smooth irreducible representation of $G(k)$ contains such a \type/ if $k$ is a non-archimedian local field of \textit{arbitrary} characteristic whose residual characteristic $p$ does not divide the order of the Weyl group of $G$. This excludes only a few residual characteristics, and we expect the restriction to be optimal in general as explained above. If $k$ has characteristic zero and $p$ is ``very large'', then Kim and Yu deduced this result already from Kim's work (\cite{Kim}).

Our approach is very different from Kim's approach. While Kim proves statements about a measure one subset of all smooth irreducible representations of $G(k)$ by matching summands of the Plancherel formula for the group and the Lie algebra, we use a more explicit approach involving the action of one parameter subgroups on the Bruhat--Tits building. This means that even though we have formulated some statements and proofs as existence results, the interested reader can use our approach to extract the input for the construction of a type 
 from a given representation.  

To indicate the rough idea of our approach, we assume from now that $p$ does not divide the order of the Weyl group of $G$, and we denote by $(\pi, V_\pi)$ an irreducible smooth representation of $G(k)$.
Recall that Moy and Prasad (\cites{MP1, MP2}) defined for every point $x$ in the Bruhat--Tits building $\sB(G,k)$ of $G$ and every non-negative real number a compact open subgroup $G_{x,r} \subset G(k)$ and a lattice $\fg_{x,r} \subset \fg$ in the Lie algebra $\fg=\Lie(G)(k)$ of $G$ such that $G_{x,r} \trianglelefteq G_{x,s}$ and $\fg_{x,r} \subseteq \fg_{x,s}$ for $r>s$.
Moy and Prasad defined the \textit{depth} of $(\pi, V_\pi)$ to be the smallest non-negative real number $r_1$ such that there exists a point $x \in \sB(G,k)$ so that the space of fixed vectors $V_\pi^{G_{x,r_1+}}$ under the action of the subgroup $G_{x,r_1+}:=\bigcup_{s>r_1}G_{x,s}$ is non-zero. In \cite{MP2} they showed that every irreducible depth-zero representation contains a \type/. A different proof using Hecke algebras was given by Morris (\cite{Morris-depth-zero}, announcement in \cite{Morris}). More generally, Moy and Prasad showed that $(\pi, V_\pi)$ contains an unrefined minimal K-type, and all unrefined minimal K-types are associates of each other. For $r_1=0$, an  \textit{unrefined minimal K-type} is a pair $(G_{x,0}, \chi)$, where $\chi$ is a cuspidal representation of the finite (reductive) group $G_{x,0}/G_{x,0+}$. If $r_1 >0$, then an unrefined minimal K-type is a pair $(G_{x,r_1}, \chi)$, where $\chi$ is a nondegenerate character of the abelian quotient $G_{x,r_1}/G_{x,r_1+}$. 
 While the work of Moy and Prasad revolutionized the study of representations of $p$-adic groups, the unrefined minimal K-type itself determines the representation only in some special cases. Our first main result in this paper (Theorem \ref{Thm-existence-of-datum}) shows that every smooth irreducible representation of $G(k)$ contains a much more refined invariant, which we call a \textit{\datum/}. A \datum/ is a tuple 
$$(x,(X_i)_{1 \leq i \leq n}, (\repzero, \Vrepzero))$$
 for some integer $n$,  where $x \in \sB(G,k)$, $X_i \in \fg^*$ for $1 \leq i \leq n$ satisfying certain conditions and $(\repzero, \Vrepzero)$ is an irreducible  representation of a finite group (which is the reductive quotient of the special fiber of the connected parahoric group scheme attached to the derived group of a twisted Levi subgroup of $G$), see Definition \ref{Def-extendeddatum} and Definition \ref{Def-datum} for the details.
  Our \datum/ can be viewed as a refinement of the unrefined minimal K-type of Moy and Prasad as follows. To a \datum/ $(x,(X_i)_{1 \leq i \leq n}, (\repzero, \Vrepzero))$ we associate a sequence of subgroups  $G \supset H_1 \supset H_2 \supset \hdots \supset H_n \supset H_{n+1}$, which are (apart from allowing $H_1=G_1$) the derived groups of twisted Levi subgroups $G = G_1 \supset G_2 \supset \hdots \supset G_n \supset G_{n+1}$, and real numbers $r_1>r_2> \hdots >r_n >0$ such that for $1 \leq i \leq n$ the element $X_i$  yields a character $\chi_i$ of  
  $$(H_i)_{x_i,r_i}/(H_i)_{x_i,r_i+} \simeq \Lie(H_i)(k)_{x_i,r_i}/\Lie(H_i)(k)_{x_i,r_i+} \subset \fg_{x,r_i}/\fg_{x,r_i+} $$
  for a suitable point $x_i \in \sB(H_i,k)$.
  If $G$ is semisimple,  for simplicity, then $(G_{x_1,r_1}, \chi_1)$ is an unrefined minimal K-type of depth $r_1$ contained in $(\pi,V_\pi)$, and the pair $((H_i)_{x_i,r_i}, \chi_i)$ is an unrefined minimal K-type of depth $r_i$ for $H_i$.

The existence of a maximal \datum/ for any irreducible representation of $G(k)$ is a key ingredient for producing the input that is needed for the construction of types as in Kim--Yu (\cite{KimYu}). 
In order to exhibit a \datum/ in a given representation, we require the elements $(X_i)_{1 \leq i \leq n}$ in the \datum/  $(x,(X_i)_{1 \leq i \leq n}, (\repzero, \Vrepzero))$ to satisfy a slightly stronger condition than the non-degeneracy necessary for an unrefined minimal K-type. We call our conditions \textit{\generic/}, see Definition \ref{Def-generic}. This condition ensures that the deduced input for Yu's construction is generic in the sense of Yu (\cite[\S~15]{Yu}) and at the same time it is crucial for the proof of the existence of the \datum/. 
 The existence of a \datum/ in $(\pi, V_\pi)$ is proved recursively, i.e. by first showing the existence of a suitable element $X_1$, then finding a compatible element $X_2$, and then $X_3$, etc., until we obtain a tuple $(X_i)_{1 \leq i \leq n}$ and finally exhibit the representation $(\repzero, \Vrepzero)$. The existence of $X_1$ can be considered as a refinement of the existence of an unrefined minimal K-type by Moy and Prasad and relies on the existence result of  \generic/ elements proved in Proposition \ref{Lemma-almoststable-generic-representative}. When constructing the remaining part of the \datum/ we need to ensure its compatibility with $X_1$ and rely on several preparatory results proved in Section \ref{Section-killing-lemmas}. At this step the imposed conditions on the elements $X_i$ become essential.

The final crucial part of this paper is concerned with deducing from the existence of a \datum/ in Theorem \ref{Thm-exhaustion-of-types} that every smooth irreducible representation of $G(k)$ contains one of the \type/s constructed by Kim--Yu and, similarly, that every irreducible supercuspidal representation of $G(k)$ arises from Yu's construction, see Theorem \ref{Thm-exhaustion-Yu}. This requires using the elements $X_i$ from the \datum/ to provide appropriate characters of the twisted Levi subgroups $G_i$ ($2 \leq  i \leq n+1$) and using $(\repzero, \Vrepzero)$ to produce a depth-zero supercuspidal representation $\pi_0$ of $G_{n+1}(k)$. We warn the reader that the depth-zero representation of $G_{n+1}(k)$ is in general not simply obtained by extending and inducing $(\repzero, \Vrepzero)$. The relationship between $\repzero$ and $\pi_0$ requires the study of Weil representations and can be found in Section \ref{Section-existence-of-type}, in particular in Lemma \ref{Lemma-repzeroK}. 
The main difficulty lies in showing that a potential candidate for the depth-zero representation $\pi_0$ of $G_{n+1}(k)$ is supercuspidal, which is the content of Lemma \ref{Lemma-cuspidal}.

 We conclude the paper by mentioning in Corollary \ref{Cor-supercuspidal-criterion} how to read off from a maximal \datum/ for $(\pi,V_\pi)$ if the representation $(\pi, V_\pi)$ is supercuspidal or not.

We would like to point out that the exhaustion of supercuspidal representations and the existence of types for arbitrary smooth irreducible representations have already been extensively studied for special classes of reductive groups for which other case-specific tools are available, e.g. a lattice theoretic description of the Bruhat--Tits building and a better understanding of the involved Hecke algebras. In 1979, Carayol (\cite{Carayol}) gave a construction of all supercuspidal representations of $\GL_n(k)$ for $n$ a prime number. In 1986, Moy (\cite{Moy-exhaustion}) proved that Howe's construction (\cite{Howe}) exhausts all supercuspidal representations of $\GL_n(k)$ if $n$ is coprime to $p$. Bushnell and Kutzko extended the construction to $\GL_n(k)$ for arbitrary $n$ and proved that every irreducible representation of $\GL_n(k)$ contains a type (\cites{BK, BK-types, BK-types-exhaustion}). As mentioned above, these results play a crucial role in the representation theory of $\GL_n(k)$. Based on the work for $\GL_n(k)$, Bushnell and Kutzko (\cite{BK-SLn}) together with Goldberg and Roche (\cite{Goldberg-Roche}) provide types for all Bernstein components for $\SL_n(k)$. For classical groups Stevens (\cite{Stevens}) has recently provided a construction of supercuspidal representations for $p \neq 2$ and proved that all supercuspidal representations arise in this way. A few years later, Miyauchi and Stevens  (\cite{Miyauchi-Stevens}) provided types for all Bernstein components in that setting. The case of inner forms of $\GL_n(k)$ was completed by Sécherre and Stevens (\cites{Secherre-Stevens, Secherre-Stevens-types}) around the same time, subsequent to earlier results of others for special cases (e.g. Zink (\cite{Zink}) treated division algebras over  non-archimedean local fields of characteristic zero and Broussous (\cite{Broussous}) treated division algebras without restriction on the characteristic). The existence of types for inner forms of $\GL_n(k)$ plays a key role in the explicit description of the local Jacquet--Langlands correspondence.

\textbf{Structure of the paper.} In Section \ref{Section-prime}, we collect some consequences of the assumption that the residual field characteristic $p$ does not divide the order of the Weyl group of $G$. Section \ref{Section-almost-stable} concerns the definition and properties of \generic/ elements and includes an existence result for  \generic/ elements. In Section \ref{Section-datum}, we introduce the notion of a \datum/ and define what it means for a representation to contain a \datum/ and for a \datum/ to be a maximal \datum/ for a representation. The proof that every smooth irreducible representation of $G(k)$ contains a \datum/ is the subject of Section \ref{Section-existence-of-datum}. Several results that are repeatedly used in this proof are shown in the preceding section, Section \ref{Section-killing-lemmas}. In Section \ref{Section-existence-of-type}, we use the result about the existence of a \datum/ to derive that every smooth irreducible representation of $G(k)$ contains one of the types constructed by Kim and Yu, and, in Section \ref{Section-exhaustion-Yu}, we prove analogously that every smooth irreducible supercuspidal representation of $G(k)$ arises from Yu's construction. 

\textbf{Conventions and notation.}
Throughout the paper, we require reductive groups to be connected and all representations are smooth complex representations unless mentioned otherwise. We do not distinguish between a representation and its isomorphism class. As explained in the introduction, by \textit{\type/} we mean an $\fs$-type for some inertial equivalence class $\fs$.

We will use the following notation throughout the paper: $k$ is a non-archimedean local field (of arbitrary characteristic) and $G$ is a reductive group over $k$ that will be assumed to split over a tamely ramified field extension of $k$. We write $\ff$ for the residue field of $k$ and denote its characteristic by $p$. We fix an algebraic closure $\ov k$ of $k$ and all field extensions of $k$ are meant to be algebraic and assumed to be contained in $\ov k$. For a field extension $F$ of $k$, we denote by $F^{{ur}}$ its maximal unramified field extension (in $\ov k$). 
We write $\cO$ for the ring of integers of $k$, $\cP$ for its maximal ideal, $\pi$ for a uniformizer, and $\val:k \ra \bZ \cup \{\infty\}$ for a valuation on $k$ with image $\bZ \cup \{\infty\}$. If $F$ is an (algebraic) field extension of $k$, then we also use $\val$ to denote the valuation on $F$ that extends the valuation on $k$. We write $\cO_F$ for the ring of integers in $F$ and $\cP_F$ for the maximal ideal of $\cO_F$.

Throughout the paper we fix an additive character $\varphi: k \ra \bC^*$ of $k$ of conductor $\cP$.

If $E$ is a field extension of a field $F$ (e.g. of $k$ or $\ff$) and $H$ is a scheme defined over the field $F$, then we denote by $H_E$ or $H \times_F E$ the base change $H \times_{\Spec F} \Spec E$. If $A$ is an $F$-module, then we write $A_E$ for $A \otimes_F E$ and $A^*$ for the $F$-linear dual of $A$. For $X \in A^*$ and $Y \in A_E$, we write $X(Y)$ for $(X \otimes 1)(Y)$, $(X \otimes 1) \in A^* \otimes_F E \simeq (A_E)^*$. If a group acts on $A$, then we let it also act on $A^*$ via the contragredient action.

In general, we use upper case roman letters, e.g. $G, H, G_i, T, \hdots$, to denote linear algebraic groups defined over a field $F$, and we denote the $F$-points of their Lie algebras by the corresponding lower case fractur  letters, e.g. $\fg, \fh, \fg_i, \ft$. The action of the group on its Lie algebra is the adjoint action, denoted by $\Ad$, unless specified otherwise. If $H$ is a reductive group over $F$, then we denote by $H^{\text{der}}$ its derived group. 
We write $\bG_a$ and $\bG_m$ for the additive and multiplicative group schemes over $\bZ$ or over the ring or field that becomes apparent from the context.  If $S$ is a split torus contained in $H$ (defined over $F$), then we write $X^*(S)=\Hom_F(S,\bG_m)$ for the characters of $S$ defined over $F$, $X_*(S)=\Hom_F(\bG_m,S)$ for the cocharacters of $S$ (defined over $F$), $\Phi(H,S) \subset X^*(S)$ for the roots of $H$ with respect to $S$, and if $S$ is a maximal torus, then $\check \Phi(H,S) \subset X_*(S)$ denotes the coroots. We might abbreviate $\Phi(H_{\ov k},T)$ by $\Phi(H)$ for a maximal torus $T$ of $H_{\ov k}$ if the choice of torus $T$ does not matter. We use the notation $\< \cdot, \cdot\>: X_*(S) \times X^*(S) \ra \bZ$ for the standard pairing, and if $S$ is a maximal torus, then we denote by $\check\alpha \in \check \Phi(H,S)$ the dual root of $\alpha \in \Phi(H,S)$.
For a subset $\Phi$ of $X^*(S) \otimes_\bZ \bR$ (or $X_*(S)\otimes_\bZ \bR$) and $R$ a subring of $\bR$, we denote by $R\Phi$ the smallest $R$-submodule of $X^*(S) \otimes_\bZ \bR$ (or $X_*(S)\otimes_\bZ \bR$, respectively) that contains $\Phi$. For $\chi \in X^*(S)$ and $\lambda \in X_*(S)$, we denote by $d\chi \in \Hom_F(\Lie(S),\Lie(\bG_m))$ and $d\lambda \in \Hom_F(\Lie(\bG_m), \Lie(S))$ the induced morphisms of Lie algebras.

If $(\pi, V)$ is a representation of a group $Q$, then we denote by $V^Q$ the elements of $V$ that are fixed by $Q$. If $Q'$ is a group containing $Q$ as a subgroup and $q' \in Q'$, then we define the representation $({^{q'}} \pi, V)$ of $q'Q{q'}^{-1}$ by ${^{q'}} \pi(q)=\pi({q'}^{-1}qq')$ for all $q \in q'Q{q'}^{-1}$.

Finally, we let $\wt \bR=\bR \cup \{ r+ \, | \, r \in \bR\}$ with its usual order, i.e. for $r$ and $s \in \bR$ with $r<s$, we have $r<r+<s<s+$.

\textbf{Acknowledgment.}
The author thanks Stephen DeBacker, Wee Teck Gan, Tasho Kaletha, Ju-Lee Kim and Loren Spice for discussions related to this paper, as well as Jeffrey Adler, Anne-Marie Aubert, Stephen DeBacker, Tasho Kaletha, Ju-Lee Kim, Gopal Prasad, Vincent Sécherre and Maarten Solleveld for feedback on some parts of an earlier version of this paper.
The author is also very grateful to the referee for a careful reading of the paper and helpful comments and suggestions.
The author thanks the University of Michigan, the Max-Planck Institut für Mathematik and the Institute for Advanced Study for their hospitality and wonderful research environment.

\section{Assumption on the residue field characteristic} \label{Section-prime}
Recall that $k$ denotes a non-archimedean local field with residual characteristic $p$ and $G$ is a connected reductive group over $k$. We assume that $G$ is not a torus. We know already all the smooth, irreducible, supercuspidal representations of a torus. They are simply the smooth characters of the torus. Moreover, Yu (\cite{Yu}) works in his construction of supercuspidal representations only with tori of $G$ that split over a tame extension. Hence we make the following assumption throughout the paper.

\begin{Assumption} \label{Assumption-p-W}
	We assume that $G$ splits over a tamely ramified extension of $k$ and $p \nmid \abs{W}$, where $W$ denotes the Weyl group $W$ of $G(\ov k)$.
\end{Assumption}

By \cite{Fi-tame-tori} Assumption \ref{Assumption-p-W} implies that all tori of $G$ are tame. For absolutely simple groups other than some non-split inner forms of split groups of type $A_n, n \geq 2, D_l, l \geq 4$ prime, or $E_6$ this assumption is also necessary (and in the excluded cases only minor modifications on the assumption on $p$ are necessary), see \cite[Theorem~2.4 and Corollary~2.6]{Fi-tame-tori} for details.

We collect a few consequences of our assumption for later use.

\begin{Lemma} \label{Lemma-restriction-on-p}
The assumption that $p \nmid \abs{W}$	implies the following
\begin{enumerate} [label=(\alph*),ref=\alph*]
	\item \label{item-p-Levi} The prime $p$ does not divide the order of the Weyl group of any Levi subgroup of (a parabolic subgroup of) $G_{\ov k}$.
	\item \label{item-bond} The prime $p$ is larger than the order of any bond of the Dynkin diagram $\Dyn(G)$ of $G_{\ov k}$, i.e. larger than the square of the ratio of two root lengths of roots in $\Phi(G)$.
 	\item  \label{item-bad-prime}  The prime $p$ is not a bad prime (in the sense of \cite[4.1]{Springer-Steinberg}) for $\check\Phi:=\check\Phi(G)$, i.e. $\bZ\check\Phi/\bZ\check\Phi_0$ has no $p$-torsion for all closed subsystems $\check\Phi_0$ in $\check \Phi$
 	\item  \label{item-torsion-prime}  The prime $p$ is not a torsion prime (in the sense of \cite[1.3~Definition]{Steinberg-torsion}) for $\Phi:=\Phi(G)$ (and hence also not for $\check\Phi(G)$), i.e. $\bZ\check\Phi/\bZ\check\Phi_0$ has no $p$-torsion for all closed subsystems $\Phi_0$ in $\Phi$.
	\item \label{item-index-of-connection} The prime $p$ does not divide the index of connection (i.e. the order of the root lattice in the weight lattice) of any root(sub)system generated by a subset of a basis of $\Phi(G)$. 
	\end{enumerate}
\end{Lemma}
\Proof
Part \eqref{item-p-Levi} is obvious, Part \eqref{item-bond}, \eqref{item-bad-prime} and \eqref{item-torsion-prime} can be read of from Table \ref{table-weyl-group}. Part \eqref{item-index-of-connection} follows from the fact that the index of connection of $\Phi(G)$ divides $\abs{W}$ (\cite[VI.2~Proposition~7]{Bourbaki-4-6}).\qed

\begin{table}[h]\footnotesize
	\begin{tabular}{|c|c|c|c|c|c|c|c|c|c|}
		\hline type  & 			$A_n \, (n \geq 1)$     & $B_n \, (n \geq 3)$  & $C_n \, (n \geq 2)$  &  $D_n \, (n \geq 3)$              & $E_6$                 &  $E_7$                      &  $E_8$ & $F_4 $ & $G_2$  \\ 
		\hline $\abs{W}$ &  $(n+1)!$  & $2^n \cdot n!$ &$2^n \cdot n!$   &  $2^{n-1} \cdot n!$ & $2^7\cdot3^4\cdot5$  &  $2^{10}\cdot3^4\cdot5\cdot7$ & $2^{14}\cdot3^5\cdot5^2 \cdot 7$ &  $2^7 \cdot 3^2$ & $2^2 \cdot 3$ \\
		\hline bad  &  -  & 2 & 2   &  2 & 2, 3 &  2, 3 & 2, 3, 5 &  2, 3 &  2, 3 \\ 
		\hline torsion &  -  & 2 & -  &  2 & 2, 3 &  2, 3 & 2, 3, 5 &  2, 3 & 2 \\  
		\hline 
	\end{tabular} 
	\caption{Order of Weyl groups (\cite[VI.4.5-VI.4.13]{Bourbaki-4-6}); bad primes (\cite[4.3]{Springer-Steinberg}) and torsion primes (\cite[1.13~Corollary]{Steinberg-torsion}) for irreducible root systems}
	\label{table-weyl-group}
\end{table}

\section{{\Almoststronglystable/} and {\generic/} elements} \label{Section-almost-stable}

Let $E$ be a field extension of $k$. We denote by $\sB(G,E)$ the (enlarged) Bruhat--Tits building of $G_E$ over $E$, and we sometimes write $\sB$ for $\sB(G,k)$. For $x \in \sB(G,E)$ and $r \in \bR_{\geq 0}$, we write $G(E)_{x,r}$ for the Moy--Prasad filtration subgroup of $G(E)$ of depth $r$, which we abbreviate to $G_{x,r}$ for $G(k)_{x,r}$, and we set $G(E)_r=\bigcup_{x \in \sB(G,E)} G(E)_{x,r}$. For $r \in \bR$, we denote by $(\fg_E)_{x,r}$ and $(\fg_E)^*_{x,r}$ the Moy--Prasad filtration of $\fg_E=\Lie(G_E)(E)$ and its dual $\fg_E^*$, respectively. We set $(\fg_E)_r=\bigcup_{x \in \sB(G,E)}(\fg_E)_{x,r}$ and  $(\fg_E)^*_r=\bigcup_{x \in \sB(G,E)}(\fg_E)^*_{x,r}$. Recall that if $X \in (\fg_E)^*_{x,r}$, then $X((\fg_E)_{x,(-r)+}) \subset \cP_E$ and  $X((\fg_E)_{x,-r}) \subset \cO_E$.   
For convenience, we define our Moy--Prasad filtration subgroups and subalgebras with respect to the valuation $\val$ of $E$ that extends the normalized valuation $\val$ of $k$, i.e. in such a way that by \cite[Proposition~1.4.1]{Adler} (which applies because $G$ splits over a tamely ramified extension) we have
\begin{equation} \label{eqn-MP-filtration}
	(\fg_E)_{x,r} \cap \fg = \fg_{x,r}
\end{equation}
for all $r \in \bR$.
We denote by $\RP_x$ the reductive quotient of the special fiber of the connected parahoric group scheme attached to $G$ at $x$, i.e. $\RP_x$ is a reductive group defined over $\ff$ satisfying $\RP_x(\ff_F)=G(F)_{x,0}/G(F)_{x,0+}$ for every unramified extension $F$ of $k$ with residue field $\ff_F$. We also refer to $\RP_x$ as ``the reductive quotient of $G$ at $x$''. For any $r \in \bR$, the adjoint action of $G(k^{ur})_{x,0}$ on $(\fg_{k^{ur}})_{x,r}$ induces a linear action of the algebraic group $\RP_x$ on $V_{x,r}:=\fg_{x,r}/\fg_{x,r+}$ and therefore also on its dual $V_{x,r}^*$.

For $X \in \fg_E^*-\{0\}$ and $x \in \sB(G,E)$, we denote by $d_E(x,X) \in \bR$ the largest real number $d$ such that $X \in (\fg_E^*)_{x,d}$, and set $d_E(x,0)=\infty$. We call $d_E(x,X)$ the \textit{depth of $X$ at $x$}. We define the \textit{depth of $X$} over $E$ to be $d_E(X):=\sup_{x \in \sB(G,E)}d_E(x,X) \in \bR \cup \infty$. If $E=k$, then we often write $d(x,X)$ for $d_k(x,X)$ and $d(X)$ for $d_k(X)$. Note that if $X \in \fg^*$, then $d_k(x,X)=d_E(x,X)$ and $d_k(X)=d_E(X)$ by our choice of normalization.

Recall that if $V$ is a finite dimensional linear algebraic representation of a reductive group $H$ defined over some field $F$, then $X \in V(F)$ is called \textit{semistable} under the action of $H$ if the Zariski-closure of the orbit $H(\ov F).X\subset V(\ov F)$ does not contain zero, and is called \textit{unstable} otherwise. We introduce two slightly stronger notions for our setting. 

\begin{Def} \label{Def-almost-stable}
	Let $X \in \fg^*$. We denote by $\ov X$ the map $V_{x,-d(x,X)}:=\fg_{x,-d(x,X)}/\fg_{x,(-d(x,X))+} \ra \ff$ induced from $X:\fg_{x,-d(x,X)} \ra \cO$.  
	\begin{itemize}
		\item We say that $X$ is \textit{\almoststable/} if the $G$-orbit of $X$ is closed. 
		\item We say that $X$ is \textit{\almoststronglystable/ at $x$} if $X$ is \almoststable/ and $\ov X \in (V_{x,-d(x,X)})^*$ is semistable under the action of $\RP_x$. 		
	\end{itemize}
	
\end{Def}

\begin{Lemma} \label{Lemma-almost-stable-depth}
	Let $X \in \fg^*-\{0\}$ be \almoststronglystable/ at $x$. Then $d(x,X)=d(X)$.
\end{Lemma}
\Proof Suppose $d(x,X)<d(X)$, and write $r=d(x,X)$. Then by \cite[Corollary~3.2.6]{Adler-DeBacker} (together with their remark at the beginning of Section~3) the coset $X+\fg^*_{x,r+}$ is degenerate, i.e. contains an unstable element. Hence $\ov X$ is unstable by \cite[4.3.~Proposition]{MP1} (while Moy and Prasad assume simply connectedness throughout their paper \cite{MP1}, it is not necessary for this claim). This contradicts that $X$ is \almoststronglystable/ and finishes the proof. \qed

\begin{Def}
	Let $H$ be a reductive group over some field $F$. A smooth, closed subgroup $H'$ of $H$ is called a \textit{twisted Levi subgroup} if there exists a finite field extension $E$ over $F$ such that $H'\times_{F} E$ is a Levi subgroup of a parabolic subgroup of $H \times_{F} E$.
\end{Def}

\begin{Lemma} \label{Lemma-stabilizer-Levi}
	Let $X \in \fg^*$ be \almoststable/ (under the contragredient of the adjoint action of $G$). Then the centralizer $\Cent_{G}(X)$ of $X$ in $G$ is a twisted Levi subgroup of $G$.
\end{Lemma}
\Proof
It suffices to show that $\Cent_{G\times_k{\ov k}}(X)$ is a Levi subgroup of $G_{\ov k}$, because $\Cent_{G}(X)\times_{k}\ov k=\Cent_{G\times_k{\ov k}}(X)$.
Since $p$ does not divide the order of the Weyl group of $G$, we can $G_{\ov k}$-equivariantly identify $\fg^*_{\ov k}$  with $\fg_{\ov k}$ (see \cite[Proposition~4.1]{Adler-Roche}, Lemma \ref{Lemma-restriction-on-p}\eqref{item-index-of-connection} and Lemma \ref{Lemma-restriction-on-p}\eqref{item-bond}). Using this identification to view $X$ in $\fg$, by \cite[14.25~Proposition]{Borel}, every $X$ is contained in the Lie algebra of a Borel subgroup $B=TU$ for $T$ a maximal torus and $U$ the unipotent radical of $B$ (defined over $\ov k$). Hence we can write $X=X_s+X_n$, where $X_s \in \Lie(T)(\ov k)$ and $X_n \in \Lie(U)(\ov k)$, and there exists a one parameter subgroup $\lambda: \bG_m \ra T \subset \Cent_{G_{\ov k}}(X_s)$ such that $\lim_{t \ra 0} \lambda(t).X_n=0$, and therefore $\lim_{t \ra 0} \lambda(t).X=X_s$. Since $X$ is \almoststable/, this implies that $X_s$ is contained in the $G(\ov  k)$-orbit of $X$, and therefore $X$ is semisimple, hence $X=X_s$. 
In other words, 
 $X$ is in the zero eigenspace in $\fg^*_{\ov k}$ of $T$.
Since $p$ is not a torsion prime for $\check \Phi(G)$ (Lemma \ref{Lemma-restriction-on-p}\eqref{item-torsion-prime}) and $p$ does not divide the index of connection of $\Phi(G)$ (Lemma \ref{Lemma-restriction-on-p}\eqref{item-index-of-connection}), we obtain by \cite[Proposition~7.1.~and~7.2]{Yu} (which is based on \cite{Steinberg-torsion}) that the centralizer $\Cent_{G_{\ov k}}(X)$ of $X$ in $G_{\ov k}$ is a connected reductive group whose root datum is given by $(X^*(T), \Phi_X, X_*(T), \check\Phi_X)$ with $\Phi_X=\{\alpha \in \Phi(G_{\ov k},T) \, | \, X(d\check\alpha(1))=0)\}$ and $\check \Phi_X=\{\check \alpha \, | \, \alpha \in \Phi_X\}$. Note that $\check\Phi_X$ is a closed subsystem of $\check\Phi$ (i.e. $\bZ \check\Phi_X \cap \check\Phi = \check\Phi_X$). Since $\bZ\check\Phi/\bZ\check\Phi_X$ is $p$-torsion free by Lemma \ref{Lemma-restriction-on-p}\eqref{item-bad-prime}, we have $\check\Phi_X = \bQ\check\Phi_X \cap \check\Phi$ and hence $\Phi_X = \bQ\Phi_X \cap \Phi$. By \cite[VI.1,~Proposition~24]{Bourbaki-4-6} there exists a basis $\Delta$ for $\Phi$ containing a basis $\Delta_X$ for $\Phi_X$. Thus $\Cent_{G_{\ov k}}(X)$ is a Levi subgroup of $G_{\ov k}$. \qed

\begin{Def} \label{Def-generic}
	We say that an element $X \in \fg^*$ is	\textit{\generic/ of depth $r$ at $x$} $\in \sB(G,k)$  if $X$ is \almoststable/ and if there exists a tamely ramified extension $E$ over $k$ and a split maximal torus $T \subset \Cent_G(X)\times_k E$ such that
	\begin{itemize}
		\item $x \in \sA(T,E) \cap \sB(G,k)$, where $\sA(T,E)$ denotes the apartment of $T$ in $\sB(G,E)$
		\item $X \in \fg^*_{x,r}$ (i.e. $X(\fg_{x,(-r)+})\subset \cP$), 
		\item for every $\alpha \in \Phi(G,T)$ we have $X(H_\alpha)=0$ or $\val(X(H_\alpha))=r$, where $H_\alpha=d\check\alpha(1)$, and
		\item if $X(H_\alpha)=0$ for all $\alpha \in \Phi(G,T)$, then $d(x,X)=r$.
	\end{itemize}
\end{Def}

Note that $H_\alpha=d\check\alpha(1) \neq 0$, because $p$ does not divide the index of connection of $\Phi(G)$ by Lemma \ref{Lemma-restriction-on-p}\eqref{item-index-of-connection}. We will see in Corollary \ref{Cor-generic-depth} below that if $X$ is \generic/ of depth $r$ at $x$, then $d(x,X)=r$.

\begin{Lemma} \label{Lemma-structure-of-X}
Let $X \in \fg^*$ be \generic/ of depth $r$ at $x$. Then for every (split) maximal torus $T \subset \Cent_G(X) \times_k \ov k$ we have
\begin{itemize}
 \item  $X(H_\alpha)=0$ for all $\alpha \in \Phi(\Cent_G(X),T)$, and
 \item  $\val(X(H_\alpha))=r$ for all $\alpha \in \Phi(G,T)-\Phi(\Cent_G(X),T)$. 
 \end{itemize}
Moreover, for all $\alpha \in \Phi(G,T)$ we have $X((\fg_{\ov k})_\alpha)=0$, where $(\fg_{\ov k})_\alpha$ denotes the $\alpha$-root subspace of $\fg_{\ov k}$.
\end{Lemma}
\Proof 
  Choose a Chevalley system $\{x_\alpha:\bG_a \ra G_{\ov k} \, | \, \alpha \in \Phi(G,T) \}$ with corresponding Lie algebra elements $\{X_{\alpha}=dx_\alpha(1) \, | \alpha \in \Phi(G,T)\}$. Since $T \subset \Cent_G(X)\times_k \ov k$, we have $X(X_\alpha)=X(\Ad(t)(X_\alpha))=\alpha(t)X(X_\alpha)$ for all $t \in T(\ov k)$, and hence $X(X_\alpha)=0$ for all $\alpha \in \Phi(G,T)$. Thus $X((\fg_{\ov k})_\alpha)=0$.

Since the split maximal tori  of $\Cent_G(X)\times_k \ov k$ are conjugate in $\Cent_G(X)\times_k \ov k$, we have $X(H_\alpha)=0$ or $\val(X(H_\alpha))=r$ for $\alpha \in \Phi(G,T)$.  By \cite[Propostion~7.1]{Yu}, we have $\alpha \in \Phi(\Cent_G(X),T)$ if and only if $X(H_\alpha)=0$, see also the proof of Lemma \ref{Lemma-stabilizer-Levi}.
 \qed

\begin{Cor} \label{Cor-generic-depth}
	Let $X \in \fg^*-\{0\}$ be \generic/ of depth $r$ at $x$. Then $d(x,X)=d(X)=r$.
\end{Cor}
\Proof Let $E$ be a tame extension of $k$ and $T$ a  split maximal torus of $\Cent_G(X) \times_k E$ such that $x\in \sA(T,E)$. By Lemma \ref{Lemma-restriction-on-p}\eqref{item-index-of-connection} the element $H_\alpha$ is of depth zero for all $\alpha \in \Phi(G,T)$. Hence $d(x,X)=d_E(x,X)=r$ by Lemma \ref{Lemma-structure-of-X} (or by definition if $X(H_\alpha)=0$ for all $\alpha \in \Phi(G,T)$). If $y \in \sB(G_E,E), s \in \bR$ and $X \in (\fg_E^*)_{y,s}$, then \cite[Lemma~8.2]{Yu} implies that $X$ restricted to $\Lie(T)(E)$ lies in $\Lie(T)^*(E)_{s}$. Since $X$ has depth $d(x,X)=r$ when restricted to $\Lie(T)(E)$, we deduce that $d_E(y,X)\leq r$. Hence $d(X)=d(x,X)=r$. \qed

\begin{Cor}\label{Cor-almoststable-and-generic-implies-almoststronglystable}
	Let $X \in \fg^*-\{0\}$ be \generic/ of depth $r$ at $x$. Then $X$ is \almoststronglystable/ at $x$.
\end{Cor}
\Proof Suppose $X$ is not \almoststronglystable/ at $x$. Then $\ov X \in \fg^*_{x, r}/\fg^*_{x, r+}$ is unstable.  Since $\ff$ is perfect, by \cite[Corollary~4.3]{Kempf} there exists a non-trivial one parameter subgroup $\ov \lambda: \bG_m \ra \RP_x$ in the reductive quotient $\RP_x$ of $G$ at $x$ (defined over $\ff$) such that $\lim_{t\ra0}\ov \lambda(t).\ov{X}=0$. Let $\bS$ be a maximal split torus of $\RP_x$ containing $\ov \lambda(\bG_m)$. 
Then there exists a split torus $\cS$ (defined over $\cO_k$) in the parahoric group scheme $\bP_x$ of $G$ at $x$ whose special fiber is $\bS$ and whose generic fiber $S$ is a split torus in $G$. This allows us to lift $\ov \lambda$ to a one parameter subgroup $\lambda: \bG_m \ra S \subset G$. Let $\sA(S, k)$ be the apartment of $S$ (i.e. the apartment of a maximal torus in $G$ that contains $S$). Then $\sA(S,k)$ contains $x$ and is the affine space underlying the real vector space $X_*(S) \otimes_{\bZ} \bR$. If $\eps>0$ is sufficiently small, we obtain $X \in \fg^*_{x+\eps\lambda, r+}$. Hence $d(X)>r$, which contradicts Corollary \ref{Cor-generic-depth}. 
\qed

\begin{Rem} \label{Rem-BT} 
Recall that if $G'$ is a Levi subgroup of (a parabolic subgroup of) $G$, then we have an embedding of the corresponding Bruhat--Tits buildings $\sB(G',k) \hookrightarrow \sB(G,k)$. Even though this embedding is only unique up to some translation, its image is unique. Since we assume that all tori of $G$ split over tamely ramified extensions of $k$, every twisted Levi subgroup of $G$ becomes a Levi subgroup over a finite tamely ramified extension of $k$. Hence using (tame) Galois descent, we obtain a well defined image of $\sB(G',k)$ in $\sB(G,k)=\sB$ for every twisted Levi subgroup $G'$ of $G$. In the sequel, we might identify $\sB(G',k)$ with its image in $\sB$.
\end{Rem}

\begin{Rem} \label{Rem-B}
	Since $p$ does not divide the index of connection of $\Phi(G)$ (Lemma \ref{Lemma-restriction-on-p}\eqref{item-index-of-connection}), Adler and Roche (\cite[Proposition~4.1]{Adler-Roche}) provide a non-degenerate, $G$-equivariant, symmetric bilinear form $B: \fg \times \fg \ra k$ such that the induced identification of $\fg$ with $\fg^*$ identifies $\fg_{x,r}$ with $\fg^*_{x,r}$ for all $x \in \sB(G,k), r \in \bR$. Moreover, $B$ stays non-degenerate when restricted to the Lie algebra of any twisted Levi subgroup of $G$.
\end{Rem} 
Using the bilinear form from Remark \ref{Rem-B} to view $(\fg')^*=(\Lie(G')(k))^*$ as a subset of $\fg^*$ for $G'$ a twisted Levi subgroup of $G$, we have the following lemma, which is a translation of a result by Kim--Murnaghan (\cite[Lemma~2.3.3]{Kim-Murnaghan}) into the dual setting.

\begin{Lemma}[Kim--Murnaghan] \label{Lemma-Kim-Murnaghan}
	Let $r\in \bR$, $x\in \sB$, and let $X\in \fg^*_{x,r}\subset \fg^*$ be \generic/ of depth $r$ at $x$. Denote $\Cent_G(X)$ by $G'$ and $\Lie(\Cent_G(X))(k)$ by $\fg'$. If $X' \in (\fg')^*_{r+} \subset \fg^*$ and $y \in \sB(G,k) - \sB(G',k)$, then $d(y,X+X')<d(X)$.
\end{Lemma}
\Proof
Suppose $X \neq 0$ as the statement is trivial otherwise. Let $E$ be a tame extension of $k$ and $T$ a split maximal torus of $\Cent_G(X) \times_k E$ such that $x\in \sA(T,E)$. By the definition of the bilinear form $B$ in the proof of \cite[Proposition~4.1]{Adler-Roche} (a sum of scalings of killing forms together with a bilinear form on the center) together with Lemma \ref{Lemma-structure-of-X}, the \generic/ element $X$ corresponds to an element $\check X$ of $\ft=\Lie(T)(k) \subset \fg$, hence of $\ft_r=\ft \cap \fg_{x,r}$. Moreover, it follows from the definition of the bilinear form that $d\alpha(\check X)=X(H_\alpha)$. Hence Lemma \ref{Lemma-structure-of-X} implies that $\check X$ is a good semisimple element of depth $r$ (see \cite[Definition~2.2.4]{Adler} for the definition of ``good semisimple element''). Since $\Cent_G(X)=\Cent_G(\check X)$ and $X'$ corresponds to an element in $\fg'_{r+}$ under the identification of $\fg^*$ with $\fg$, the lemma follows from \cite[Lemma~2.3.3]{Kim-Murnaghan}, because $B$ preserves depth. (Note that Kim and Murnaghan impose in \cite{Kim-Murnaghan} much stronger conditions on $G$ and $k$, in particular that $k$ has characteristic zero. However the required Lemma~2.3.3 holds also in our setting by the same proof and observing that Corollary~2.2.3 and Lemma~2.2.4 in \cite{Kim-Murnaghan} (which are used in the proof of \cite[Lemma~2.3.3]{Kim-Murnaghan}) follow from results of Adler and Roche (\cite{Adler-Roche}) that are valid in our situation by our Lemma \ref{Lemma-restriction-on-p}.) \qed

To state the following main result in this section more conveniently, we fix a $G$-equivariant distance function $\mathrm{d}:\sB(G,k) \times \sB(G,k) \ra \bR_{\geq 0}$ on the building $\sB(G,k)$, which is the restriction of a distance function $\mathrm{d}_E:\sB(G,E) \times \sB(G,E) \ra \bR_{\geq 0}$ for some tame extension $E$ of $k$ over which $G$ splits that satisfies $\abs{\alpha(x-y)}\leq \mathrm{d}_E(x,y)$ for all maximal split tori $T_E$ of $G_E$, all $x, y \in \sA(T_E,E)$ and all $\alpha \in \Phi(G_E,T_E)$. (This normalization will only become relevant in the proof of Theorem \ref{Thm-existence-of-datum} below.)

\begin{Prop} \label{Lemma-almoststable-generic-representative}
	Let $r \in \bR$ and $x \in \sB$. If $X \in \fg^*$ is \almoststronglystable/ at $x$ with $d(x,X)=r$, then for every $\eps > 0$ there exists  $x' \in \sB$ with $\mathrm{d}(x,x')<\eps$ such that $X \in \fg^*_{x',r}$, the coset $X+\fg^*_{x',r+}$ contains an element $\wt X$ that is \generic/ of depth $r$ at $x'$, and the points $x$ and $x'$ are contained in $\sB(\Cent_G(\wt X),k) \subset \sB$.
\end{Prop}
\Proof 
Let $T$ be a  maximal torus of $\Cent_G(X)$ and $E$ a tame extension of $k$ over which $T$ splits. Choose a point $y$ in $\sA(T,E) \cap \sB(G,k)$. If $\alpha \in \Phi(G,T_E)$ and $X_\alpha \in (\fg_E)_\alpha$, then $X(X_\alpha)=X(\Ad(t)X_\alpha)=\alpha(t)X(X_\alpha)$ for all $t \in T(E)$, hence $X(X_\alpha)=0$. Thus the depth of $X$ at $y$ is equal to the depth of $X$ restricted to $\ft=\Lie(T)(k)$. On the other hand, by \cite[Lemma~8.2]{Yu}, the assumption that $X \in \fg^*_{x,r}$ implies that $X$ restricted to $\ft$ lies in $\ft^*_r$. Hence $d(y,X) \geq r$. Since $r=d(X)$ by Lemma \ref{Lemma-almost-stable-depth}, we deduce that $d(y,X)=r$.

\textbf{Claim.} $X+\fg^*_{y,r+}$ contains a \generic/ element of depth $r$ at $y$.

\textbf{Proof of claim.}
 Let $\check\Phi_0 \subset \check\Phi:=\check\Phi(G,T_{E})$ be the collection of coroots $\check \alpha$ for which $\val(X(H_\alpha))>r$. Note that $\check\Phi_0$ is a closed subsystem of $\check\Phi$ (i.e. $\bZ \check\Phi_0 \cap \check\Phi = \check\Phi_0$). Since $\bZ\check\Phi/\bZ\check\Phi_0$ is $p$-torsion free by Lemma \ref{Lemma-restriction-on-p}\eqref{item-bad-prime}, we also have $\check\Phi_0 = \bQ\check\Phi_0 \cap \check\Phi$. Moreover, since $X$ and $T$ are defined over $k$, the set $\check\Phi_0$ is stable under the action of the  Galois group $\Gal(E/k)$. Let $Y \subset \fg_{E}=\Lie(G)(E)$ be the $E$-subspace spanned by $\{H_\alpha \, | \, \check\alpha \in \check\Phi_0\}$. By the above observations about $\check \Phi_0$, the subspace $Y$ is $\Gal(E/k)$-stable, and if $H_\alpha \in Y$, then $\check\alpha \in \check\Phi_0$. Define 
$$Y_T^\perp:=\left\{ Z \in \Lie(T)({E}) \, | \, d\alpha(Z)=0 \, \forall \check\alpha \in \check\Phi_0 \right\}.$$
Then $Y_T^\perp$ is a $\Gal(E /k)$-stable complement to $Y$ in $\Lie(T)({E})$, and we set
$$Y^\perp := Y_T^\perp \oplus \bigoplus_{\alpha \in \Phi(G,T)} (\fg_{E})_\alpha .$$
Then $Y^\perp$ is a $\Gal(E /k)$-stable complement to $Y$ in $\fg_{E}$, and we define $X' \in \fg_{E}^*$ by  
	$$ X'(Z+Z^\perp)=X(Z) \, \text{ for all } Z \in Y, Z^\perp \in Y^\perp . $$
	Since $Y$ and $Y^\perp$  are $\Gal(E/k)$-stable and $X$ is defined over $k$, the linear functional $X'$ is $\Gal(E/k)$-invariant and hence defined over $k$, i.e. we can view $X'$ as an element of $\fg^*$. 

Let $\check \Delta_0$ be a basis for $\check \Phi_0$, and $\check \Delta$ a basis for $\check \Phi$ containing $\check \Delta_0$ (such a $\check \Delta$ exists by \cite[VI.1,~Proposition~24]{Bourbaki-4-6}). For $\check\alpha \in \check\Delta$, we denote by $\check\omega_\alpha \in \bQ\check\Phi$ the fundamental coweight corresponding to $\alpha$, i.e. $\<\check\omega_\alpha,\alpha\>=1$ and $\<\check\omega_\alpha,\beta\>=0$ for $\check \beta \in \check\Delta-\{\check\alpha\}$. Similarly, for $\check\alpha \in \check\Delta_0$, let $\check\omega^0_\alpha \in \bQ\check\Phi_0$ be the fundamental coweight with respect to the (co-)root system $\check \Phi_0$. By Lemma \ref{Lemma-restriction-on-p}\eqref{item-index-of-connection}, we have $\check \omega_\alpha \in \bZ\left[\frac{1}{\abs{W}}\right]\check\Phi$ and $\check \omega^0_\alpha \in \bZ\left[\frac{1}{\abs{W}}\right]\check\Phi_0$.
Denote by $H_{\check\omega_\alpha}$ ($\check\alpha \in \check\Delta$) and $H_{\check\omega^0_{\alpha'}}$  ($\check{\alpha'} \in \check\Delta_0$)  the image of $\check\omega_\alpha$ and  $\check\omega^0_{\alpha'}$ under the linear map $\bZ\left[\frac{1}{\abs{W}}\right]\check\Phi \rightarrow \Lie(T)(E) $ obtained by sending $\check{\alpha''}$ to $H_{\alpha''}$ (${\alpha''} \in \Phi$). 
Then we have
$$H_{\check\omega_\alpha}\equiv \left\{ \begin{array}{rl} 0  \mod Y^\perp & \text{ for } \check\alpha \in \check\Delta-\check\Delta_0 \\ H_{\check\omega^0_\alpha} \mod Y^\perp & \text{ for } \check \alpha \in \check \Delta_0  . \end{array} \right. $$
For $\beta \in \Phi$, we have $\check \beta=\sum_{\check \alpha \in \check \Delta}\<\check\beta, \alpha \> {\check\omega_\alpha}$, and hence we obtain
\begin{equation}\label{equation-X-beta} 
	H_\beta= \sum_{\check \alpha \in \check \Delta}\<\check\beta, \alpha \> H_{\check\omega_\alpha} \equiv \sum_{\check \alpha \in \check \Delta_0}\<\check\beta, \alpha \> H_{\check\omega^0_\alpha} \mod Y^\perp .
\end{equation}
Recall that $\<\check\beta, \alpha \>$ are integers for $\check\alpha \in \check \Delta$ and that the index of the coroot lattice $\bZ\check\Phi_0$ in the coweight lattice  $\bZ[\check\omega_{\alpha}\, | \, \check\alpha \in \check\Phi_0]$ is coprime to $p$ by Lemma \ref{Lemma-restriction-on-p}\eqref{item-index-of-connection}. Hence $\sum_{\check \alpha \in \check \Delta_0}\<\check\beta, \alpha \> H_{\check\omega^0_\alpha}$ is contained in the $\cO_{\ov k}$-span of $\{ H_\alpha \, | \, \check\alpha \in \check \Phi_0\}$. Thus, by the definition of $\check\Phi_0$, we obtain $\val(X'(H_\beta))>r$ for all $\check\beta \in \check\Phi$. In addition, $X'$ vanishes on the center of $\fg_E$ and on $\bigoplus_{\alpha \in \Phi(G,T)} (\fg_{E})_\alpha$, because these subspaces are contained in $Y^\perp$. Hence, by Lemma \ref{Lemma-restriction-on-p}\eqref{item-index-of-connection}, we have $\val(X'((\fg_{E})_{y,0})) \subset \bR_{>r}$.

 Using that the Moy--Prasad filtration behaves well with respect to base change (Equation \eqref{eqn-MP-filtration}), we obtain 
$$\val(X'(\fg_{y,{-r}})) \subset \val(X'((\fg_{E})_{y,-r}))  \subset \bR_{\geq -r} + \val(X'((\fg_{E})_{y,0})) \subset \bR_{>0} .$$
Thus $X'\in \fg^*_{y,r+}$, and $\wt X=X-X' \in X+\fg^*_{y,r+}$ with $\val(\wt X(H_\alpha))=r$ for $\check \alpha \not \in \check\Phi_0$ and $\wt X(H_\alpha)=0 $ for $\check\alpha \in \check\Phi_0$. 

In order to prove the claim, it remains to show that the orbit of $\wt X$ is closed. Since $p \nmid \abs{W}$ we can $G$-equivariantly identify $\fg^*$ with $\fg$ as in Remark \ref{Rem-B}. Since $T$ is in $\Cent_G(X)$ and acts trivially on $X'$, the torus $T$ also centralizes $\wt X=X-X'$, and hence $\wt X \in \Lie(T)(E)$ under the identification of $\fg^*$ with $\fg$. Thus $\wt X$ is semisimple, and therefore its $G$-orbit is closed (\cite[9.2]{Borel}). Hence $\wt X \in X+\fg^*_{y,r+}$ is \generic/ of depth $r$ at $y$.

To finish the proof of the proposition, recall that $d(x,\wt X + X')=d(x,X)=r=d(y,\wt X)=d(\wt X)$ (by Corollary \ref{Cor-generic-depth}). We write $G'=\Cent_G(\wt X)$ and $\fg'=\Lie(G')(k)$. Since $X'$ has depth greater than $r$ at $y$ and vanishes on $\bigoplus_{\alpha \in \Phi(G,T)} (\fg_{E})_\alpha$, it lies in $(\fg')^*_{r+} \subset \fg^*$. Hence we deduce from Lemma \ref{Lemma-Kim-Murnaghan} that $x \in \sB(G',k)$. Thus there exists a maximal torus $\wt T$ in $G'\subset G$ with $x \in \sA(\wt T)$, and, by Lemma \ref{Lemma-structure-of-X}, the element $\wt X$ is \generic/ of depth $r$ at $x$. If $\wt X \in X+\fg^*_{x,r+}$, then we are done by choosing $x'=x$ and observing that $\ov X =\ov{\wt X}$. 

Hence it remains to consider the case that $\wt X \notin X+\fg^*_{x,r+}$. Then $d(x,X')=d(x,X-\wt X)=r<d(y,X')\leq d(X')$. Viewing these as depths for $\sB(G',k)$, we deduce from \cite[Corollary~3.2.6]{Adler-DeBacker} (together with their remark at the beginning of Section~3) that the coset $X'+(\fg')^*_{x,r+}$ is degenerate, i.e. contains an unstable element. Hence $\ov{X'} \in ({\fg'}_{x,-r}/{\fg'}_{x,(-r)+})^*$ is unstable by \cite[4.3.~Proposition]{MP1}. Since $\ff$ is perfect, by \cite[Corollary~4.3]{Kempf} there exists a non-trivial one parameter subgroup $\ov \lambda: \bG_m \ra \RP'_x$ in the reductive quotient $\RP'_x$ of $G'$ at $x$ (defined over $\ff$) such that $\lim_{t\ra0}\ov \lambda(t).\ov{X'}=0$. As in the proof of Corollary \ref{Cor-almoststable-and-generic-implies-almoststronglystable}, we let $\bS$ be a maximal split torus of $\RP'_x$ containing $\ov \lambda(\bG_m)$, 
and $\cS$ a split torus (defined over $\cO_k$) in the parahoric group scheme $\bP'_x$ of $G'$ at $x$ whose special fiber is $\bS$ and whose generic fiber $S$ is a split torus in $G'$. This allows us to consider $\ov \lambda$ as an element $\lambda$ of $X_*(S)$. 
  Let $\sA(S, k)$ be the apartment of $S$ (i.e. the apartment of a maximal (maximally split) torus $T_S \subset G'$ containing $S$). Then $\sA(S,k)$ contains $x$ and is the affine space underlying the real vector space $X_*(S) \otimes_{\bZ} \bR$. If $\eps>0$ is sufficiently small, then $X' \in \fg^*_{x+\eps\lambda, r+}$ and $X \equiv \wt X \mod \fg^*_{x+\eps\lambda,r+}$. Let $E'$ be a tamely ramified extension of $k$ over which $T_S$ splits. Then $x':=x+\eps\lambda \in \sA(T_S,E') \cap \sB(G,k)$, and since  $T_S \subset G'=\Cent_G(\wt X)$, the element $\wt X$ is \generic/ at $x'$ of depth $r$ (by Lemma \ref{Lemma-structure-of-X}). 
 \qed

\begin{Aside} 
	The claim proved within the proof of Proposition \ref{Lemma-almoststable-generic-representative} is the dual statement of \cite[Theorem~3.3]{Fi-tame-tori} and could be deduced from the latter as well. We decided to give an independent (but analogous) proof so that the reader has the option to see what assumptions on $p$ enter the claim at which point and observe that in many cases slightly weaker assumptions on $p$ suffice.
\end{Aside}

\section{The datum} \label{Section-datum}
In this section we define the notion of a \datum/ of $G$ and what it means for a \datum/ to be contained in a smooth irreducible representation of $G(k)$. In Section \ref{Section-existence-of-datum} (Theorem \ref{Thm-existence-of-datum}) we will show that every irreducible representation contains such a \datum/. From this result we will deduce in Section \ref{Section-existence-of-type} (Theorem \ref{Thm-exhaustion-of-types}) and Section \ref{Section-exhaustion-Yu} (Theorem \ref{Thm-exhaustion-Yu}) that every irreducible representation contains a \type/ of the form constructed by Kim--Yu (\cite{KimYu}) based on Yu's construction of supercuspidal representations (\cite{Yu}) and that Yu's construction yields all supercuspidal representations.

\begin{Def} \label{Def-extendeddatum}
	Let $n \in \bZ_{\geq 0}$.  An \textit{\extendeddatum/} of $G$ of length $n$ is a tuple 
	$$(x, (r_i)_{1 \leq i \leq n}, (X_i)_{1 \leq i \leq n}, (G_i)_{1 \leq i \leq n+1},(\repzero, \Vrepzero))$$
		where
	\begin{enumerate}[label=(\alph*),ref=\alph*]
		\item \label{Con-a} $r_1 > r_2 > \hdots > r_n >0$ are real numbers
		\item $X_i \in \fg^*_{x,-r_i} \setminus \fg^*_{x,(-r_i)+}$ for $1 \leq i \leq n$ 
		\item $G=G_1 \supseteq G_2 \supsetneq G_3 \supsetneq \hdots \supsetneq G_{n+1}$ are twisted Levi subgroups of $G$
		\item \label{Con-d} $x \in \sB(G_{n+1},k)\subset \sB(G,k)$
		\item $(\repzero, \Vrepzero)$ is an irreducible representation of $(G_{n+1}^{\text{der}})_{x,0}/(G_{n+1}^{\text{der}})_{x,0+}$
	\end{enumerate}
	satisfying the following conditions for all $1 \leq i \leq n$
	\begin{enumerate}[label=(\roman*),ref=\roman*]
		\item \label{Cond-1} $X_i \in \fg_i^*:=\Lie(G_i)(k)^* \subset \fg^*$
		\item $X_i$ is \generic/ of depth $-r_i$ at $x \in \sB(G_i,k)$ as element of $\fg_i^*$ (under the action of $G_i$) 
		\item \label{Cond--1} $G_{i+1}=\Cent_{G_i}(X_i)$
	\end{enumerate}
A \textit{\truncatedextendeddatum/} of $G$ of length $n$ is a tuple 
	$$(x, (r_i)_{1 \leq i \leq n}, (X_i)_{1 \leq i \leq n}, (G_i)_{1 \leq i \leq n+1})$$
	 of data as above satisfying \eqref{Con-a} through \eqref{Con-d} and (\ref{Cond-1}) through (\ref{Cond--1}).
\end{Def}

Note that a \truncatedextendeddatum/ of $G$ of length 0 consists only of a point $x \in \sB(G,k)$ and the group $G$, and an \extendeddatum/ of $G$ of length 0 consists only of a point $x \in \sB(G,k)$, the group $G$, and an irreducible representation of $(G^{\text{der}})_{x,0}/(G^{\text{der}})_{x,0+}$. We are mainly interested in \extendeddata/ of positive length.

\begin{Def} \label{Def-datum}
	Let $n \in \bZ_{\geq 0}$. A \textit{\datum/} of $G$ of length $n$ is a tuple
		$(x, (X_i)_{1 \leq i \leq n}, (\repzero, \Vrepzero))$ consisting of a point $x \in \sB(G,k)$, elements $X_i \in \fg^*$ for $1 \leq i \leq n$ and an irreducible representation $ (\repzero, \Vrepzero)$ of $(\Cent_{G}(\sum_{i=1}^n X_i))^{\text{der}}(k) \cap G_{x,0}/ (\Cent_{G}(\sum_{i=1}^n X_i))^{\text{der}}(k) \cap G_{x,0+}$ for which there exist real numbers $r_1 > r_2 > \hdots > r_n >0$ and a sequence of twisted Levi subgroups $G=G_1 \supseteq G_2 \supsetneq G_3 \supsetneq \hdots \supsetneq G_{n+1}$ of $G$ such that 
$(x, (r_i)_{1 \leq i \leq n}, (X_i)_{1 \leq i \leq n}, (G_i)_{1 \leq i \leq n+1},(\repzero, \Vrepzero))$
is an \extendeddatum/.

A \textit{\truncateddatum/} of $G$ of length $n$ is a tuple
$(x, (X_i)_{1 \leq i \leq n})$ consisting of a point $x \in \sB(G,k)$ and elements $X_i \in \fg^*$ for $1 \leq i \leq n$ for which there exist real numbers $r_1 > r_2 > \hdots > r_n >0$ and a sequence of twisted Levi subgroups $G=G_1 \supseteq G_2 \supsetneq G_3 \supsetneq \hdots \supsetneq G_{n+1}$ of $G$ such that 
$(x, (r_i)_{1 \leq i \leq n}, (X_i)_{1 \leq i \leq n}, (G_i)_{1 \leq i \leq n+1})$
is a \truncatedextendeddatum/.
\end{Def}

Given a  \truncateddatum/  $(x, (X_i)_{1 \leq i \leq n})$ or a \datum/ $(x, (X_i)_{1 \leq i \leq n}, (\repzero, \Vrepzero))$  of $G$ of length $n$, we denote by $(x, (r_i)_{1 \leq i \leq n}, (X_i)_{1 \leq i \leq n}, (G_i)_{1 \leq i \leq n+1})$ the unique \truncatedextendeddatum/ containing it or by $(x, (r_i)_{1 \leq i \leq n}, (X_i)_{1 \leq i \leq n}, (G_i)_{1 \leq i \leq n+1},(\repzero, \Vrepzero))$ the unique \extendeddatum/ containing it, respectively, as in Definition \ref{Def-datum}.

\begin{Rem}
	There are two main differences between a \datum/ and the input for Yu's construction in \cite{Yu}. The first difference is that we only work with elements $X_i \in \fg^*$ and not with characters of $G_{i+1}(k)$. The second difference is that $(\repzero, \Vrepzero)$ is an irreducible representation of $(G_{n+1}^{\text{der}})_{x,0}/(G_{n+1}^{\text{der}})_{x,0+}$ that might not be cuspidal. The representation $(\repzero, \Vrepzero)$ is more a place holder at this point that appears in some sense naturally in Section \ref{Section-existence-of-datum}, and from which we have to extract a cuspidal representation that forms the input for Yu's construction in Section \ref{Section-existence-of-type}  (Lemma \ref{Lemma-repzeroK} and Lemma \ref{Lemma-cuspidal}).
	Thus our \datum/ can be viewed as a skeleton of the input for Yu's construction. 
\end{Rem}

For later convenience, we note the following lemma.
\begin{Lemma} \label{Lemma-truncateddatum-with-y}
	If $(x, (X_i)_{1 \leq i \leq n})$ is a \truncateddatum/ of $G$, and $y$ is a point of $\sB(G_{n+1},k) \subset \sB(G,k)$, then $(y, (X_i)_{1 \leq i \leq n})$ is also a \truncateddatum/ of $G$.
\end{Lemma}
\Proof
This follows from Lemma \ref{Lemma-structure-of-X}. 
\qed

In order to relate a \truncateddatum/  $(x, (X_i)_{1 \leq i \leq n})$ or a \datum/ $(x, (X_i)_{1 \leq i \leq n}, (\repzero, \Vrepzero))$ to representations of $G(k)$, we introduce the following associated groups for $1 \leq i \leq n+1$: 
\begin{itemize}
	\item $H_1:=G_1$ if $G_1=G_2$ and $H_1:=G_1^{\text{der}}$ if $G_1 \neq G_2$
	\item $H_i:=G_i^{\text{der}}$ for $i>1$
	\item $(H_{i})_{x,\wt r}:=G_{x,\wt r}\cap H_i(k)=(H_i)_{x_i,\wt r}$ for $\wt r \in \wt \bR_{\geq 0}:=\bR_{\geq 0} \cup \{ r+ \, | \, r \in \bR_{\geq 0} \}$, 
\end{itemize}
where $x_i$ denotes the image of $x \in \sB(G_i,k)$ in $\sB(H_i,k)$. 
In order to define another subgroup $(H_{i})_{x,\wt r,\wt r'}$ of $G(k)$ for $\wt r \geq \wt r' \geq \frac{\wt r}{2}>0$ ($\wt r, \wt r' \in \wt \bR$) and $1 \leq i \leq n$, we choose a maximal torus $T$ of $G_{i+1}$ such that $x \in \sA(T,E)$, where $E$ denotes a finite tamely ramified extension of $k$ over which $T$ splits. Then we define
$$ (G_{i})_{x,\wt r,\wt r'} := G(k) \cap \<T(E)_{\wt r}, U_\alpha(E)_{x,\wt r}, U_\beta(E)_{x,\wt r'} \, | \, \alpha \in \Phi(G_i, T)\subset\Phi(G,T), \beta \in \Phi(G_i, T)-\Phi(G_{i+1}, T) \, \>,$$
where 
$U_\alpha(E)_{x,r}$ denotes the Moy--Prasad filtration subgroup of depth $r$ (at $x$) of the root group $U_\alpha(E) \subset G(E)$ corresponding to the root $\alpha$, 
and 
$$ (H_{i})_{x,\wt r, \wt r'} := H_i(k) \cap (G_{i})_{x,\wt r, \wt r'} $$
Note that $(G_{i})_{x,\wt r, \wt r'}$ is denoted $(G_{i+1},G_i)(k)_{x_i,\wt r, \wt r'}$ in \cite{Yu}. Yu (\cite[p.~585 and p.~586]{Yu}) shows that this definition is independent of the choice of $T$ and $E$.
We define the subalgebras $\fh_{i}$, $(\fh_{i})_{x,\wt r}$ and $(\fh_{i})_{x,\wt r, \wt r'}$ of $\fg$ analogously. For convenience, we also set $r_{n+1}=0$.

\begin{Def} \label{Def-datum-contained}
	Let $(\pi, V_\pi)$ be a smooth irreducible representation of $G(k)$. 
	A \datum/ 
	$(x, (X_i)_{1 \leq i \leq n}, (\repzero, \Vrepzero))$
	 of $G$ is said to be \textit{contained in} $(\pi, V_\pi)$ if 
$V_\pi^{\cup_{1 \leq i \leq n+1}((H_{i})_{x,r_i+})}$ contains a subspace $V'$ such that
\begin{itemize}
	\item $(\pi|_{(H_{n+1})_{x,0}}, V')$ is isomorphic to $(\repzero, \Vrepzero)$ as a representation of $(H_{n+1})_{x,0}/(H_{n+1})_{x,0+}$ and 
	\item $(H_{i})_{x,r_i,\frac{r_i}{2}+} /(H_{i})_{x,r_i+} \simeq (\fh_{i})_{x,r_i,\frac{r_i}{2}+} /(\fh_{i})_{x,r_i+}$ acts on $V'$ via the character $\varphi \circ X_i$ for $1\leq i \leq n$, 
	\end{itemize} 
where we recall that $\varphi: k \ra \bC^*$ is an additive character of $k$ of conductor $\cP$ that is fixed throughout the paper. 

	Similarly,  $(\pi,V_\pi)$ is said to contain a \truncateddatum/ 
	$(x, (X_i)_{1 \leq i \leq n})$ if 
	$V_\pi^{\cup_{1 \leq i \leq n}((H_{i})_{x,r_i+})}$ contains a one dimensional subspace on which
$(H_{i})_{x,r_i,\frac{r_i}{2}+} /(H_{i})_{x,r_i+} \simeq (\fh_{i})_{x,r_i,\frac{r_i}{2}+} /(\fh_{i})_{x,r_i+}$ acts via 	$\varphi \circ X_i$ for $1 \leq i \leq n$.
		
\end{Def}

The \data/ that we are going to use to extract a type from a given representation are the following.

\begin{Def} \label{Def-datum-for}
	Let $(\pi, V_\pi)$ be a smooth, irreducible representation of $G(k)$. 
	We say that a tuple
	$(x, (X_i)_{1 \leq i \leq n}, (\repzero, \Vrepzero))$
	\textit{is a maximal \datum/ for} $(\pi, V_\pi)$ if  $(x, (X_i)_{1 \leq i \leq n}, (\repzero, \Vrepzero))$ is a \datum/ of $G$ that is contained in $(\pi,V_\pi)$ such that if $(x', (X_i)_{1 \leq i \leq n'}, (\repzero', \Vrepzero'))$ is another \datum/ of $G$ contained in $(\pi,V_\pi)$, then the dimension of the facet of $\sB(G_{n+1},k)$ that contains $x$ is at least the dimension of the facet of $\sB(G_{n+1},k)$ that contains $x'$. 
\end{Def}

\section{Some results used to exhibit \data/} \label{Section-killing-lemmas}
In order to prove that every irreducible representation of $G(k)$ contains a \datum/, we first prove a lemma and derive some corollaries that we are going to repeatedly use in the proof of the existence of a \datum/ in Section \ref{Section-existence-of-datum}. (The reader might skip this section at first reading and come back to it when the results are used in the proof of Theorem \ref{Thm-existence-of-datum}.)

\begin{Lemma} \label{Lemma-kill-C} 
	Let $x \in \sB(G,k)$, $r \in \bR_{>0}$, and let $X \in \fg^*-\{0\}$ be \generic/ of depth $-r$ at $x$. Write $G'=\Cent_G(X)$, $\fg'= \Lie(G')(k)$, and let $T$ be a maximal torus of $G'$ that splits over a tamely ramified extension $E$ of $k$ and such that $x \in \sA(T,E) \cap \sB(G,k)$. We set $\ft=\Lie(T)(k)$, $\fr'=\fg \cap \left(\bigoplus_{\alpha \in \Phi(G',T_E)} (\fg_E)_\alpha \right) \subset \fg'$ and $\fr''=\fg \cap \left(\bigoplus_{\alpha \in \Phi(G,T_E)-\Phi(G',T_E)} (\fg_E)_\alpha \right)$, and we denote by $\fj^*$ the subspace of elements in $\fg^*$ that vanish on $\ft \oplus \fr' =\fg'$. Then
	\begin{enumerate}[label=(\alph*),ref=\alph*]
		\item \label{item-fX} The map $f_X: \fr'' \ra \fg^*$ defined by $Y \mapsto (Z \mapsto X([Y,Z]))$ is a vector space isomorphism of $\fr''$ onto $\fj^*$, and $f_X(\fr'' \cap \fg_{x,r'}) = \fj^*\cap\fg^*_{x,r'-r}$ for $r' \in \bR$. 
		\item \label{item-kill-C} Let $d$ be a real number such that $\frac{r}{2}\leq d < r$. For every $0<\eps<\frac{r-d}{2}$, if $C \in \fj^* \cap \fg^*_{x,-(d+\eps)}$, then there exists $g \in G_{x,r-d-\eps} \subset G_{x,0+}$ such that
		\begin{enumerate}[label=(\roman*),ref=\roman*]
			\item \label{item-kill-C-i} $\Ad(g)(X+C)|_{\fg_{x,r}} =  X|_{\fg_{x,r}}$,
			\item  $\Ad(g)(X+C)|_{\fr'' \cap \fg_{x,d+}} =  0 =  X|_{\fr'' \cap \fg_{x,d+}}$,
			\item \label{item-kill-C-iii} if $(\pi, V_\pi)$ is a representation of $G(k)$ and $V'$ is a subspace of $V_\pi$ 
			on which the group 
			$$ G(k) \cap \<T(E)_{r+}, U_\alpha(E)_{x,(d+\eps)+}, U_\beta(E)_{x,r+} \, | \, \alpha \in \Phi(G,T_E)-\Phi(G',T_E), \beta \in \Phi(G',T_E) \> $$
			acts trivially and 
			that is stable under the action of a subgroup $H$ of $G'^{\text{der}}(k) \cap G_{x,(2d-r+2\eps)+}$,
			then $g^{-1}Hg$ preserves $V'$ and 
			$(^{g}\pi|_{H},V') = (\pi|_{H},V'),$ 
		\end{enumerate}
		where $\Ad$ denotes the contragredient of the adjoint action. 
	\end{enumerate}
\end{Lemma}
\Proof 
\eqref{item-fX} 
Let $Y \in \fr''$. Recall that $[(\fg_E)_\alpha, (\fg_E)_\beta] \subset (\fg_E)_{\alpha+\beta}$ for $\alpha, \beta \in \Phi(G, T_E)$, $\alpha \neq -\beta$ (where $(\fg_E)_{\alpha+\beta}=\{0\}$ if $\alpha+\beta \notin \Phi(G,T_E)$). Hence, if $Z \in \fr'$, then $[Y,Z] \in \fr''$, and $X([Y,Z])=0$ by Lemma \ref{Lemma-structure-of-X}. Similarly, if $Z \in \ft$, then $[Y,Z] \in \fr''$, and $X([Y,Z])=0$. Thus the image of the linear map $f_X$ is contained in $\fj^*\subset \fg^*$.

Choose a Chevalley system $\{x_\alpha:\bG_a \ra G_{E} \, | \, \alpha \in \Phi(G,T_E) \}$  for $G_E$ with corresponding Lie algebra elements $\{X_{\alpha}=dx_\alpha(1) \, | \alpha \in \Phi(G,T_E)\}$. Then for $\alpha \in \Phi(G, T_E)-\Phi(G',T_E)$, we have $[X_\alpha,X_{-\alpha}]=H_\alpha=d\check\alpha(1)$. Hence, extending $f_X$ linearly to $\fr'' \otimes_k E$, the element $f_X(X_\alpha)$ in $\fj^* \otimes_k E$ is a map that sends $X_\beta$ to $c\delta_{-\alpha,\beta}$ for $\beta \in \Phi(G,T_E)-\Phi(G',T_E)$ for some constant $c \in E$ with $\val(c)=\val(X(H_\alpha))=-r$, by Lemma \ref{Lemma-structure-of-X}. From this description, we see that $d(x,f_X(X_\alpha))=-r-d(x,X_{-\alpha})=-r-(-\alpha(x))=\alpha(x)-r$ while $d(x,X_\alpha)=\alpha(x)$. Thus $f_X(\fr'' \otimes_k E \cap (\fg_E)_{x,r'})=\fj^* \otimes_k E \cap (\fg^*_E)_{x,r'-r}$, and hence $f_X(\fr''  \cap \fg_{x,r'})=\fj^* \cap (\fg^*)_{x,r'-r}$ because $E$ is tamely ramified over $k$. In particular, $f_X: \fr'' \ra \fj^*$ is a vector space isomorphism.

\eqref{item-kill-C}
By \eqref{item-fX}, there exists $Y \in \fr''$ of depth $\geq r-d-\eps> 0$ such that $C=X([Y,\_])$. Let $\exp$ denote a mock exponential function from $\fg_{x,r-d-\eps}$ to $G_{x,r-d-\eps}$ as defined in \cite[Section~1.5]{Adler}, i.e. if $Y=\sum_{\alpha \in \Phi(G,T_E)-\Phi(G',T_E)} a_\alpha X_{\alpha}$ for some $a_\alpha \in E$, then 
$$\exp(Y) \equiv \prod_{\alpha \in \Phi(G,T_E)-\Phi(G',T_E)} x_\alpha(a_\alpha) \mod G(E)_{x,2(r-d-\eps)}$$ (viewing $\exp(Y) \in G(k)$ inside $G(E)$) for some fixed (arbitrarily chosen) order of the roots $\Phi(G,T_E)-\Phi(G',T_E)$. We set $g=(\exp(-Y))^{-1}$. Note that $Y \in \fg_{x,r-d-\eps}$ implies that $x_\alpha(-a_\alpha) \in G(E)_{x,r-d-\eps}$.

Let $Z \in \fg_{x,r'}$ for some $r'\in \bR$. Then by \cite[Proposition~1.6.3]{Adler}, we have
$$ \Ad(g^{-1})(Z) \equiv Z + [-Y,Z] \mod \fg_{x,r'+2(r-d-\eps)} .$$
Hence, using that $g \in G_{x,r-d-\eps}, X \in \fg^*_{x,-r}, C \in \fg^*_{x,-(d+\eps)}$, $Y \in \fg_{x,r-d-\eps}$ and $\eps < \frac{r-d}{2}$, we obtain
$$\Ad(g)(X+C)(Z)=X(\Ad(g^{-1})(Z))+C(\Ad(g^{-1})(Z))=X(Z) \quad \text{for $Z \in  \fg_{x,r}$, and}$$
$$\Ad(g)(X+C)(Z)=X(\Ad(g^{-1})(Z))+C(\Ad(g^{-1})(Z))=X(Z+[-Y,Z])+C(Z)=X(Z)$$
for $Z \in  \fr'' \cap \fg_{x,d+}$ (using $C=X([Y,\_])$).

To prove the remaining claim, observe that for $h \in H \subset G'(k) \cap G_{x,(2d-r+2\eps)+}$, $\alpha \in \Phi(G,T_E)-\Phi(G',T_E)$ and $x_\alpha(a_\alpha) \in G(E)_{x,r-d-\eps}$, we have
$x_\alpha(a_\alpha) h x_\alpha(a_\alpha)^{-1} = h \cdot u $ for some $u$ in $\<U_\alpha(E)_{x,(d+\eps)+}\, | \, \alpha \in \Phi(G,T_E)-\Phi(G',T_E) \>$,
and hence
$$g h g^{-1} \equiv h \cdot u' \mod G(E)_{x,r+}  $$
for some $u'$ in $\<U_\alpha(E)_{x,(d+\eps)+}\, | \, \alpha \in \Phi(G,T_E)-\Phi(G',T_E) \>$.
Thus $\pi(g h g^{-1})$ and $\pi(h)$ agree on $V'$.
\qed

\begin{Cor} \label{Cor-kill-C-new} 
	Let $n$ be a positive integer and $(x,(X_i)_{1 \leq i \leq n})$ a \truncateddatum/ of length $n$ (with corresponding \truncatedextendeddatum/ $(x,(r_i)_{1 \leq i \leq n}, (X_i)_{1 \leq i \leq n}, (G_i)_{1 \leq i \leq n+1})$). Let $1 \leq j \leq n$, choose a maximal torus $T_j$ of $\Cent_{H_j}(X_j) \subset H_j$ that splits over a tame extension $E$ of $k$ and such that $x \in \sA(T_j,E) \cap \sB(H_j,k)$. Write $\ft_j=\Lie(T_j)(k)$. Let $(\pi, V_\pi)$ be a representation of $G(k)$. Let $d, \eps \in \bR$ such that $\frac{r_j}{2} \leq d < r_j$ and  $\frac{r_j-d}{2}>\eps>0$. Suppose that $V'$ is a nontrivial subspace of $V_\pi^{\bigcup_{1 \leq i \leq j}(H_i)_{x,r_i+}}$ on which
	$(H_i)_{x,r_i,\frac{r_i}{2}+}/(H_i)_{x,r_i+}\simeq (\fh_i)_{x,r_i,\frac{r_i}{2}+}/(\fh_i)_{x,r_i+}$ acts via $\varphi \circ X_i$ for $1 \leq i \leq j-1$ 
	and 
	$(H_j)_{x,r_j-\eps,d+}/(H_j)_{x,r_j+}\simeq (\fh_j)_{x,r_j-\eps,d+}/(\fh_j)_{x,r_j+}$ acts via $\varphi \circ (X_j+C)$ for some $C \in (\fh_j^*)_{x,-(d+\eps)}$ that is trivial on $\ft_j + \fh_{j+1}$. Then there exists $g \in (H_j)_{x,r_j-d-\eps}$ such that 
	\begin{enumerate}[label=(\roman*),ref=\roman*]
		\item $V'':=\pi(g)V' \subset V_\pi^{\bigcup_{1 \leq i \leq j}(H_i)_{x,r_i+}}$
		\item  \label{item-blubb2} $(H_i)_{x,r_i,\frac{r_i}{2}+}/(H_i)_{x,r_i+}\simeq (\fh_i)_{x,r_i,\frac{r_i}{2}+}/(\fh_i)_{x,r_i+}$ acts on $ V''$  via $\varphi \circ X_i$ for $1 \leq i \leq j-1$ 
		\item  $(H_j)_{x,r_j,d+}/(H_j)_{x,r_j+}\simeq (\fh_j)_{x,r_j,d+}/(\fh_j)_{x,r_j+}$ acts on $ V''$  via $\varphi \circ X_j$ 
		\item any subgroup $H$ of $(H_{j+1})_{x,(2d-r_j+2\eps)+}$ that stabilizes $V'$ also stabilizes $V''$ and $(\pi|_{H},V'') \simeq (\pi|_{H},V')$.
	\end{enumerate}
\end{Cor}	
\Proof
Let $g \in (H_j)_{x,r_j-d-\eps}$ be as constructed in the proof of Lemma \ref{Lemma-kill-C}\eqref{item-kill-C} applied to the group $H_j$ with  \generic/ element $X_j$ of depth $r_j$. Hence $(\pi|_{H},V'') \simeq (\pi|_{H},V')$ by Lemma \ref{Lemma-kill-C}\eqref{item-kill-C}\eqref{item-kill-C-iii}. Note that $g^{-1}((H_i)_{x,r_i+})g=(H_i)_{x,r_i+}$ for $1 \leq i \leq j$. Thus $V'' \subset V_\pi^{\bigcup_{1 \leq i \leq j}(H_i)_{x,r_i+}}$. To show \eqref{item-blubb2}, recall that $g^{-1}=\prod_{\alpha \in \Phi(H_j,(T_j)_E)-\Phi(\Cent_{H_j}(X_j),(T_j)_E)} x_\alpha(-a_\alpha)g'$ with $g'\in H_j(E)_{x,2r_j-2d-2\eps}$ and $x_\alpha(-a_\alpha) \in H_j(E)_{x,r_j-d-\eps}$. Hence for $h \in (H_i)_{x,r_i,\frac{r_i}{2}+}$ for $1 \leq i <j$ we have $g^{-1}hg \equiv h \mod (H_i)_{x,r_i+,(\frac{r_i}{2}+r_j-d-\eps)+}$ with $r_j-d-\eps>0$, and therefore $(H_i)_{x,r_i,\frac{r_i}{2}+}/(H_i)_{x,r_i+}\simeq (\fh_i)_{x,r_i,\frac{r_i}{2}+}/(\fh_i)_{x,r_i+}$ acts on $ V''$  via $\varphi \circ X_i$. In addition, we have $g^{-1}(H_j)_{x,r_j,d+}g\subset(H_j)_{x,r_j-\eps,d+}$. Since  $(H_j)_{x,r_i-\eps,d+}/(H_j)_{x,r_j+}\simeq (\fh_j)_{x,r_j-\eps,d+}/(\fh_j)_{x,r_j+}$ acts via $\varphi \circ (X_j+C)$ on $V'$, we obtain from Lemma \ref{Lemma-kill-C}\eqref{item-kill-C} that $(H_j)_{x,r_j,d+}/(H_j)_{x,r_j+}\simeq (\fh_j)_{x,r_j,d+}/(\fh_j)_{x,r_ji+}$ acts on $ V''$  via $\varphi \circ X_j$. \qed

\begin{Cor} \label{Cor-kill-C-type1}
		Let $n$ be a positive integer and $(x,(X_i)_{1 \leq i \leq n})$ a \truncateddatum/ of length $n$ (with corresponding \truncatedextendeddatum/ $(x,(r_i)_{1 \leq i \leq n}, (X_i)_{1 \leq i \leq n}, (G_i)_{1 \leq i \leq n+1})$). Let $T$ be a maximal torus of $G_{n+1}$ such that $x\in \sA(T,E)$, set $\ft=\Lie(T)(k)$, and let $(\pi, V_\pi)$ be a representation of $G(k)$. Let $0<\eps < \frac{r_n}{4}$ such that $(H_{i})_{x,r_i-\eps,\frac{r_i}{2}+}=(H_{i})_{x,r_i,\frac{r_i}{2}+}$ for all $1 \leq i \leq n$. Suppose 	
		 that $V'$ is a nontrivial subspace of $V_\pi^{\bigcup_{1 \leq i \leq n}(H_i)_{x,r_i+}}$ on which the action of 
		$(H_{i})_{x,r_i,\frac{r_i}{2}+} /(H_{i})_{x,r_i+} \simeq (\fh_{i})_{x,r_i,\frac{r_i}{2}+} /(\fh_{i})_{x,r_i+}$ via $\pi$ is given by
		$\varphi \circ \left(C_i+ X_i\right)$
		for some $C_i \in (\fh_i)^*_{x,-(\frac{r_i}{2}+\eps)}$ that is trivial on $(\ft \cap \fh_i) + \fh_{i+1}$ for all $1 \leq i \leq n$.
		Then there exists $g \in G_{x,\frac{r_n}{2}-\eps}\subset G_{x,0+}$ such that 
		\begin{enumerate}[label=(\roman*),ref=\roman*]
			\item $V'':=\pi(g)V' \subset V_\pi^{\bigcup_{1 \leq i \leq n}(H_i)_{x,r_i+}}$
			\item  $(H_i)_{x,r_i,\frac{r_i}{2}+}/(H_i)_{x,r_i+}\simeq (\fh_i)_{x,r_i,\frac{r_i}{2}+}/(\fh_i)_{x,r_i+}$ acts on $ V''$  via $\varphi \circ X_i$ for $1 \leq i \leq n$ 
			\item any subgroup $H$ of $(H_{n+1})_{x,2\eps+}$ that stabilizes $V'$ also stabilizes $V''$ and $(\pi|_{H},V'') \simeq (\pi|_{H},V')$.
		\end{enumerate}
\end{Cor}
\Proof 
Since $(H_{i})_{x,r_i-\eps,\frac{r_i}{2}+}=(H_{i})_{x,r_i,\frac{r_i}{2}+}$ for all $1 \leq i \leq n$, we can apply Corollary \ref{Cor-kill-C-new} successively for $j=1, 2, \hdots, n$  with $d=\frac{r_1}{2}, \frac{r_2}{2}, \hdots , \frac{r_n}{2}$, respectively. We obtain $g=g_{n} \cdot \hdots \cdot g_1\in (H_{n})_{x,\frac{r_{n}}{2}-\eps} \cdot \hdots \cdot (H_1)_{x,\frac{r_1}{2}-\eps}\subset G_{x,\frac{r_n}{2}-\eps}$ such that $\pi(g)V' \subset V_\pi^{\cup_{1 \leq i \leq j-1}((H_i)_{x,r_i+})}$ and the action of $(H_{i})_{x,r_i,\frac{r_i}{2}+} /(H_{i})_{x,r_i+} \simeq (\fh_{i})_{x,r_i,\frac{r_i}{2}+} /(\fh_{i})_{x,r_i+}$ on $\pi(g)V'$ via $\pi$ is given by
$\varphi \circ X_i$ for $1\leq i \leq n$,
and the action of any subgroup $H$ of $(H_{n+1})_{x,2\eps+}$ that stabilizes $V'$ also stabilizes $V''$ and $(\pi|_{H},V'') \simeq (\pi|_{H},V')$. \qed

\begin{Cor} \label{Cor-kill-C-type2}
		Let $n$ be a positive integer and $(x,(X_i)_{1 \leq i \leq n})$ a \truncateddatum/ of length $n$ (with corresponding \truncatedextendeddatum/ $(x,(r_i)_{1 \leq i \leq n}, (X_i)_{1 \leq i \leq n}, (G_i)_{1 \leq i \leq n+1})$). Let $(\pi, V_\pi)$ be a representation of $G(k)$. Let $0<\eps < \frac{r_n}{4}$ such that $(H_{n})_{x,r_n-2\eps}=(H_{n})_{x,r_n}$. Suppose 	
		that $V'$ is a nontrivial subspace of $V_\pi^{\bigcup_{1 \leq i \leq n}(H_i)_{x,r_i+}}$ on which the action of 
		$(H_{i})_{x,r_i,\frac{r_i}{2}+} /(H_{i})_{x,r_i+} \simeq (\fh_{i})_{x,r_i,\frac{r_i}{2}+} /(\fh_{i})_{x,r_i+}$ via $\pi$ is given by
		$\varphi \circ X_i$
		for all $1 \leq i \leq n-1$, and
		the action of 
		$(H_{n})_{x,r_n} /(H_{n})_{x,r_n+} \simeq (\fh_{n})_{x,r_n} /(\fh_{n})_{x,r_n+}$ via $\pi$ is given by
	 	$\varphi \circ X_n$.
		Then there exists a subspace $V'' \subset V_\pi^{\bigcup_{1 \leq i \leq n}(H_i)_{x,r_i+}}$ such that 
			 $(H_i)_{x,r_i,\frac{r_i}{2}+}/(H_i)_{x,r_i+}\simeq (\fh_i)_{x,r_i,\frac{r_i}{2}+}/(\fh_i)_{x,r_i+}$ acts on $ V''$  via $\varphi \circ X_i$ for $1 \leq i \leq n$. 
\end{Cor}
\Proof
Let $T$ be a maximal torus of $G_{n+1}$ with Lie algebra $\ft=\Lie(T)(k)$. Let $d=\max(\frac{r_n}{2},r_n-3\eps)$. 
Note that the commutator $[\prod_{1 \leq i \leq n-1} (H_i)_{x,r_i,\frac{r_i}{2}+}(H_{n})_{x,r_n},(H_{n})_{x,r_n,d+}]$ 
acts trivially on $V'$ 
and hence we can replace $V'$ without loss of generality by $\pi((H_{n})_{x,r_n,d+})V'$.
Since $(H_{n})_{x,r_n-2\eps}=(H_{n})_{x,r_n}$, the action of $(H_{n})_{x,r_n-\eps,d+}$ on $V'$ factors through $(H_{n})_{x,r_n-\eps,d+}/(H_{n})_{x,r_n+}$ and, after replacing $V'$ by a subspace if necessary, 
is given by $\varphi \circ \left(X_n + C_3 \right)$ for some $C_3 \in (\fh_n)^*_{x,-(r_n-2\eps)}$ that is trivial on $\ft \cap \fh_n + \fh_{n+1}$. Applying Corollary \ref{Cor-kill-C-new} for $j=n$ and $d=\max(\frac{r_n}{2},r_n-3\eps)$, we obtain $g_3 \in (H_n)_{x,\eps}$ such that $\pi(g_3)V'\subset V_\pi^{\bigcup_{1 \leq i \leq n}(H_i)_{x,r_i+}}$, the group $(H_{i})_{x,r_i,\frac{r_i}{2}+} /(H_{i})_{x,r_i+} \simeq (\fh_{i})_{x,r_i,\frac{r_i}{2}+} /(\fh_{i})_{x,r_i+}$ acts on $\pi(g_3)V'$ via $\varphi \circ X_i$ for $1 \leq i \leq n-1$ and $(H_{n})_{x,r_n,d+} /(H_{n})_{x,r_n+} \simeq (\fh_{n})_{x,r_n, d+} /(\fh_{n})_{x,r_n+}$ acts on $\pi(g_3)V'$ via $\varphi \circ X_n$. Replacing $V'$ by  $\pi(g_3)V'$ and using the same reasoning, we can apply Corollary \ref{Cor-kill-C-new} for $j=n$ repeatedly with $d= r_n - 4 \eps , r_n - 5 \eps, r_n-6\eps , \hdots, r_n - (N-1) \cdot \eps, r_n - N \cdot \eps,\frac{r_n}{2}$ (and replacing $V'$ at each step if necessary), where $N$ is the largest integer for which $N \cdot \eps < \frac{r_n}{2}$. 
After the final step we obtain a subspace $V'' \subset V_\pi^{\bigcup_{1 \leq i \leq n}(H_i)_{x,r_i+}}$ on which $(H_i)_{x,r_i,\frac{r_i}{2}+}/(H_i)_{x,r_i+}\simeq (\fh_i)_{x,r_i,\frac{r_i}{2}+}/(\fh_i)_{x,r_i+}$ acts via $\varphi \circ X_i$ for $1 \leq i \leq n$. \qed

\section{Every irreducible representation contains a \datum/} \label{Section-existence-of-datum}

\begin{Thm} \label{Thm-existence-of-datum} Let $(\pi, V_\pi)$ be a smooth irreducible representation of $G(k)$.  Then $(\pi, V_\pi)$ contains a \datum/.
	\end{Thm}

\Proof 
The strategy of this proof consists of recursively extending the length of a \truncateddatum/ contained in $(\pi , V_\pi)$ until a certain function $f$ defined below is zero. We then show how to turn this \truncateddatum/ into a \datum/ that is contained in $(\pi, V_\pi)$.

More precisely, for the recursion step, we let $j$ be a positive integer such that $(\pi,V_\pi)$ contains a \truncateddatum/ $(x_{j-1}, (X_i)_{1 \leq i \leq j-1})$
  of $G$ of length $j-1$. We will then show that  $(\pi,V_\pi)$ contains a \truncateddatum/ $(x_{j}, (X_i)_{1 \leq i \leq j})$ of $G$ of length $j$ and repeat the recursion or that $(\pi,V_\pi)$ contains a \datum/ $(x_{j-1}, (X_i)_{1 \leq i \leq j-1}, (\repzero, \Vrepzero))$ and the proof is finished. Since $G_i \subsetneq G_{i-1}$ for $i>2$, the recursion has to terminate after finitely many steps. 
  
The base case of the recursion is given by $j=1$, in which case 
 we let $x_0=x_{j-1}$ be an arbitrary point of $\sB=\sB(G,k)$ and denote by $r_0=r_{j-1}$ the depth of $(\pi,V_\pi)$ at $x_0$.

To perform the recursion step, let $j$ be a positive integer such that $(\pi,V_\pi)$ contains a \truncateddatum/ $(x_{j-1}, (X_i)_{1 \leq i \leq j-1})$
of $G$ of length $j-1$, and write $\sB_j:=\sB(G_j,k) \subset \sB$, where the inclusion of Bruhat--Tits buildings is as explained in Remark \ref{Rem-BT}. 
We define a function $f:\sB_j \ra \bR_{\geq 0} \cup \infty$ as follows:
For $y \in \sB_j$, we set $f(y)$ to be the smallest non-negative real number $r_j$ such that 
\begin{itemize}
	\item the \truncateddatum/  
$(y, (X_i)_{1 \leq i \leq j-1})$ is contained in $(\pi,V_\pi)$ 
\item there exists $X_j \in (\fg_j)^*_{y,-r_j}$ \almoststable/, where $\fg_j=\Lie(G_j)(k)$, and
\item there exists $V_{j-1} \subset V_\pi^{\cup_{1 \leq i \leq j-1}({(H_{i})_{y,r_i+}})}$
\end{itemize}
 satisfying the following two properties
\begin{enumerate}[label=(\alph*),ref=\alph*]
\item \label{item-in-pfa} for $1\leq i \leq j-1$ the group $$(H_{i})_{y,r_i,\frac{r_i}{2}+} /(H_{i})_{y,r_i+} \simeq (\fh_{i})_{y,r_i,\frac{r_i}{2}+} /(\fh_{i})_{y,r_i+}$$ 
 associated to 
$(y, (X_i)_{1 \leq i \leq j-1})$ 
acts 
on $V_{j-1}$ via  $\varphi \circ X_i$ (this condition is automatically satisfied for $j=1$), and 
\item \label{item-in-pfb} $V_{j-1}^{(H_j)_{y,r_j+}}$ contains a nontrivial subspace $V_j'$ such that if $r_j>0$, then  $(H_j)_{y,r_j}/(H_j)_{y,r_j+} \simeq (\fh_j)_{y,r_j}/(\fh_j)_{y,r_j+}$ acts on $V_j'$ via $\varphi \circ X_j$.
\end{enumerate}
If such a real number $r_j$ does not exist, then we set $f(y)=\infty$. 

Note that $(y, (X_i)_{1 \leq i \leq j-1})$ is a \truncateddatum/ of $G$ by Lemma \ref{Lemma-truncateddatum-with-y}. Moreover, $f$ is well defined, because the Moy--Prasad filtration is semi-continuous and for every $r \in \bR$ every $(\fg_j)^*_{y,r+}$-coset contains an \almoststable/ element (e.g. take an element dual to a semisimple element under the non-degenerate bilinear form $B$ provided by \cite{Adler-Roche}, see Remark \ref{Rem-B}).
In addition, by our assumption, $f(x_{j-1}) \leq r_{j-1}$ (because if $j>1$, we could take $X_j=0$ for $r_j=r_{j-1}$). In the case $j=1$, the real number $f(y)$ is simply the depth of $(\pi, V_\pi)$ at $y$.

\begin{Lemmasub} \label{Lemma-Aj}
	\begin{enumerate}[label=(\roman*),ref=\roman*]
		\item \label{Lemma-Aj-ii} $f(g.x)=f(x)$ for all $x \in \sB_j$ and $g \in G_j(k)$
		\item \label{Lemma-Aj-iii} The subset $f^{-1}(\bR_{\geq 0})$ of $\sB_j$  is open in $\sB_j$ and the function $f:\sB_j \ra \bR_{\geq 0} \cup \infty$ is continuous on $f^{-1}(\bR_{\geq 0})$.
		\item \label{Lemma-Aj-iv} The subset $f^{-1}(\bR_{\geq 0})$ of $\sB_j$  is closed in $\sB_j$, hence equal to $\sB_j$.		 
	\end{enumerate}
\end{Lemmasub}
\textbf{Proof of Lemma \ref{Lemma-Aj}.}\\
\textit{Proof of part (\ref{Lemma-Aj-ii}).} Observe that $X_i$ ($1\leq i <j$) and $G_i$ ($1 \leq i \leq j$) are stabilized by $G_j(k)$, hence the $G_j(k)$-invariance of $f$ follows. 

\textit{Proof of part (\ref{Lemma-Aj-iii}).}
If $j=1$, then $f(x)$ is the depth of $\pi$ at $x$, and the claim is true. Hence we assume $j>1$.
Let $(x, (X_i)_{1 \leq i \leq j-1})$ be a \truncateddatum/ contained in $(\pi,V_\pi)$, $X_j \in (\fg_j)^*_{x,-f(x)}$ \almoststable/ and $V_{j}' \subset V_{j-1} \subset V_\pi^{\cup_{1 \leq i \leq j-1}((H_i)_{x,r_i+})}$ satisfying the conditions \eqref{item-in-pfa} and \eqref{item-in-pfb} above. If $f(x)>0$, then set $r_j=f(x)$, otherwise let $0<r_j \leq r_{j-1}$ be arbitrary. For $1 \leq i \leq j-1$, let $d_{i}<r_i$ be a positive real number such that $(G_i)_{x,r_i-d_i}=(G_i)_{x,r_i}$.
Note that $d_i>0$ exists for $1 \leq i \leq j-1$ by the semi-continuity of the Moy--Prasad filtration.

Let $\min\{\frac{r_j}{4},\frac{d_i}{2} \, | \, 1 \leq i <j-1\}>\eps>0$ and let $y \in \sB_j$ with $\mathrm{d}(x,y)<\eps$. 
Let $T$ be a maximal torus of $G_j$ that splits over a tamely ramified extension $E$ of $k$ such that $x$ and $y$ are contained in $\sA(T_E,E)$. Then $(y,(X_i)_{1 \leq i \leq j-1})$ is a \truncateddatum/ by Lemma \ref{Lemma-truncateddatum-with-y}.  By the normalization of the distance $\mathrm{d}$ on the building $\sB$, we have $\abs{\alpha(x-y)}\leq \mathrm{d}(x,y)<\eps$ for all $\alpha \in \Phi(G,T_E)$. Hence, since $X_i$ vanishes on $\fg_i \cap \bigoplus_{\alpha \in \Phi(G_i,T_E)} \fg(E)_{\alpha}$ for $1 \leq i \leq j-1$, we have  $V_j'\subset V_\pi^{\cup_{1 \leq i \leq j-1}((H_i)_{y,r_i+})}$ and the commutator $\left[\prod_{1 \leq i \leq j-1} (H_{i})_{y,r_i-\frac{d_i}{2},\frac{r_i}{2}+}(H_j)_{y,r_j+\eps},\right.$ $\left. \prod_{1 \leq i \leq j-1} (H_{i})_{y,r_i-\frac{d_i}{2},\frac{r_i}{2}+}(H_j)_{y,r_j+\eps}\right]$ is contained in $\prod_{1 \leq i \leq j-1} \ker (\varphi \circ X_i)|_{(H_{i})_{x,r_i,\frac{r_i}{2}+}}(H_j)_{x,r_j+}$. Hence, adjusting $V_j'$ if necessary (to a subspace of $\pi(\prod_{1 \leq i \leq j-1} (H_{i})_{y,r_i-\frac{d_i}{2},\frac{r_i}{2}+})V_j'$), the action of 
$$(H_{i})_{y,r_i-\frac{d_i}{2},\frac{r_i}{2}+} /(H_{i})_{y,r_i+} \simeq (\fh_{i})_{y,r_i-\frac{d_i}{2},\frac{r_i}{2}+} /(\fh_{i})_{y,r_i+}$$ via $\pi$ on $V_j'$ is given by 
$$\varphi \circ \left(C_i+ X_i\right)$$
for some $C_i \in (\fh_i)^*_{y,-(\frac{r_i}{2}+\eps)}$ being trivial on $(\ft \cap \fh_i) + \fh_{i+1}$ for all $1 \leq i \leq  j-1$. Moreover, $X_j \in (\fg_j)^*_{x,-r_j} \subset (\fg_j)^*_{y,-r_j-\eps}$, and the action of $(H_j)_{y,r_j+\eps}$ on $V_{j-1}'$ factors through $(H_j)_{y,r_j+\eps}/(H_j)_{y,(r_j+\eps)+}\simeq (\fh_j)_{y,r_j+\eps}/(\fh_j)_{y,(r_j+\eps)+}$, on which it is given by $\varphi \circ X_j$ (which, as an aside, yields the trivial action).
By Corollary \ref{Cor-kill-C-type1}, there exists $g\in G_{y,\frac{r_{j-1}}{2}-\eps}\subset G_{y,0+}$ such that $\pi(g)V'_j \subset V_\pi^{\cup_{1 \leq i \leq j-1}((H_i)_{y,r_i+})}$, the action of $(H_{i})_{y,r_i,\frac{r_i}{2}+} /(H_{i})_{y,r_i+} \simeq (\fh_{i})_{y,r_i,\frac{r_i}{2}+} /(\fh_{i})_{y,r_i+}$ on $\pi(g)V_j'$ via $\pi$ is given by
 $\varphi \circ  X_i$ for $1 \leq i \leq j-1$,
 and the action of $(H_j)_{y,r_j+\eps}\subset (H_j)_{y,2\eps+}$ on $\pi(g)V_j'$ factors through $(H_j)_{y,r_j+\eps}/(H_j)_{y,(r_j+\eps)+}\simeq (\fh_j)_{y,r_j+\eps}/(\fh_j)_{y,(r_j+\eps)+}$ and is given by $\varphi \circ X_j$.
 Thus $f(y)\leq r_j+\eps$. Hence the set $f^{-1}(\bR_{\geq 0})$ is open in $\sB_j$.
 
 Moreover, if $f(x)=0$, then this implies that $f$ is continuous on $f^{-1}(\bR_{\geq 0})$, because $f(y) \geq 0$ and $r_j>0$ can be chosen arbitrarily small in this case.
 
 It remains to prove continuity around $x$ in the case $f(x)=r_j >0$. Suppose $f(y)< r_j-\eps$, and let $X'_j \in (\fg_j)^*_{y,-(r_j-\eps)+}$ be almost stable satisfying condition \eqref{item-in-pfb} above. Note that $(G_i)_{y,r_i-\frac{d_i}{2}}=(G_i)_{y,r_i}$ for $1 \leq i \leq j-1$. Hence, by the same reasoning as above (switching $x$ and $y$), we deduce that $f(x) < r_j$, a contradiction. 
 Thus $f(y) \geq r_j-\eps$ and $f$ is continuous on $f^{-1}(\bR_{\geq 0})$.
 
\textit{Proof of part (\ref{Lemma-Aj-iv}).}
Suppose $y \in \sB_j$ is in the closure of $f^{-1}(\bR_{\geq 0})$, and let $d>0$ be sufficiently small such that for all $r \in \bR_{\geq d}$ with  $G_{y,r} \neq G_{y,r+}$ we have $G_{y,{r-d}}=G_{y,r}$.
Let $\frac{d}{8} > \eps >0$ and $x \in f^{-1}(\bR_{\geq 0})$ with $\mathrm{d}(x,y)<\eps$. Then $G_{x,r} \neq G_{x,r+}$ implies $G_{x,{r-d+2\eps}}=G_{x,r}$ (if $r \in \bR_{\geq d-2\eps}$), hence $G_{x,{r-\frac{d}{2}}}=G_{x,r}$, and $G_{x,0+}=G_{x,\frac{d}{2}}$.  Thus we can apply the proof of part \eqref{Lemma-Aj-iii} to deduce that $f(y)$ is finite.

\qed \textsubscript{Lemma \ref{Lemma-Aj}}

Since $f$ is $G_j(k)$-equivariant, continuous, bounded below by zero, and the fundamental domain for the action of $G_j(k)$ on $\sB_j$ is bounded, there exists a point $x_j \in \sB_j$ such that $f(x_j) \leq f(x)$ for all $x \in \sB_j$. Define $r_j=f(x_j)$ and note that $r_j \leq f(x_{j-1}) \leq r_{j-1}$. 

We distinguish two cases.

\textbf{Case 1: $r_j>0$.}

 Let $(x_j, (X_i)_{1 \leq i \leq j-1})$ be a \truncateddatum/ contained in $(\pi,V_\pi)$, $X_j \in (\fg_j)^*_{x,r_j}$ \almoststable/ and $V_{j}' \subset V_{j-1} \subset V_\pi^{\cup_{1 \leq i \leq j-1}((H_i)_{x,r_i+})}$ satisfying the conditions \eqref{item-in-pfa} and \eqref{item-in-pfb} above. 

\begin{Lemmasub} \label{Lemma-Pf1} The element $X_j$ of $\fg_j^*$ is \almoststronglystable/ at $x_j \in \sB_j$.
\end{Lemmasub}
\textbf{Proof of Lemma \ref{Lemma-Pf1}.} Suppose $X_j$ is not \almoststronglystable/. Since $X_j$ is \almoststable/, this implies that $\ov{X_j} \in ((\fg_j)_{x_j,r_j}/(\fg_j)_{x_j,r_j+})^*$ is unstable. Thus, by \cite[Corollary~4.3]{Kempf} there exists a non-trivial one parameter subgroup $\ov \lambda: \bG_m \ra {(\RP_j)}_{x_j}$ in the reductive quotient ${(\RP_j)}_{x_j}$ of $G_j$ at $x_j$ such that $\lim_{t \ra 0}\ov\lambda(t).\ov{X_j}=0$. This means $\ov{X_j}$ is trivial on the root spaces corresponding to roots $\alpha$ with $\<\alpha,\lambda\><0$.  Let $\sS$ be a split torus of the parahoric group scheme ${(\bP_j)}_{x_j}$ of $G_j$ such that $\sS_{\ff}$ is a maximal split torus of $({\RP_j})_{x_j}$ containing $\ov \lambda(\bG_m)$ and such that $\sS_k$ is contained in a maximal torus $T_j \subset G_j$ which splits over a tame extension $E$ of $k$ and whose apartment $\sA(T_E,E) \cap \sB_j$ contains $x_j$. Let $\lambda: \bG_m \ra \sS_{k}$ be the one parameter subgroup corresponding to $\ov \lambda$. Then for $\eps>0$ small enough, we have $(H_j)_{x_j+\eps \lambda,r_j}\subset (H_j)_{x_j,r_j}$ and $X_j \in (\fg_j)_{x_j+\eps \lambda,-r_j}^*$, where $x_j+\eps \lambda \in \sA(T_E,E) \cap \sB_j$. Moreover, the image of $X_j$ in $(\fg_j)_{x_j+\eps \lambda,-r_j}^*/(\fg_j)_{x_j+\eps \lambda,-r_j+}^*$ is trivial. Let $r_j > \delta >0$ such that the subgroup $(H_j)_{x_j+\eps \lambda,r_j-\delta}$ equals $(H_j)_{x_j+\eps \lambda,r_j}$ and therefore acts trivially on $V_j'$. Analogously to the first part of the proof of Lemma \ref{Lemma-Aj}\eqref{Lemma-Aj-iii}, for $\eps$ sufficiently small,
there exist $d_i>0$ for $1 \leq i \leq j-1$ such that we have $V_j'\subset V_\pi^{\cup_{1 \leq i \leq j-1}((H_i)_{x_j+\eps \lambda,r_i+})}$ and (after potentially adjusting $V_j'\subset V_\pi^{\cup_{1 \leq i \leq j-1}((H_i)_{x_j+\eps \lambda,r_i+})}$) the action of $(H_{i})_{x_j+\eps \lambda,r_i-\frac{d_i}{2},\frac{r_i}{2}+} /(H_{i})_{x_j+\eps \lambda,r_i+} \simeq (\fh_{i})_{x_j+\eps \lambda,r_i-\frac{d_i}{2},\frac{r_i}{2}+} /(\fh_{i})_{x_j+\eps \lambda,r_i+}$ via $\pi$ on $V_j'$ is given by
$\varphi \circ \left(C_i + X_i\right)$
for some $C_i \in (\fh_i)^*_{x_j+\eps \lambda,-(\frac{r_i}{2}+\eps)}$ being trivial on $(\ft \cap \fh_i) + \fh_{i+1}$ for all $1 \leq i < j$, the group $(H_j)_{x_j+\eps\lambda,r_j-\delta}$ acts trivially on $V_j'$, and $(x_j+\eps\lambda, (X_i)_{1\leq i \leq j-1})$ is a \truncateddatum/. Assuming $\eps$ is sufficiently small and applying Corollary \ref{Cor-kill-C-type1} (or if $j=1$, set $g=1$), we obtain $g \in G_{x_j+\eps\lambda,0+}$ such that $\pi(g)V'_j \subset V_\pi^{\cup_{1 \leq i \leq j-1}((H_i)_{x_j+\eps \lambda,r_i+})}$, the action of $(H_{i})_{x_j+\eps \lambda,r_i,\frac{r_i}{2}+} /(H_{i})_{x_j+\eps \lambda,r_i+} \simeq (\fh_{i})_{x_j+\eps \lambda,r_i,\frac{r_i}{2}+} /(\fh_{i})_{x_j+\eps \lambda,r_i+}$ on $\pi(g)V_j'$ via $\pi$ is given by $\varphi \circ  X_i$ for $1 \leq i \leq j-1$, and the action of $(H_j)_{x_j+\eps \lambda,r_j-\delta}$ on $\pi(g)V_j'$ is trivial. Hence $f(x_j+\eps\lambda)\leq r_j-\delta < r_j=f(x_j)$, which contradicts the choice of $x_j$. Thus $X_j$ is \almoststronglystable/. \qed\textsubscript{Lemma \ref{Lemma-Pf1}}

Now we can show that, after changing $x_j$ and $X_j$ if necessary, we obtain a \truncateddatum/ of $G$ of length $j$ that is contained in $(\pi, V_\pi)$.

\begin{Lemmasub} \label{Lemma-Pf2} There exists a choice of $x_j$ and $X_j$ as above such that  $(x_j, (X_i)_{1 \leq i \leq j})$  is a  \truncateddatum/ contained in $(\pi, V_\pi)$.
\end{Lemmasub}
\textbf{Proof of Lemma \ref{Lemma-Pf2}.} 
Let $x_j$ and $X_j$ be as in Lemma \ref{Lemma-Pf1}. Let $\eps>0$ be sufficiently small (as specified later). By Proposition \ref{Lemma-almoststable-generic-representative} (applied to $G_j$) there exists $y \in \sB_j \subset \sB$ and $\wt X \in X_j + (\fg_j)_{y,(-r_j)+}^*$ such that $\mathrm{d}(x_j,y)< \eps$, the element $\wt X$ is \generic/ of depth $-r_j$ at $y$, and $x_j$ and $y$ are contained in $\sB(\Cent_{G_j}(\wt X),k)\subset \sB_j$. Note that for $\eps$ sufficiently small, we have $(H_j)_{y,r_j} \subset (H_j)_{x_j,r_j}$ and the action of $(H_j)_{y,r_j}$ on $V_j'$ factors through $(H_j)_{y,r_j}/(H_j)_{x_j,r_j+}$ on which it is given by $X_j$. Since $X_j-\wt X \in (\fh_j)_{y,(-r_j)+}^*$, this difference is trivial on $(\fh_j)_{y,r_j}$. Therefore the action of $(H_j)_{y,r_j}/(H_j)_{x_j,r_j+}$ on $V_j'$ is also given by $\wt X$, and, in particular, it factors through $(H_j)_{y,r_j}/(H_j)_{y,r_j+}$. Moreover, 
 the tuple $(y,(X_i)_{1 \leq i \leq j-1})$ is a \truncateddatum/ by Lemma \ref{Lemma-truncateddatum-with-y}. Substituting $X_j$ by $\wt X$ and applying Corollary \ref{Cor-kill-C-type1} (if $j>1$) as in the proofs of Lemma \ref{Lemma-Aj} and Lemma \ref{Lemma-Pf1} and possibly substituting $V_j'$ by $\pi(g)V_j'$ for some $g \in G_{y,0+}$, we can achieve that \eqref{item-in-pfa} and \eqref{item-in-pfb} above are satisfied at the point $y$. This implies that $f(y)=r_j$.

 Note that $(y, (X_i)_{1 \leq i \leq j})$ is a \truncateddatum/. (If $j>2$, then $G_j \neq \Cent_{G_j}(X_j)$, because otherwise $X_j(\fh_j)=0$ and hence $f(x_j)$ would not be minimal.) By Corollary \ref{Cor-kill-C-type2} this \truncateddatum/ is contained in $(\pi, V_\pi)$.  \qed\textsubscript{Lemma \ref{Lemma-Pf2}}

This finishes the recursion step. Since $G_j \subsetneq G_{j-1}$ for $j>2$, after repeating this construction finitely many times we obtain an integer $n$ and a \truncateddatum/ $(x_n, (X_i)_{1 \leq i \leq n})$ contained in $(\pi,V_\pi)$ with $r_{n+1}=0$, i.e. we move to the second case.

\textbf{Case 2: $r_j=0$.}
 Let $V_j'$ be the maximal subspace of $V_\pi^{\cup_{1 \leq i \leq j}((H_i)_{x_j,r_i+})}$ satisfying \eqref{item-in-pfa} and \eqref{item-in-pfb} above. Note that $(H_j)_{x_j,0}$ stabilizes $V_j'$, because $(H_j)_{x_j,0}$ centralizes $X_i$ and stabilizes $(H_i)_{x_j,r_i+}$ and $(H_i)_{x_j,r_i,\frac{r_i}{2}+}$ for $1 \leq i <j$. Let $(\repzero, \Vrepzero)$ be an irreducible  $(H_j)_{x_j,0}/(H_j)_{x_j,0+}$-subrepresentation of $V_j'$ viewed as a representation of $(H_j)_{x_j,0}/(H_j)_{x_j,0+}$. Then $(x_j, (X_i)_{1 \leq i \leq j-1},$ $ (\repzero,\Vrepzero))$ is a \datum/ contained in $(\pi, V_\pi)$. 
 \qed

 \begin{Rem} \label{Rmk-uniqueness}
 	In the next section we will use the existence of a maximal \datum/ for a given representation $(\pi, V_\pi)$ to deduce the existence of a \type/ for $(\pi, V_\pi)$. Note however that a \datum/ itself might not determine the Bernstein component, i.e. a given \datum/ might be a maximal \datum/ for representations in different Bernstein components. If one is interested in determining the Bernstein component uniquely, one has to enhance the \datum/ slightly (to a representation of $(M_{n+1})_x$, where $M_{n+1}$ is a Levi subgroup of $G_{n+1}$ that we are going to attached to $x$ and $G_{n+1}$ in Section \ref{Section-existence-of-type}, page \pageref{page-Levi}). Such an enhancement determines the Bernstein component uniquely by \cite[10.3~Theorem]{KimYu}, which is based on the work of Hakim--Murnaghan \cite{Hakim-Murnaghan} for supercuspidal representations. The assumption required in Hakim--Murnaghan was removed by Kaletha in \cite[Corollary~3.5.5]{Kaletha}.
 \end{Rem}

\section{From a {\datum/} to types} \label{Section-existence-of-type}
 Let $(\pi,V_\pi)$ be a smooth irreducible representation of $G(k)$, and let  $(x, (X_i)_{1 \leq i \leq n}, (\repzero, \Vrepzero))$ be a maximal \datum/ for $(\pi, V_\pi)$, which exists by Theorem \ref{Thm-existence-of-datum}.  In this section we show how to use this \datum/ in order to exhibit a \type/ contained in $(\pi, V_\pi)$. In order to do so we will define characters $\phi_i:G_{i+1}(k) \ra \bC$ of depth $r_i$ for $1 \leq i \leq n$ and a depth-zero representation of a compact open subgroup $K_{G_{n+1}}$ of $G_{n+1}(k)$ that contains $(G_{n+1})_{x,0}$. We will prove that these objects satisfy all necessary conditions imposed by Kim and Yu (\cite{KimYu}) so that Yu's construction (\cite{Yu}) yields a \type/. In Theorem \ref{Thm-exhaustion-of-types} we will conclude that the resulting \type/ is contained in $(\pi, V_\pi)$.

Recall that Moy and Prasad (\cite[6.3 and 6.4]{MP2}) attach to $x$ and $G_{n+1}$ a Levi subgroup $M_{n+1}$ \label{page-Levi} of $G_{n+1}$ such that $x \in \sB(M_{n+1},k) \subset \sB(G_{n+1},k)$ and $(M_{n+1})_{x,0}$ is a maximal parahoric subgroup of $M_{n+1}(k)$ with $(M_{n+1})_{x,0}/(M_{n+1})_{x,0+} \simeq (G_{n+1})_{x,0}/(G_{n+1})_{x,0+}$. We denote by $(M_{n+1})_{x}$ the stabilizer of $x \in \sB(M_{n+1},k)$ in $M_{n+1}(k)$. Then we define following Kim and Yu (\cite[7.1 and 7.3]{KimYu}) the group $K_{G_{n+1}}$ to be the group generated by $(M_{n+1})_x$ and $(G_{n+1})_{x,0}$. 

Let $V'$ be a subspace of $V_\pi$ as provided by Definition \ref{Def-datum-contained}  for the \datum/  $(x, (X_i)_{1 \leq i \leq n}, (\repzero, \Vrepzero))$ contained in $(\pi, V_\pi)$, and let $\wt V$ \label{page-V-tilda} be the irreducible $K_{G_{n+1}}$-subrepresentation of $V_\pi$ containing $V'$. 
Note that any  $g \in K_{G_{n+1}} \subset (G_{n+1})_x$ centralizes $X_i$ for $1 \leq i \leq n$ and hence stabilizes $(H_i)_{x,r_i+}$ ($1 \leq i \leq n+1$). Thus $\wt V$ is contained in $V_\pi^{\cup_{1 \leq i \leq n+1}((H_i)_{x,r_i+})}$.

 Moreover, let $T$ be a maximal torus of $M_{n+1} \subset G_{n+1}$ whose apartment contains $x$. Then, for $t \in T(k)_{0+}$ and $g \in K_{G_{n+1}}$, we have $tgt^{-1}g^{-1} \in (H_{n+1})_{x,0+}$. Hence, if $v \in \wt V$ is an element such that $T(k)_{0+}$ preserves $\bC \cdot v$, then $T(k)_{0+}$ also preserves $\bC \cdot gv$ and acts on both spaces via the same character. Since $\wt V$ is an irreducible $K_{G_{n+1}}$-representation, we deduce that $T(k)_{0+}$ acts on $\wt V$ via some character $\phi_T$ (times identity).

Before using $\phi_T$ to define the characters $\phi_i$, we recall Lemma~3.1.3 of \cite{Kaletha}.
\begin{Lemma}[\cite{Kaletha}] \label{Lemma-Kaletha}
	If $r \in \bR_{>0}$ and $1 \ra A \ra B \ra C \ra 1 $ is an exact sequence of tori that are defined over $k$ and split over a tamely ramified extension of $k$, then
	$$ 1 \ra A(k)_r \ra B(k)_r \ra C(k)_r \ra 1$$
	is an exact sequence.
\end{Lemma}

\begin{Cor} \label{Cor-generated}
	 Let $r \in \bR_{>0}$ and $1 \leq j \leq n+1$. Then
	$(G_{j})_{x,r}$ is generated by $T(k)_{r}$ and $(H_{j})_{x,r}$. 
\end{Cor}
\Proof Note that $T\cap H_{j}$ is a maximal torus of $H_{j}$ (\cite[Example~2.2.6]{ConradSGA3}). Then by Lemma \ref{Lemma-Kaletha} the map $T(k)_r \ra (T/T\cap H_{j})(k)_r = (G_{j}/H_{j})(k)_r$ is surjective, and hence $T(k)_r$ also surjects onto $G_{j}(k)_{x,r}/H_{j}(k)_{x,r} \subset (G_{j}/H_{j})(k)_r$. (That $G_{j}(k)_{x,r}$ maps to $(G_{j}/H_{j})(k)_r$ can be seen by considering a tame extension over which $T$ splits.) \qed 

Now we define $\phi_i$ recursively, first for $i=n$, then $i=n-1, n-2, \hdots, 1$. Suppose we have already defined $\phi_{n}, \hdots, \phi_{j+1}$ of depth $r_n, \hdots, r_{j+1}$ for some $1 \leq j \leq n$ ($j=n$ meaning no character has been defined yet) such that
 $$\pi|_{T(k) \cap (H_{j+1})_{x,0+}}  = \phi_n|_{T(k)\cap (H_{j+1})_{x,0+}} \cdot \hdots  \cdot \phi_{j+1}|_{T(k)\cap (H_{j+1})_{x,0+}}  \cdot \Id_{\wt V} \, \text{ on } \wt V .$$
Then we let $\phi_j'=\phi_T \cdot \phi_n|_{T(k)_{0+}}^{-1} \cdot \hdots \phi_{j+1}|_{T(k)_{0+}}^{-1}$, which is trivial on $T(k) \cap (H_{j+1})_{x,0+}=(T\cap H_{j+1})(k)_{0+}$ and on $T(k) \cap (H_{j})_{x,r_j+}=(T\cap H_{j})(k)_{r_j+}$. By Lemma \ref{Lemma-Kaletha} we have $(G_{j+1}/H_{j+1})(k)_{0+}\simeq T/(T\cap H_{j+1})(k)_{0+} \simeq T(k)_{0+}/(T\cap H_{j+1})(k)_{0+}$. Hence $\phi_j'$ defines a character of $(G_{j+1}/H_{j+1})(k)_{0+}$ that extends via Pontryagin duality to some character $\wt \phi_j'$ of $(G_{j+1}/H_{j+1})(k)$. Restricting $\wt \phi_j'$ to the image of $G_{j+1}(k)$, we obtain a character of $G_{j+1}(k)$ that we also denote by $\wt \phi_j'$. Note that $\wt \phi_j'$ is trivial on $H_{j+1}(k)$ and  $\wt \phi_j'|_{T(k)_{0+}}$ coincides with $\phi_j'$. Similarly, the character $\phi_j'$ gives rise to a character $\hat \phi_j'$  of $(G_{j}/H_{j})(k)_{r_j+}$ that can be extended and composed to yield a character (also denoted by $\hat \phi_j'$) of $G_j(k)$ that is trivial on $H_j(k)$ and coincides with $\wt \phi_j'$ on $T(k)_{r_j+}$. If $j=1$, we may and do choose $\hat \phi_j'$ to be the trivial character. We define $\phi_j=\wt \phi_j' \cdot (\hat \phi_j')^{-1}|_{G_{j+1}(k)}$. Then $\phi_j$ has depth $r_j$ (because by considering a tame extension that splits the torus $T$, we see that $G_{j+1}(k)_{x,r_j+}$ maps to $(G_{j+1}/H_{j+1})(k)_{r_j+}\simeq T(k)_{r_j+}/(T\cap H_{j+1})(k)_{r_j+}$; or use Corollary \ref{Cor-generated}) and
$$\pi|_{T(k) \cap (H_{j})_{x,0+}}  = \phi_n|_{T(k)\cap (H_{j})_{x,0+}} \cdot \hdots  \cdot \phi_{j}|_{T(k)\cap (H_{j})_{x,0+}}  \cdot \Id_{\wt V} \, \text{ on } \wt V .$$

\begin{Lemma} \label{Lemma-phi-properties}
For $1 \leq j \leq n$ the character $\phi_j:G_{j+1}(k) \ra \bC^*$ satisfies the following properties:
  \begin{enumerate}[label=(\roman*),ref=\roman*]
  	\item \label{item-character-extension-0}  $\phi_j$ is trivial on ${(G_{j+1})_{x,r_j+}}$ and on $H_{j+1}(k)$, 
  	\item \label{item-character-extension-1} $\phi_j|_{(H_j)_{x,r_j}\cap G_{j+1}(k)}$ factors through 
  	$$((H_j)_{x,r_j}\cap G_{j+1}(k))/((H_j)_{x,r_j+}\cap G_{j+1}(k)) \simeq ((\fh_j)_{x,r_j}\cap \fg_{j+1})/((\fh_j)_{x,r_j+}\cap \fg_{j+1})$$ and is given by $\varphi \circ X_j|_{(\fh_j)_{x,r_j}\cap \fg_{j+1}}$, 
  	\item \label{item-character-extension-generic} $\phi_j$ is $G_j$-generic of depth $r_j$ (in the sense of \cite[\S~9]{Yu}) relative to $x$, and
  	\item \label{item-character-extension-2} the group $(G_{n+1})_{x,0+}$ acts on $\wt V$ via $\prod_{1 \leq i \leq n} \phi_i|_{(G_{n+1})_{x,0+}}$.
  \end{enumerate}
\end{Lemma} 
\Proof   Part \eqref{item-character-extension-0} follows immediately from the above construction.

For Part \eqref{item-character-extension-1}, note that using Corollary \ref{Cor-generated} we see that $(H_j)_{x,r_j} \cap G_{j+1}(k)=H_j(k) \cap (G_{j+1})_{x,r_j}$ is generated by $T(k)_{r_j} \cap H_j(k)$ and $(H_{j+1})_{x,r_j}$. Since $\phi_j$ is trivial on $H_{j+1}(k)$ and coincides with $\phi_T$ on $T(k)_{r_j} \cap H_j(k)$, the claim follows from the properties of $V'$ in Definition \ref{Def-datum-contained}.

For Part \eqref{item-character-extension-generic}, note that by Part \eqref{item-character-extension-1} there exists $Y \in (\fg_{j+1})^*_{x,-r_j}$ such that $Y$ is trivial on $\fh_{j}$ and $\phi_j|_{(G_{j+1})_{x,r_j}}$ is given by the character of $(G_{j+1})_{x,r_j}/(G_{j+1})_{x,r_j+} \simeq (\fg_{j+1})_{x,r_j}/(\fg_{j+1})_{x,r_j+}$ arising from $\varphi \circ (Y+X_j)$. Since $Y$ is trivial on $\fh_j$, the element $Y$ is fixed under the dual of the adjoint action of $G_{j+1}$ on $\fg_{j+1}$. Hence by the definition of $G_{j+1}$, the group $G_{j+1}$ centralizes $Y+X_j$. Moreover, if $T$ is a maximal torus of $G_{j+1}$, and $\alpha \in \Phi(G_{j},T_{\ov k})-\Phi(G_{j+1},T_{\ov k})$, then $H_\alpha \in (\fh_j)_{\ov k}$ and hence 
$$ \val((Y+X_j)(H_\alpha))= \val(X_j(H_\alpha)) = -r_j,$$
where the last equality follows from Lemma \ref{Lemma-structure-of-X}. Since $p$ is not a torsion prime for the dual root datum of $G_j$ by Lemma \ref{Lemma-restriction-on-p}\eqref{item-p-Levi}, \eqref{item-torsion-prime} and \eqref{item-index-of-connection} (applied to the dual root datum of $G_j$), the character $\phi_j$ is $G_j$-generic of depth $r_j$ by \cite[Lemma~8.1]{Yu}.
 
  Part \eqref{item-character-extension-2} follows from the observation above that $\phi_1|_{T(k)_{0+}}=\phi_T\cdot \phi_n|_{T(k)_{0+}}^{-1} \cdot \hdots \phi_{2}|_{T(k)_{0+}}^{-1}$ and that $(H_{n+1})_{x,0+}$ acts trivially on $\wt V$ together with Corollary \ref{Cor-generated}. \qed

\begin{Cor} \label{Cor-repzero}
	The irreducible representation  $( \prod_{1 \leq i \leq n} \phi_i^{-1}|_{K_{G_{n+1}}}  \cdot \pi|_{K_{G_{n+1}}}, \wt V )$ of $K_{G_{n+1}}$ is trivial on $(G_{n+1})_{x,0+}$ and its restriction to $(H_{n+1})_{x,0}$ contains $(\repzero,\Vrepzero)$ as an irreducible subrepresentation. 
\end{Cor}

\Proof 
This is an immediate consequence of Lemma \ref{Lemma-phi-properties} \eqref{item-character-extension-0} and \eqref{item-character-extension-2} and the definition of $\wt V$. 
\qed

Recall our convention that by ``\type/'' we mean an \stype/ for some inertial equivalence class $\fs \in \fI$.
In order to obtain a type for our representation $(\pi,V_\pi)$ of $G(k)$ using the construction of Kim and Yu in \cite{KimYu} we denote by $r_\pi$ the depth of the representation $(\pi, V_\pi)$, i.e. $r_\pi=r_1$ if $n \geq 1$ and $r_\pi=0$ if $n=0$, and we make the following definitions: \label{page-Yu-datum}
\begin{eqnarray*} 
	\vec G &=& \left\{ \begin{array}{ll} (G_{n+1}, G_n, \hdots, G_2, G_1=G) & \text{ if } G_2 \neq G_1 \text{ or } n=0 \\
(G_{n+1}, G_n, \hdots, G_3, G_2=G) & \text{ if } G_2=G_1 \end{array}\right. 
\\
 \vec r & = & \left\{ \begin{array}{ll} (r_{n}, r_{n-1}, \hdots, r_2, r_1, r_\pi) & \text{ if } G_2 \neq G_1 \text{ or } n=0 \\
 	 (r_{n}, r_{n-1}, \hdots, r_2, r_1)  & \text{ if } G_2=G_1 \end{array}\right.
 	\\
 \vec \phi & = & \left\{ \begin{array}{ll} (\phi_{n}, \phi_{n-1}, \hdots, \phi_2, \phi_1, 1) & \text{ if } G_2 \neq G_1 \text{ or } n=0 \\
 	(\phi_{n}, \phi_{n-1}, \hdots, \phi_2, \phi_1)  & \text{ if } G_2=G_1 \end{array}\right. 	
 \\
 K &=& K_{G_{n+1}}(G_{n})_{x,\frac{r_n}{2}}\hdots(G_1)_{x,\frac{r_1}{2}} 
  \\
   K_{0+} &=& (G_{n+1})_{x,0+}(G_{n})_{x,\frac{r_n}{2}}\hdots(G_1)_{x,\frac{r_1}{2}}
  \\
  K_+ &=& (G_{n+1})_{x,0+}(G_{n})_{x,\frac{r_n}{2}+}\hdots(G_1)_{x,\frac{r_1}{2}+} \\
  K^H_{0+} &=& (H_{n+1})_{x,0+}(H_{n})_{x,\frac{r_n}{2}}\hdots(H_1)_{x,\frac{r_1}{2}}
     \\
     K^H_+ &=& (H_{n+1})_{x,0+}(H_{n})_{x,\frac{r_n}{2}+}\hdots(H_1)_{x,\frac{r_1}{2}+} 
\end{eqnarray*}

\begin{Lemma} \label{Lemma-Ks}
We have the following identities: 
	 \begin{equation*} 
	 	\begin{aligned}	K_+ & =   (G_{n+1})_{x,0+}(H_{n})_{x,\frac{r_n}{2}+}\hdots(H_1)_{x,\frac{r_1}{2}+}=(G_{n+1})_{x,0+}K^H_+ \\
	 	 \quad    K_{0+} & = (G_{n+1})_{x,0+}(H_{n})_{x,\frac{r_n}{2}}\hdots(H_1)_{x,\frac{r_1}{2}}= (G_{n+1})_{x,0+} K^H_{0+}   \\ 
	   K^H_{0+} &= (H_{n+1})_{x,0+}(H_{n})_{x,r_n,\frac{r_n}{2}}\hdots(H_1)_{x,r_1,\frac{r_1}{2}}	\\
	   K^H_+ &= (H_{n+1})_{x,0+}(H_{n})_{x,r_n,\frac{r_n}{2}+}\hdots(H_1)_{x,r_1,\frac{r_1}{2}+} \, . 
	 	\end{aligned}
	 \end{equation*}
\end{Lemma}
\Proof
The first two lines follow from Corollary \ref{Cor-generated}.
It is clear that 
$$	K^H_{0+} \supset (H_{n+1})_{x,0+}(H_{n})_{x,r_n,\frac{r_n}{2}}\hdots(H_1)_{x,r_1,\frac{r_1}{2}} .$$
 In order to prove the fourth identity, it remains to show that  
 $$(H_{i})_{x, \frac{r_i}{2}+} \subset (H_{n+1})_{x,0+}(H_{n})_{x,r_n,\frac{r_n}{2}+}\hdots(H_1)_{x,r_1,\frac{r_1}{2}+}$$
  for all $n+1 \geq i \geq 1$, where we recall that $r_{n+1}=0$. We show this by induction. For $i=n+1$ the statement is obvious, so assume $n \geq i \geq 1$ and that the statement holds for $i+1$. Then $(H_{i+1})_{x, \frac{r_i}{2}+} \subset (H_{i+1})_{x, \frac{r_{i+1}}{2}+} \subset (H_{n+1})_{x,0+}(H_{n})_{x,r_n,\frac{r_n}{2}+}\hdots(H_1)_{x,r_1,\frac{r_1}{2}+}$, and it suffices to prove that $(H_{i})_{x, \frac{r_i}{2}+} = (H_{i})_{x, r_i, \frac{r_i}{2}+}  (H_{i+1})_{x, \frac{r_i}{2}+}$, or, equivalently,
$$(\fh_{i})_{x, \frac{r_i}{2}+} / (\fh_{i})_{x, r_i} = ((\fh_{i})_{x, r_i, \frac{r_i}{2}+} + (\fh_{i+1})_{x, \frac{r_i}{2}+})/ (\fh_{i})_{x, r_i} . $$
This follows by taking $\Gal(E/k)$-invariants of the following equality (of abelian groups with $\Gal(E/k)$-action)
$$ (\fh_{i}(E))_{x, \frac{r_i}{2}+} / (\fh_{i}(E))_{x, r_i} = (\fh_{i}(E))_{x, r_i, \frac{r_i}{2}+}/ (\fh_{i}(E))_{x, r_i} \oplus (\fh_{i+1}(E))_{x, \frac{r_i}{2}+}/(\fh_{i+1}(E))_{x, r_i}, $$ 
where $E$ is a tamely ramified extension of $k$ over which $G_{i+1}$ and $G_i$ split. The third identity is proved analogously. \qed

Let $\rhoYu$ be an irreducible representation of $K_{G_{n+1}}$ such that $\rhoYu|_{(G_{n+1})_{x,0}}$ factors through $(G_{n+1})_{x,0}/(G_{n+1})_{x,0+}$ and contains a cuspidal representation of $(G_{n+1})_{x,0}/(G_{n+1})_{x,0+}$ .   By Lemma \ref{Lemma-phi-properties}\eqref{item-character-extension-generic}
  and \cite[7.3~Remark]{KimYu}\footnote{Remark 7.3. in \cite{KimYu} explains how to get the 5-tuple $\Sigma$ (using the notation from \cite{KimYu}) from our 5-tuple. The authors mention in this remark that as a last step one ``can then extend/modify $\iota$ to a family $\{\iota\}$ which is $\vec s$-generic''. However, by doing so one might have to change our point $x$ (which is denoted by $y$ in \cite{KimYu}) to a nearby point in the building. In order to keep working with $x$ we will not perform this last modification. As a consequence the  requirement of $\{\iota\}$ being $\vec s$-generic in Condition D2 of \cite[7.2]{KimYu} might not be satisfied. However, we can still carry out Yu's construction with our tuple.} the tuple $(\vec G, x, \vec r, \rhoYu, \vec \phi)$ satisfies Conditions D1, D3, D4 and D5\footnote{Ju-Lee Kim confirmed that ''relative to $x$ for all $x \in \sB(G')$'' in Condition D5 in \cite[7.2]{KimYu} should be ``relative to $y$'' (using the notation of \cite{KimYu}). 
  } in \cite[7.2]{KimYu}.
 Using this tuple we can carry out Yu's construction 
  (\cite[\S4]{Yu}) as explained in \cite[7.4]{KimYu} to obtain a representation of $K$ that we denote by $(\piK,\VpiK)$.
  (Note that $p \neq 2$ by our assumption that $p \nmid \abs{W}$ and that $G$ is not a torus.)
  By construction, the representation $\piK$ is of the form $\rhoYu \otimes \kappa_{\vec \phi}$, where $\rhoYu$ also denotes the extension of $\rhoYu$ from $K_{G_{n+1}}$ to $K$ that is trivial on $(G_n)_{x,\frac{r_n}{2}}\hdots (G_1)_{x,\frac{r_1}{2}}$, and $(\kappa_{\vec \phi}, V_\kappa)$ \label{page-kappa} is a representation of  $K$ that depends only on $(\vec G, \vec r, \vec \phi)$, i.e. not on the choice of $\rhoYu$ (\cite{Yu}[\S4] or \cite[12.4]{Kim}). In particular, $(\piK|_{K_+},\VpiK)$ does not depend on $\rhoYu$. 

We denote by $\hat \phi_i$ ($1 \leq i \leq n$) the character of $K_{G_{n+1}}(G_{i+1})_{x,0}G_{x,\frac{r_i}{2}+}$ defined in \cite[\S~4]{Yu}, i.e. the unique character of $K_{G_{n+1}}(G_{i+1})_{x,0}G_{x,\frac{r_i}{2}+}$ satisfying
\begin{itemize}
	\item $\hat \phi_i|_{K_{G_{n+1}}(G_{i+1})_{x,0}}=\phi_i|_{K_{G_{n+1}}(G_{i+1})_{x,0}}$, and 
	\item $\hat \phi_i|_{G_{x,\frac{r_i}{2}+}} $ factors through 
	\begin{eqnarray*} 
		G_{x,\frac{r_i}{2}+}/G_{x,r_i+} & \simeq& \fg_{x,\frac{r_i}{2}+}/\fg_{x,r_i+} = (\fg_{i+1} \oplus \fr'')_{x,\frac{r_i}{2}+}/(\fg_{i+1} \oplus \fr'')_{x,r_i+} \\
		& \ra & (\fg_{i+1})_{x,\frac{r_i}{2}+}/(\fg_{i+1})_{x,r_i+} \simeq 
		(G_{i+1})_{x,\frac{r_i}{2}+}/(G_{i+1})_{x,r_i+},
	\end{eqnarray*}
	on which it is induced by $\phi_i$. Here $\fr''$ is as defined in Lemma \ref{Lemma-kill-C}, i.e. $\fr''=\fg \cap \bigoplus_{\alpha \in \Phi(G,T_E)-\Phi(G_{i+1},T_E)} (\fg_E)_\alpha$ for some maximal torus $T$ of $G_{i+1}$ that splits over a tame extension $E$ of $k$, and the map $\fg_{i+1} \oplus \fr'' \ra \fg_{i+1}$ sends $\fr''$ to zero.
\end{itemize}
Then Yu  proves in \cite[Proposition~11.4]{Yu} that  $(G_i)_{x,r_i,\frac{r_i}{2}}/\left((G_i)_{x,r_i,\frac{r_i}{2}+} \cap \ker(\hat \phi_i)\right)$ is a Heisenberg $p$-group with center $(G_i)_{x,r_i,\frac{r_i}{2}+} / \left((G_i)_{x,r_i,\frac{r_i}{2}+} \cap \ker(\hat \phi_i)\right)$. Let $(\omega_i, V_{\omega_i})$ denote the Heisenberg representation of this Heisenberg $p$-group with central character $\hat \phi_i|_{(G_i)_{x,r_i,\frac{r_i}{2}+} }$. Then we observe from the construction of $(\kappa_{\vec \phi}, V_\kappa)$ and \cite[Theorem~11.5]{Yu} that $(\kappa_{\vec \phi}|_{K_{0+}}, V_\kappa)$ is irreducible and that the underlying vector space $V_\kappa$ is  $\bigotimes_{i=1}^n V_{\omega_i}$. If $n=0$, then the empty tensor product is meant to be a one dimensional vector space. In that case $(\kappa_{\vec \phi}|_{K_{0+}}, V_\kappa)$ is the trivial one dimensional representation. 
The restriction of $(\kappa_{\vec \phi}, V_\kappa)$ to $(G_i)_{x,r_i,\frac{r_i}{2}}$ for $1  \leq i \leq n$ is given by letting $(G_i)_{x,r_i,\frac{r_i}{2}}$ act via the Heisenberg representation $\omega_i$ of $(G_i)_{x,r_i,\frac{r_i}{2}}/\left((G_i)_{x,r_i,\frac{r_i}{2}+} \cap \ker(\hat \phi_i)\right)$ with central character $\hat \phi_i|_{(G_i)_{x,r_i,\frac{r_i}{2}+} }$  on $V_{\omega_i}$ and via $\hat \phi_j|_{(G_i)_{x,r_i,\frac{r_i}{2}}}$ on $V_{\omega_j}$ for $j \neq i$.

\begin{Lemma} \label{Lemma-K+}
There exists an irreducible $K_+$-subrepresentation 
	 of $(\pi|_{K_+}, \wt V)$ that is isomorphic to any one-dimensional $K_+$-subrepresentation of $(\piK|_{K_+},\VpiK)$.
\end{Lemma}
\Proof  By \cite[Proposition~4.4]{Yu}, the representation $(\piK|_{K_+},\VpiK)$ is $\theta:=\prod_{1 \leq i \leq n} \hat \phi_i|_{K_+}$-isotypic.
Let $(\pi|_{K_+},V'')$ be an irreducible $K_+$-subrepresentation of $(\pi|_{K_+}, V') \subset (\pi|_{K_+}, \wt V)$. By Lemma \ref{Lemma-phi-properties}\eqref{item-character-extension-2}, the group $(G_{n+1})_{x,0+}$ acts on $V''$ via $\theta$. Moreover, by Lemma \ref{Lemma-phi-properties}\eqref{item-character-extension-0} and \eqref{item-character-extension-1} the restriction $\theta|_{(H_{i})_{x,r_i,\frac{r_i}{2}+}}$ for $1 \leq i \leq n$ factors through $(H_{i})_{x,r_i,\frac{r_i}{2}+} /(H_{i})_{x,r_i+} \simeq (\fh_{i})_{x,r_i,\frac{r_i}{2}+} /(\fh_{i})_{x,r_i+}$, where it is given by $\varphi \circ X_i$ (by the last line of Lemma \ref{Lemma-structure-of-X}). Hence the group $(H_{i})_{x,r_i,\frac{r_i}{2}+}$ acts on $V''$ via $\theta$ for $1 \leq i \leq n$. Since $(G_{n+1})_{x,0+}$ together with $(H_{i})_{x,r_i,\frac{r_i}{2}+}$, $1 \leq i \leq n$, generate $K_+$ by Lemma \ref{Lemma-Ks}, we are done. \qed

We denote by $N^H$ the kernel in $K^H_+$ of $\theta|_{K^H_+}=\prod_{1 \leq i \leq n} \hat \phi_i|_{K^H_+}$.
\begin{Lemma} \label{Lemma-Heisenberg}
  If $n>0$, 
   then $K^H_{0+}/N^H$ is a Heisenberg $p$-group with center $K_+^H/N^H$.\\
  If $n=0$,
   then $K^H_{0+}/N^H=K^H_{+}/N^H$.
\end{Lemma}
\Proof
Note that $[K^H_{0+}, K^H_{0+}] \subset K^H_+$ and $[K^H_{0+}, K^H_{+}] \subset (H_{n+1})_{x,0+}(H_{n})_{x,r_n+,\frac{r_n}{2}+}\hdots(H_1)_{x,r_1+,\frac{r_1}{2}+} \subset N^H$. Thus the center of $K^H_{0+}/N^H$ contains $K_+^H/N^H$ and we have a pairing $(a,b)= \theta(aba^{-1}b^{-1})$ on $K^H_{0+}/K^H_{+} \times K^H_{0+}/K^H_{+}$. 
Note that 
$$K^H_{0+}/K^H_{+} \simeq (H_1)_{x,r_1,\frac{r_1}{2}}/(H_1)_{x,r_1,\frac{r_1}{2}+} \oplus \hdots \oplus (H_n)_{x,r_n,\frac{r_n}{2}}/(H_n)_{x,r_n,\frac{r_n}{2}+} $$
and it is easy to check (as done in the proof of \cite[Proposition~18.1]{Kim}) that 
$( \cdot, \cdot)$ is the sum of the pairings $(\cdot, \cdot)_i$ on $(H_i)_{x,r_i,\frac{r_i}{2}}/(H_i)_{x,r_i,\frac{r_i}{2}+} $ defined by 
$(a,b)_i=\hat \phi_i(aba^{-1}b^{-1})$. By \cite[Lemma~11.1]{Yu}, the pairing $(\cdot,\cdot)_i$ is non-degenerate $1 \leq i \leq n$, and hence the pairing $(\cdot,\cdot)$ is non-degenerate. Thus the center of $K^H_{0+}/N^H$ is contained in $K_+^H/N^H$, and therefore equals $K_+^H/N^H$. Moreover, the image of $\theta|_{K^H_+}$ is $\{ c \in \bC \, | \, c^p=1\}$, which implies that $K_+^H/N^H$ has order $p$. The remainder of the proof works completely analogous to Yu's proof (\cite[Proposition~11.4]{Yu}) that the group $(H_i)_{x,r_i,\frac{r_i}{2}}/\left((H_i)_{x,r_i,\frac{r_i}{2}+} \cap \ker(\hat \phi_i)\right) $ is a Heisenberg $p$-group with center $(H_i)_{x,r_i,\frac{r_i}{2}+}/\left((H_i)_{x,r_i,\frac{r_i}{2}+} \cap \ker(\hat \phi_i)\right)$ for $1 \leq i \leq n$. We outline the proof as a convenience for the reader and refer to \cite[Proposition~11.4]{Yu} for details: We first prove the statement over a tame extension $E$ over which $G_{n+1}$ is split, and denote by $K^H_{0+}(E), K^H_0(E)$ and $N(E)$ the corresponding groups constructed over $E$. By \cite[Lemma~10.1]{Yu} and the above observations (over $E$), it suffices to exhibit subgroups $W_1$
and $W_2$ of $K^H_{0+}(E)/N^H(E)$ that have trivial intersection with the center and whose image in $K^H_{0+}(E)/K^H_0(E)$ form a complete polarization. This can be achieved by using positive and negative root groups, respectively. To conclude that $K^H_{0+}/N^H$ is a Heisenberg $p$-group, we then embed $K^H_{0+}/N^H$ into $K^H_{0+}(E)/N^H(E)$, observe that by above its image $K^H_{0+}/K^H_0$ in $K^H_{0+}(E)/K^H_0(E)$ is a non-degenerate subspace, and apply \cite[Lemma~10.3]{Yu}.

The second half of the lemma follows immediately from the definition of $K^H_{0+}$ and $K^H_+$. 
\qed

Let $(\pi|_{K}, \widehat V)$ \label{page-V-hat} be the irreducible $K$-subrepresentation of $(\pi, V_\pi)$ that contains $\wt V$. 

\begin{Lemma} \label{Lemma-repzeroK}
	There exists an irreducible representation $(\repzeroK,\VrepzeroK)$ of $K$ that is trivial on $K_{0+}$ such that $( \repzeroK \otimes \kappa_{\vec \phi} ,  \VrepzeroK \otimes V_\kappa) \simeq (\pi|_{K}, \widehat V)$.
\end{Lemma}
\Proof
Since $K_{0+}=G_{x,0+}K^H_{0+}$ (Lemma \ref{Lemma-Ks}) and $G_{x,0+}\subset K_+$ acts on $\VpiK$ via $\theta|_{G_{x,0+}}$ (times identity) by \cite[Proposition~4.4]{Yu}, we deduce from the irreducibility of   $(\kappa_{\vec \phi}|_{K_{0+}}, V_\kappa)$ mentioned above that also its restriction $(\kappa_{\vec \phi}|_{K^H_{0+}}, V_\kappa)$ to $K^H_{0+}$ is irreducible. Recall that $(\kappa_{\vec \phi}|_{K^H_{0+}}, V_\kappa)$ factors through $K^H_{0+}/N^H$ and $K^H_+$ acts via the character $\theta|_{K^H_+}$ (times identity). 
By Lemma \ref{Lemma-Heisenberg} and the theory of Heisenberg representations there exists a unique irreducible representation of $K^H_{0+}$ factoring through $K^H_{0+}/N^H$ and having $K^H_+/N^H$ act via the character $\theta|_{K^H_+}$ (times identity). On the other hand, Lemma \ref{Lemma-K+} and the observation that $[K^H_{0+},K_{+}]\subset N^H$ imply that $(\pi|_{K^H_{0+}}, \widehat V)$ contains an irreducible $K^H_{0+}$-subrepresentation on which $K_{+}$ acts via the character $\theta|_{K_+}$ (times identity), and which therefore is isomorphic to  $(\kappa_{\vec \phi}|_{K^H_{0+}}, V_\kappa)$ as a $K^H_{0+}$-representation. Moreover, since $K_{0+}=K_+K^H_{0+}$, we deduce from the $K_+$-action that $(\pi|_{K_{0+}}, \widehat V)$ contains an irreducible $K_{0+}$-subrepresentation isomorphic to  $(\kappa_{\vec \phi}|_{K_{0+}}, V_\kappa)$. Hence, by \cite[Proposition~18.5]{Kim} (or rather the analogous statement in our setting that is proved in the same way), the irreducible representation $(\pi|_K,\widehat V)$ of $K$ that extends  $(\kappa_{\vec \phi}|_{K_{0+}}, V_\kappa)$  is of the form $(\repzeroK \otimes \kappa_{\vec \phi},  \VrepzeroK \otimes V_\kappa)$ for some irreducible representation $(\repzeroK,\VrepzeroK)$ of $K$ that is trivial on $K_{0+}$. \qed

\begin{Cor} \label{Cor-action-of-H}
	The subspace $\widehat V$ is contained in $V_\pi^{\cup_{1 \leq i \leq n+1}((H_i)_{x,r_i+})}$ and the action of the group $(H_i)_{x,r_i,\frac{r_i}{2}+}/(H_i)_{x,r_i+}\simeq (\fh_i)_{x,r_i,\frac{r_i}{2}+}/(\fh_i)_{x,r_i+}$ on $\widehat V$ via $\pi$ is given by the character $\varphi \circ X_i$ for $1 \leq i \leq n$.
\end{Cor}
\Proof
Let $1 \leq i \leq n$. We have $(H_i)_{x,r_i,\frac{r_i}{2}+} \subset K_+$ and $(H_{n+1})_{x,0+} \subset K_+$ and by Lemma \ref{Lemma-repzeroK} the representation $(\pi|_{K_+}, \widehat V)$  is $\theta$-isotypic. As we saw in the proof of Lemma \ref{Lemma-K+}, the character $\theta|_{(H_i)_{x,r_i,\frac{r_i}{2}+}}$ factors through $(H_i)_{x,r_i,\frac{r_i}{2}+}/(H_i)_{x,r_i+}\simeq (\fh_i)_{x,r_i,\frac{r_i}{2}+}/(\fh_i)_{x,r_i+}$, on which it is given by $\varphi \circ X_i$, and $\theta|_{(H_{n+1})_{x,0+}}$ is trivial by Lemma \ref{Lemma-phi-properties}\eqref{item-character-extension-0}. \qed

\begin{Lemma} \label{Lemma-cuspidal}
	The irreducible components of the representation $(\repzeroK|_{(G_{n+1})_{x,0}}, \VrepzeroK)$ provided by Lemma \ref{Lemma-repzeroK} are cuspidal representations of $(G_{n+1})_{x,0}/(G_{n+1})_{x,0+}$. 
\end{Lemma}

\begin{Rem} 
	Readers familiar with Kim's work may expect that we could mainly cite \cite{Kim} for the proof of Lemma \ref{Lemma-cuspidal}. However, contrary to the claim in \cite[Proposition~17.2.(2)]{Kim}, the representation $\rho|_{(G_{n+1})_{x,0}} \otimes \kappa_{\vec \phi}|_{(G_{n+1})_{x,0}} \otimes \prod_{1 \leq i \leq n} \phi_i^{-1}|_{(G_{n+1})_{x,0}}$ might not necessarily be cuspidal when viewed as a representation of $(G_{n+1})_{x,0}/(G_{n+1})_{x,0+}$. Since the proof of the above mentioned proposition in \cite{Kim} is not correct, we provide a different and independent proof of Lemma \ref{Lemma-cuspidal}.

\end{Rem}	

\textbf{Proof of Lemma \ref{Lemma-cuspidal}}
Suppose $(\rho', V_{\rho'})$ is an irreducible subrepresentation of $(\repzeroK|_{(G_{n+1})_{x,0}}, \VrepzeroK)$ that is not cuspidal (viewed as a representation of $(G_{n+1})_{x,0}/(G_{n+1})_{x,0+}$). 
Then there exists (the $\ff$-points of) a unipotent radical $\Uff$ of a (proper) parabolic subgroup of the reductive group (with $\ff$-points) $(G_{n+1})_{x,0}/(G_{n+1})_{x,0+}=(\RP_{n+1})_x(\ff)$ such that $\rho'|_{\Uff}$ contains the trivial representation of $\Uff$. Denote by $V_{\rho'}''$ a subspace of $V_{\rho'}$ on which $\Uff$ acts trivially.
 By \cite[Corollary~2.2.5 and Proposition~2.2.9]{Pseudoreductive2} there exists a one parameter subgroup $\ov \lambda: \bG_m \ra (\RP_{n+1})_x$ such that $\Uff=\{g \in (\RP_{n+1})_x(\ff) \, | \lim_{t \ra 0}\ov\lambda(t).g=1 \,\}$. 
  Let $\lambda: \bG_m \ra G_{n+1}$ denote a lift of $\ov \lambda$ that factors through a maximally split maximal torus $T$ of $G_{n+1}$ whose apartment $\sA(T)$ contains $x$ (see the proof of Lemma \ref{Lemma-Pf1} for more details about such a lift). 
Let $\kappa'={\kappa_{\vec \phi}}|_{(G_{n+1})_{x,0}} \otimes \prod_{1 \leq i \leq n} \phi_i^{-1}|_{(G_{n+1})_{x,0}}$, which is trivial on $(G_{n+1})_{x,0+}$ (by either combining Corollary \ref{Cor-repzero} with Lemma \ref{Lemma-K+} or by using the proof of Lemma \ref{Lemma-K+}) and therefore can also be regarded as a representation of $(G_{n+1})_{x,0}/(G_{n+1})_{x,0+}$ and hence of $\Uff$.
Recall that $(\omega_i, V_{\omega_i})$ denotes the Heisenberg representation of $(G_i)_{x,r_i,\frac{r_i}{2}}/\left((G_i)_{x,r_i,\frac{r_i}{2}+} \cap \ker(\hat \phi_i)\right)$ with central character (the restriction of) $\hat \phi_i$, and that the vector space $V_\kappa$ underlying the representation of $\kappa'$ is $\bigotimes_{i=1}^n V_{\omega_i}$.
By the construction of Yu (\cite{Yu}[\S4, p.~592 and Theorem~11.5]), the representation $\kappa'$ is defined by letting $(G_{n+1})_{x,0}/(G_{n+1})_{x,0+}$ act on each of the tensor product factors $V_{\omega_i}$ in $\bigotimes_{i=1}^nV_{\omega_i}$ by mapping $(G_{n+1})_{x,0}/(G_{n+1})_{x,0+}$ to the symplectic group  $\Sp(V_i)$ of the corresponding symplectic space $V_i:=(G_i)_{x,r_i,\frac{r_i}{2}}/(G_i)_{x,r_i,\frac{r_i}{2}+}$ with pairing defined by 
$(a,b)_i=\hat \phi_i(aba^{-1}b^{-1})$ and composing with a Weil representation. 
The map from $(G_{n+1})_{x,0}/(G_{n+1})_{x,0+}$ to $\Sp(V_i)$ is induced by the conjugation action of $(G_{n+1})_{x,0}$ on $(G_i)_{x,r_i,\frac{r_i}{2}}$. 

Let $E$ be a tamely ramified extension of $k$ over which $T$ splits, and define for $1 \leq i \leq n$ the space
$V_i^+$ to be the image 
of $ G(k) \cap \<U_\alpha(E)_{x,\frac{r_i}{2}} \, | \, \alpha \in \Phi(G_i, T)-\Phi(G_{i+1}, T), \lambda(\alpha)>0 \> $ in $V_i$, 
the space
$V_i^0$ to be the image 
of $ G(k) \cap \<U_\alpha(E)_{x,\frac{r_i}{2}} \, | \, \alpha \in \Phi(G_i, T)-\Phi(G_{i+1}, T), \lambda(\alpha)=0 \> $ in $V_i$, 
and 
$V_i^-$ to be the image 
of $ G(k) \cap \<U_\alpha(E)_{x,\frac{r_i}{2}} \, | \, \alpha \in \Phi(G_i, T)-\Phi(G_{i+1}, T), \lambda(\alpha)<0 \> $ in $V_i$.
Then $V_i=V_i^+ \oplus V_i^{0} \oplus V_i^-$, the subspaces $V_i^+$ and $V_i^-$ are both totally isotropic, the orthogonal complement of $V_i^+$ is $V_i^+ \oplus V_i^0$, and $V_i^0$ is a non-degenerate subspace of $V_i$. 
Let $P_i \subset \Sp(V_i)$ be the (maximal) parabolic subgroup of $\Sp(V_i)$ that preserves the subspace $V_i^+$. Note that the image of $\Uff$ in $\Sp(V_i)$ is contained in $P_i$. Let $\Uiff$ be the image 
of 
$$U_i:= G(k) \cap \<U_\alpha(E)_{x,\frac{r_i}{2}} \, | \, \alpha \in \Phi(G_i, T)-\Phi(G_{i+1}, T), \lambda(\alpha) >0 \> $$
 in the Heisenberg group $(G_i)_{x,r_i,\frac{r_i}{2}}/\left((G_i)_{x,r_i,\frac{r_i}{2}+} \cap \ker(\hat \phi_i)\right)$. Then by Yu's construction of the special isomorphism $$j_i:(G_i)_{x,r_i,\frac{r_i}{2}}/\left((G_i)_{x,r_i,\frac{r_i}{2}+} \cap \ker(\hat \phi_i)\right) \ra V_i^\sharp$$
 in \cite[Proposition~11.4]{Yu}, where
 $V_i^\sharp$ is the group $V_i \ltimes \bF_p$ with group law $(v,a).(v',a')=(v+v', a+a'+\frac{1}{2}(v,v')_i)$, and since $\lambda(\bG_m) \subset T$,
  we have $j_i(U_{i,\ff})=V_i^+ \ltimes 0$.
By \cite[Theorem~2.4.(b)]{Gerardin}\footnote{As Loren Spice pointed out, the statement of \cite[Theorem~2.4.(b)]{Gerardin} contains a typo. From the proof provided by \cite{Gerardin} one can deduce that the stated representation of $P(E_+,j)H(E_+^\perp,j)$ (i.e. the pull-back to $P(E_+,j)H(E_+^\perp,j)$ of a representation of $SH(E_0,j_0)$ as in part (a')) should be tensored with $\chi^{E_+} \ltimes 1$ before inducing it to $P(E_+,j)H(E,j)$ in order to define $\pi_+$ (using the notation of \cite{Gerardin}).}   the restriction of the Weil--Heisenberg representation $V_{\omega_i}$ (via $j_i^{-1}$) to $P_i \ltimes \Uiff$ contains a subrepresentation $V_{\omega_i}'$ on which $\Uiff$ acts trivially and on which the action of $P_i$ is as follows: By \cite[Lemma~2.3.(c)]{Gerardin} there exist surjections $p_i^1: P_i \twoheadrightarrow \GL(V_i^+)$ and $p_i^2: P_i \twoheadrightarrow \Sp(V_i^0)$. Then the action of $P_i$ on $V_{\omega_i}'$  is the tensor product of $p_i^1$ composed with a (quadratic) character $\chi$ of $\GL(V_i^+)$ and $p_i^2$ composed with a Weil representation of $\Sp(V_i^0)$.
Note that the image of $\Uff$ in $\GL(V_i^+)$ (by composing $\Uff \rightarrow P_i$ with $p_i^1: P_i \twoheadrightarrow \GL(V_i^+)$) is unipotent and hence contained in the commutator subgroup of $\GL(V_i^+)$. Thus $\chi \circ p_i^1$ is trivial on the image of $\Uff$. Moreover, the image of $\Uff$ in $\Sp(V_i^0)$ (by composing $\Uff \rightarrow P_i$ with the surjection $p_i^2: P_i \twoheadrightarrow \Sp(V_i^0)$) is contained in a minimal parabolic subgroup of $\Sp(V_i^0)$  (\cite[3.7.~Corollaire]{Borel-Tits-unipotent}) and hence also in a parabolic subgroup $P_i^0$ of $\Sp(V_i^0)$ that fixes a maximal totally isotropic subspace of $V_i^0$. By \cite[Theorem~2.4.(b)]{Gerardin} the Weil representation $V_{\omega_i}'$ restricted to $P_i^0$ contains a one dimensional subrepresentation $V_{\omega_i}''$ on which the action of $P_i^0$ factors through a character of $P_i^0/U(P_i^0)$ where $U(P_i^0)$  denotes the unipotent radical  of $P_i^0$.  
Since the image of $\Uff$ is unipotent and hence its image in $P_i^0/U(P_i^0)$ (which is isomorphic to a general linear group)
is contained in the commutator subgroup of  $P_i^0/U(P_i^0)$, the group $\Uff$ acts trivially on $V_{\omega_i}''$.

Let $V_{\kappa}''$ denote the subspace $\otimes_{1 \leq i \leq n} V_{\omega_i}''$ of $\otimes_{1 \leq i \leq n} V_{\omega_i}=V_\kappa$. Let $U_{n+1}^H$ be the preimage of $\Uff$ in $(H_{n+1})_{x,0}$ under the surjection $(H_{n+1})_{x,0} \twoheadrightarrow (H_{n+1})_{x,0}/(H_{n+1})_{x,0+}$. Since $\phi_i$ is trivial on $H_{n+1}(k)$ for all $1 \leq i \leq n$ (Lemma \ref{Lemma-phi-properties}\eqref{item-character-extension-0}), the action of the group $U_{n+1}^H$ via $\rho \otimes \kappa_{\vec \phi}$ on the subspace $V_{\rho'}'' \otimes V_{\kappa}''$ of $V_\rho \otimes V_\kappa$ is the trivial action. Moreover, recall that the restriction of $(\kappa_{\vec \phi}, V_\kappa)$ to $(G_i)_{x,r_i,\frac{r_i}{2}}$ for $1  \leq i \leq n$ is given by letting $(G_i)_{x,r_i,\frac{r_i}{2}}$ act via the Heisenberg representation $\omega_i$ on $V_{\omega_i}$ and via $\hat \phi_j|_{(G_i)_{x,r_i,\frac{r_i}{2}}}$ on $V_{\omega_j}$ for $j \neq i$. By Lemma \ref{Lemma-phi-properties}\eqref{item-character-extension-0} and the definition of $\hat \phi_j$, the character $\hat \phi_j$ is trivial on ${(H_i)_{x,r_i,\frac{r_i}{2}}}$ for $j \neq i$. Hence $U_i$ (which is contained in $(H_i)_{x,r_i,\frac{r_i}{2}}$) acts trivially via $\rho \otimes \kappa_{\vec \phi}$ on $V_{\rho'}'' \otimes V_{\kappa}''$.

If $\eps>0$ is sufficiently small, then we have 
$$(H_{n+1})_{x+\eps \lambda,0+}\subset \<(H_{n+1})_{x,0+}, U_{n+1}^H\> \quad \text{ and } \quad (H_i)_{x+\eps \lambda,r_i,\frac{r_i}{2}+}\subset\<(H_i)_{x,r_i,\frac{r_i}{2}+},U_i\>$$
 for $1 \leq i \leq n$, where $x + \eps \lambda$ arises from the action of $\eps \lambda$ on $x \in \sA(T)$. Since $\rho \otimes \kappa_{\vec \phi}$ is by the definition of $\rho$ isomorphic to a $K$-subrepresentation of $(\pi|_K, V_\pi)$, we obtain a non-trivial subspace $V''$ of $V_\pi$ on which $(H_i)_{x+\eps \lambda, r_i, \frac{r_i}{2}+}/(H_i)_{x+\eps \lambda,r_i+} \simeq (\fh_i)_{x+\eps \lambda, r_i, \frac{r_i}{2}+}/(\fh_i)_{x+\eps \lambda,r_i+}$  acts via $\varphi \circ X_i$ for $1 \leq i \leq n$ and that is fixed by $(H_{n+1})_{x+\eps \lambda,0+}$. Since $x+\eps\lambda \in \sA(T)$, the tuple $(x+\eps\lambda,(X_i)_{1 \leq i \leq n})$ is a \truncateddatum/ by Lemma \ref{Lemma-truncateddatum-with-y}, and by the same arguments as in Case 2 of the proof of Theorem \ref{Thm-existence-of-datum}, we can extend it to a \datum/ $(x+\eps\lambda,(X_i)_{1 \leq i \leq n}, (\repzero',\Vrepzero'))$ contained in $(\pi, V_\pi)$. 
 However, since $\Uff$ was non-trivial (and $\eps>0$ sufficiently small), the dimension of the facet of $\sB(G_{n+1},k)$ that contains $x+\eps\lambda$ is larger than the dimension of the facet of $\sB(G_{n+1},k)$ that contains $x$. This is a contradiction to the choice of $(x, (X_i)_{1 \leq i \leq n}, (\repzero,\Vrepzero))$, i.e. to the assumption that $(x, (X_i)_{1 \leq i \leq n}, (\repzero, \Vrepzero))$ is a maximal \datum/ for $(\pi, V_\pi)$ (as in Definition \ref{Def-datum-for}).
 \qed

In order to prove that $(\pi,V_\pi)$ contains a type as constructed by Kim and Yu in \cite{KimYu}, we introduce some additional notation following \cite[2.4]{KimYu}. We denote by $Z_s(M_{n+1})$ the maximal split torus in the center of $M_{n+1}$ and by $M_i$ the centralizer of $Z_s(M_{n+1})$ in $G_i$ for $1 \leq i \leq n$. We say (compare \cite[3.5.~Definition]{KimYu}) that the resulting commutative diagram of embeddings  (where the embeddings are chosen as explained in Remark \ref{Rem-BT})
\begin{eqnarray} \label{diagram-generic}  
\xymatrix{
	\sB(M_{n+1},k) \ar@{^{(}->}[r] \ar@{^{(}->}[d] &   \sB(M_{n},k) \ar@{^{(}->}[r] \ar@{^{(}->}[d] & \hdots \ar@{^{(}->}[r] \ar@{^{(}->}[d] & \sB(M_1,k)  \ar@{^{(}->}[d] \\
	\sB(G_{n+1},k) \ar@{^{(}->}[r] &   \sB(G_{n},k) \ar@{^{(}->}[r] &  \hdots \ar@{^{(}->}[r] & \sB(G_1,k) \\
} 
\end{eqnarray} 
is
\textit{$(\frac{r_{n+1}}{2}, \frac{r_n}{2}, \hdots, \frac{r_1}{2})$-generic relative to $x$} if 
$$ \sum_{i=1}^{n} \left(\dim((G_i)_{x,\frac{r_i}{2}}/(G_i)_{x,\frac{r_i}{2}+}) - \dim((M_i)_{x,\frac{r_i}{2}}/(M_i)_{x,\frac{r_i}{2}+}) \right)=0 ,  $$
where we recall that $r_{n+1}=0$.
Note that this property is independent of the choice of embeddings in Diagram \eqref{diagram-generic}\footnote{While our point $x$ is a point of $\sB(G,k)$ that is viewed as a point of $\sB(M_{i},k)$ and $\sB(G_i,k)$ via the above embeddings, Kim and Yu (\cite{KimYu}) fix a point in $\sB(M_{n+1},k)$ and consider its image in $\sB(M_{i},k)$ and $\sB(G_i,k)$. Hence the genericity property in  \cite[3.5.~Definition]{KimYu} does depend on the embeddings.}.

\begin{Thm} \label{Thm-exhaustion-of-types}
	Let $(\pi, V_\pi)$ be a smooth irreducible representation of $G(k)$. Then $(\pi,V_\pi)$ contains one of the \type/s constructed by Kim--Yu in \cite{KimYu}.
\end{Thm}	
\Proof By Theorem \ref{Thm-existence-of-datum}, the representation $(\pi, V_\pi)$ contains a \datum/. Let $(x, (X_i)_{1 \leq i \leq n}, (\repzero, \Vrepzero))$ be a maximal \datum/ for $(\pi, V_\pi)$ such that the non-negative number $\sum_{i=1}^{n} \left(\dim((G_i)_{x,\frac{r_i}{2}}/(G_i)_{x,\frac{r_i}{2}+})\right.$ $\left.- \dim((M_i)_{x,\frac{r_i}{2}}/(M_i)_{x,\frac{r_i}{2}+}) \right)$
is minimal among all possible choices of maximal \data/ for $(\pi, V_\pi)$. Performing the constructions above (page \pageref{page-Yu-datum} and Lemma \ref{Lemma-repzeroK}) we obtain a tuple $(\vec G, x, \vec r, \repzeroK|_{K_{G_{n+1}}}, \vec \phi)$ and an associated representation $(\piK,\VpiK)=(\repzeroK \otimes \kappa_{\vec \phi}, \VrepzeroK \otimes V_\kappa)$ as constructed by Kim and Yu that is contained in $(\pi,V_\pi)$. It remains to show that  $(K, \piK)$ is a type, i.e. that all the requirements that Kim and Yu impose on the tuple $(\vec G, x, \vec r, \repzeroK|_{K_{G_{n+1}}}, \vec \phi)$ for the construction of types are satisfied. By Lemma \ref{Lemma-phi-properties}\eqref{item-character-extension-generic} and Lemma \ref{Lemma-cuspidal} it therefore remains to show that Diagram \eqref{diagram-generic} is $(\frac{r_{n+1}}{2}, \frac{r_n}{2}, \hdots, \frac{r_1}{2})$-generic relative to $x$. 
Suppose that this is not the case. Then, by \cite[3.6~Lemma~(b)]{KimYu} and the definition of the Moy--Prasad filtration, there exists $\lambda \in X_*(Z_s(M_{n+1}))$ such that if $\eps>0$ is sufficiently small, then Diagram \eqref{diagram-generic} is $(\frac{r_{n+1}}{2}, \frac{r_n}{2}, \hdots, \frac{r_1}{2})$-generic relative to $x+\eps \lambda$ and $(G_i)_{x+\eps\lambda,\frac{r_i}{2}} \subseteq (G_i)_{x,\frac{r_i}{2}}$ and  $(G_i)_{x+\eps\lambda,r_i} \subseteq (G_i)_{x,r_i}$ for $1 \leq i \leq n+1$.
Note that this implies that  $(G_{n+1})_{x+\eps\lambda,0} = (G_{n+1})_{x,0}$ and $(G_{n+1})_{x+\eps\lambda,0+} = (G_{n+1})_{x,0+}$ because $\lambda \in X_*(Z_s(M_{n+1}))$ and $(M_{n+1})_{x,0}/(M_{n+1})_{x,0+} \simeq (G_{n+1})_{x,0}/(G_{n+1})_{x,0+}$ by definition of $M_{n+1}$. Using the notation of the proof of Lemma \ref{Lemma-cuspidal} the image of $(G_i)_{x+\eps\lambda,r_i,\frac{r_i}{2}+}$ in the Heisenberg group $(G_i)_{x,r_i,\frac{r_i}{2}}/((G_i)_{x,r_i,\frac{r_i}{2}+}\cap \ker(\hat\phi_i))=j_i^{-1}(V_i \ltimes \bF_p)$ is $j_i^{-1}(V_i^+ \ltimes \bF_p)$, where $V_i^+$ is the totally isotropic subspace $((G_i)_{x+\eps\lambda,r_i,\frac{r_i}{2}+}(G_i)_{x,r_i,\frac{r_i}{2}+})/(G_i)_{x,r_i,\frac{r_i}{2}+}$ of $V_i=(G_i)_{x,r_i,\frac{r_i}{2}}/(G_i)_{x,r_i,\frac{r_i}{2}+}$. For $1 \leq i \leq n$,  let $V_{\omega_i}'$ be a subspace of the Heisenberg representation $V_{\omega_i}$ on which $V_i^+$ acts trivially, and denote by $V_\kappa'$ the subspace $\otimes_{1 \leq j \leq n} V_{\omega_j}'$ of $\otimes_{1 \leq j \leq n} V_{\omega_j}=V_\kappa$. Then the action of $(H_i)_{x+\eps\lambda,r_i,\frac{r_i}{2}+}$ on $\VrepzeroK \otimes V_\kappa'$ factors through $(H_i)_{x+\eps\lambda,r_i,\frac{r_i}{2}+}/(H_i)_{x+\eps\lambda,r_i+} \simeq (\fh_i)_{x+\eps\lambda,r_i,\frac{r_i}{2}+}/(\fh_i)_{x+\eps\lambda,r_i+} $ on which it is given by the character $\varphi \circ X_i$ for $1 \leq i \leq n$. Moreover, $(H_{n+1})_{x+\eps\lambda,0+} = (H_{n+1})_{x,0+}$ acts trivially on $\VrepzeroK \otimes V_\kappa'$. Hence (by  Lemma \ref{Lemma-truncateddatum-with-y} and the same arguments as those in Case 2 of the proof of Theorem \ref{Thm-existence-of-datum}) we obtain a maximal \datum/ $(x+\eps\lambda, (X_i)_{1 \leq i \leq n}, (\repzero', \Vrepzero'))$ for $(\pi, V_\pi)$ with $$\sum_{i=1}^{n} \left(\dim((G_i)_{x+\eps\lambda,\frac{r_i}{2}}/(G_i)_{x+\eps\lambda,\frac{r_i}{2}+}) - \dim((M_i)_{x+\eps\lambda,\frac{r_i}{2}}/(M_i)_{x+\eps\lambda,\frac{r_i}{2}+}) \right)=0. $$ This contradicts that $0 < \sum_{i=1}^{n} \left(\dim((G_i)_{x,\frac{r_i}{2}}/(G_i)_{x,\frac{r_i}{2}+}) - \dim((M_i)_{x,\frac{r_i}{2}}/(M_i)_{x,\frac{r_i}{2}+}) \right)$ was minimal among all possible choices of maximal \data/ for $(\pi, V_\pi)$. 
 \qed

\begin{Rem}
	Theorem \ref{Thm-exhaustion-of-types} has been derived by Kim and Yu (\cite[9.1~Theorem]{KimYu}) from the result about exhaustion of Yu's supercuspidal representations by Kim (\cite{Kim}) under much more restrictive assumptions than our Assumption \ref{Assumption-p-W}. First of all they require the local field $k$ to have characteristic zero, and secondly their assumption on the residual characteristic $p$ is much stronger than ours, i.e. far from optimal, see \cite[\S~3.4]{Kim}. 
\end{Rem}

\section{Exhaustion of supercuspidal representations} \label{Section-exhaustion-Yu} 
Recall that we assume throughout the paper that $G$ splits over a tame extension and $p \nmid \abs{W}$. Under these assumptions, we obtain the following corollary of Section \ref{Section-existence-of-type}.

\begin{Thm} \label{Thm-exhaustion-Yu} 
	Every smooth irreducible supercuspidal representation of $G(k)$ arises from the construction of Yu (\cite{Yu}).
\end{Thm}
\Proof
Let $(\pi, V_\pi)$ be a smooth irreducible supercuspidal representation of $G(k)$. By Section \ref{Section-existence-of-type}, in particular Theorem \ref{Thm-exhaustion-of-types},  we can associate to $(\pi, V_\pi)$  a tuple $(\vec G, x, \vec r, \repzeroK|_{K_{G_{n+1}}}, \vec \phi)$ such that $(\pi,V_\pi)$ contains the \type/ $(K, \piK)$ associated to it by Kim--Yu following Yu's construction. Let $M_{n+1}$ be the Levi subgroup of $G_{n+1}$ attached to $x$ and $G_{n+1}$ as in Section \ref{Section-existence-of-type}, page \pageref{page-Levi}. 
We recall that $Z_S(M_{n+1})$ denotes the maximal split torus of the center $Z(M_{n+1})$ of $M_{n+1}$, and that $M_1$ is the Levi subgroup of $G$ that is the centralizer $\Cent_{G}(Z_S(M_{n+1}))$ of $Z_S(M_{n+1})$ in $G$. Kim and Yu (\cite[7.5~Theorem]{KimYu}) show that the type $(K, \piK)$ is a cover of a type for the group $M_1$. Hence, since $(\pi, V_\pi)$ is supercuspidal, we have $M_1=G$. This implies that $Z_S(M_{n+1})$ is contained in the center of $G$. Hence $Z(G_{n+1})/Z(G)$ is anisotropic, where $Z(G_{n+1})$ and $Z(G)$ denote the centers of $G_{n+1}$ and $G$, respectively, and $M_{n+1}=G_{n+1}$. Instead of working with $K_{G_{n+1}}=(G_{n+1})_x$ in Section \ref{Section-existence-of-type}, we could have equally well performed all constructions for the stabilizer $(G_{n+1})_{[x]}$ of the image $[x]$ of $x$ in the reduced Bruhat--Tits building of $G_{n+1}$ (by replacing $(M_{n+1})_x$ by $(M_{n+1})_{[x]}$ everywhere) to obtain a representation $(\repzeroKYu, \VrepzeroKYu)$ of $\KYu=(G_{n+1})_{[x]}(G_{n})_{x,\frac{r_n}{2}}\hdots(G_1)_{x,\frac{r_1}{2}} $ such that
the representation $(\piKYu,\VpiKYu)$ of $\KYu$ associated to $(\vec G, x, \vec r, \repzeroKYu|_{(G_{n+1})_{[x]}}, \vec \phi)$ by Yu is contained in $(\pi|_{\KYu}, V_\pi)$. Since $M_{n+1}=\Cent_{G_{n+1}}(Z_S(M_{n+1}))=G_{n+1}$, the compactly induced representation $\ind_{(G_{n+1})_{[x]}}^{G(k)}\repzeroKYu|_{(G_{n+1})_{[x]}}$ is irreducible supercuspidal (by \cite[Proposition~6.6]{MP2}).
Hence $(\vec G, x, \vec r, \repzeroKYu|_{(G_{n+1})_{[x]}}, \vec \phi)$ satisfies all the conditions that Yu requires for his construction of supercuspidal representations (\cite[\S~3]{Yu}), and $\ind_{\KYu}^{G(k)}\piKYu$ is the corresponding irreducible supercuspidal representations (\cite[Proposition~4.6]{Yu}). By Frobenius reciprocity, we obtain a non-trivial morphism from $(\ind_{\KYu}^{G(k)}\piKYu, \ind_{\KYu}^{G(k)}\VpiKYu)$ to $(\pi, V_\pi)$, and hence these two irreducible representations are isomorphic. \qed

\begin{Rem}
	The exhaustion of supercuspidal representations by Yu's construction has been known under the assumption that $k$ has characteristic zero and $p$ is a sufficiently large prime number thanks to Kim (\cite{Kim}). We refer the reader to \cite[\S~3.4]{Kim} for the precise conditions for $p$ being ``sufficiently large''. These assumptions are much stronger than $p \nmid \abs{W}$. 
\end{Rem}

The proof of Theorem \ref{Thm-exhaustion-Yu} also shows how to recognize if a representation is supercuspidal by only considering a maximal \datum/ for this representation.

\begin{Cor} \label{Cor-supercuspidal-criterion}
	Let $(\pi, V_\pi)$ be a smooth irreducible representation of $G(k)$, and let $(x, (X_i)_{1 \leq i \leq n},$ $(\repzero, \Vrepzero))$ be a maximal \datum/ for $(\pi,V_\pi)$. Then $(\pi,V_\pi)$ is supercuspidal if and only if $x$ is a facet of minimal dimension in $\sB(G_{n+1},k)$ and $Z(G_{n+1})/Z(G)$ is anisotropic, where $G_{n+1}=\Cent_{G}(\sum_{i=1}^n X_i)$.
\end{Cor}	
\Proof
The point $x$ is a facet of minimal dimension in $\sB(G_{n+1},k)$ if and only if $M_{n+1}=G_{n+1}$. Hence we have seen in the proof of Theorem \ref{Thm-exhaustion-Yu} that $(\pi, V_\pi)$ being supercuspidal implies the other two conditions in the corollary. The proof of Theorem \ref{Thm-exhaustion-Yu} also shows that the other two conditions are sufficient to prove that $(\pi, V_\pi)$ is supercuspidal. \qed

%\newpage

\bibliography{Jessicasbib}

\begin{bibdiv}
\begin{biblist}

\bib{Adler-DeBacker}{article}{
      author={Adler, Jeffrey~D.},
      author={DeBacker, Stephen},
       title={Some applications of {B}ruhat-{T}its theory to harmonic analysis
  on the {L}ie algebra of a reductive {$p$}-adic group},
        date={2002},
        ISSN={0026-2285},
     journal={Michigan Math. J.},
      volume={50},
      number={2},
       pages={263\ndash 286},
         url={https://doi.org/10.1307/mmj/1028575734},
      review={\MR{1914065}},
}

\bib{Adler}{article}{
      author={Adler, Jeffrey~D.},
       title={Refined anisotropic {$K$}-types and supercuspidal
  representations},
        date={1998},
        ISSN={0030-8730},
     journal={Pacific J. Math.},
      volume={185},
      number={1},
       pages={1\ndash 32},
}

\bib{Adler-Roche}{article}{
      author={Adler, Jeffrey~D.},
      author={Roche, Alan},
       title={An intertwining result for {$p$}-adic groups},
        date={2000},
        ISSN={0008-414X},
     journal={Canad. J. Math.},
      volume={52},
      number={3},
       pages={449\ndash 467},
         url={http://dx.doi.org/10.4153/CJM-2000-021-8},
      review={\MR{1758228}},
}

\bib{Bernstein}{incollection}{
      author={Bernstein, Joseph~N.},
       title={Le ``centre'' de {B}ernstein},
        date={1984},
   booktitle={Representations of reductive groups over a local field},
      series={Travaux en Cours},
   publisher={Hermann, Paris},
       pages={1\ndash 32},
        note={Edited by P. Deligne},
      review={\MR{771671}},
}

\bib{BK}{book}{
      author={Bushnell, Colin~J.},
      author={Kutzko, Philip~C.},
       title={The admissible dual of {${\rm GL}(N)$} via compact open
  subgroups},
      series={Annals of Mathematics Studies},
   publisher={Princeton University Press, Princeton, NJ},
        date={1993},
      volume={129},
        ISBN={0-691-03256-4; 0-691-02114-7},
         url={https://doi.org/10.1515/9781400882496},
      review={\MR{1204652}},
}

\bib{BK-SLn}{article}{
      author={Bushnell, Colin~J.},
      author={Kutzko, Philip~C.},
       title={The admissible dual of {${\rm SL}(N)$}. {II}},
        date={1994},
        ISSN={0024-6115},
     journal={Proc. London Math. Soc. (3)},
      volume={68},
      number={2},
       pages={317\ndash 379},
         url={https://doi.org/10.1112/plms/s3-68.2.317},
      review={\MR{1253507}},
}

\bib{BK-types}{article}{
      author={Bushnell, Colin~J.},
      author={Kutzko, Philip~C.},
       title={Smooth representations of reductive {$p$}-adic groups: structure
  theory via types},
        date={1998},
        ISSN={0024-6115},
     journal={Proc. London Math. Soc. (3)},
      volume={77},
      number={3},
       pages={582\ndash 634},
         url={https://doi.org/10.1112/S0024611598000574},
      review={\MR{1643417}},
}

\bib{BK-types-exhaustion}{article}{
      author={Bushnell, Colin~J.},
      author={Kutzko, Philip~C.},
       title={Semisimple types in {${\rm GL}_n$}},
        date={1999},
        ISSN={0010-437X},
     journal={Compositio Math.},
      volume={119},
      number={1},
       pages={53\ndash 97},
         url={https://doi.org/10.1023/A:1001773929735},
      review={\MR{1711578}},
}

\bib{Borel}{book}{
      author={Borel, Armand},
       title={Linear algebraic groups},
     edition={Second},
      series={Graduate Texts in Mathematics},
   publisher={Springer-Verlag, New York},
        date={1991},
      volume={126},
        ISBN={0-387-97370-2},
}

\bib{Bourbaki-4-6}{book}{
      author={Bourbaki, Nicolas},
       title={Lie groups and {L}ie algebras. {C}hapters 4--6},
      series={Elements of Mathematics (Berlin)},
   publisher={Springer-Verlag, Berlin},
        date={2002},
        ISBN={3-540-42650-7},
         url={http://dx.doi.org/10.1007/978-3-540-89394-3},
        note={Translated from the 1968 French original by Andrew Pressley},
      review={\MR{1890629}},
}

\bib{Broussous}{article}{
      author={Broussous, Paul},
       title={Extension du formalisme de {B}ushnell et {K}utzko au cas d'une
  alg\`ebre \`a division},
        date={1998},
        ISSN={0024-6115},
     journal={Proc. London Math. Soc. (3)},
      volume={77},
      number={2},
       pages={292\ndash 326},
         url={https://doi-org.proxy.lib.umich.edu/10.1112/S0024611598000471},
      review={\MR{1635145}},
}

\bib{Borel-Tits-unipotent}{article}{
      author={Borel, Armand},
      author={Tits, Jacques},
       title={{\'E}l\'ements unipotents et sous-groupes paraboliques de groupes
  r\'eductifs. {I}},
        date={1971},
        ISSN={0020-9910},
     journal={Invent. Math.},
      volume={12},
       pages={95\ndash 104},
         url={https://doi.org/10.1007/BF01404653},
      review={\MR{0294349}},
}

\bib{Carayol}{article}{
      author={Carayol, Henri},
       title={Repr\'esentations supercuspidales de {${\rm GL}_{n}$}},
        date={1979},
        ISSN={0151-0509},
     journal={C. R. Acad. Sci. Paris S\'er. A-B},
      volume={288},
      number={1},
       pages={A17\ndash A20},
      review={\MR{522009}},
}

\bib{Pseudoreductive2}{book}{
      author={Conrad, Brian},
      author={Gabber, Ofer},
      author={Prasad, Gopal},
       title={Pseudo-reductive groups},
     edition={Second},
      series={New Mathematical Monographs},
   publisher={Cambridge University Press, Cambridge},
        date={2015},
      volume={26},
        ISBN={978-1-107-08723-1},
         url={https://doi.org/10.1017/CBO9781316092439},
      review={\MR{3362817}},
}

\bib{ConradSGA3}{incollection}{
      author={Conrad, Brian},
       title={Reductive group schemes},
        date={2014},
   booktitle={Autour des sch\'emas en groupes. {V}ol. {I}},
      series={Panor. Synth\`eses},
      volume={42/43},
   publisher={Soc. Math. France, Paris},
       pages={93\ndash 444},
      review={\MR{3362641}},
}

\bib{Fi-tame-tori}{article}{
      author={Fintzen, Jessica},
       title={{T}ame {T}ori in $p$-{A}dic {G}roups and {G}ood {S}emisimple
  {E}lements},
        date={2019},
        ISSN={1073-7928},
     journal={International Mathematics Research Notices},
  eprint={https://academic.oup.com/imrn/advance-article-pdf/doi/10.1093/imrn/rnz234/30622725/rnz234.pdf},
         url={https://doi.org/10.1093/imrn/rnz234},
        note={rnz234},
}

\bib{Fi}{article}{
      author={Fintzen, Jessica},
       title={On the {M}oy-{P}rasad filtration},
     journal={to appear in the Journal of the European Mathematical Society
  (JEMS)},
        note={Available at \url{https://arxiv.org/pdf/1511.00726.pdf}.},
}

\bib{FR}{article}{
      author={Fintzen, Jessica},
      author={Romano, Beth},
       title={Stable vectors in {M}oy-{P}rasad filtrations},
        date={2017},
        ISSN={0010-437X},
     journal={Compos. Math.},
      volume={153},
      number={2},
       pages={358\ndash 372},
         url={https://doi.org/10.1112/S0010437X16008228},
      review={\MR{3705228}},
}

\bib{Gerardin}{article}{
      author={G\'erardin, Paul},
       title={Weil representations associated to finite fields},
        date={1977},
        ISSN={0021-8693},
     journal={J. Algebra},
      volume={46},
      number={1},
       pages={54\ndash 101},
         url={https://doi.org/10.1016/0021-8693(77)90394-5},
      review={\MR{0460477}},
}

\bib{Goldberg-Roche}{article}{
      author={Goldberg, David},
      author={Roche, Alan},
       title={Types in {${\rm SL}_n$}},
        date={2002},
        ISSN={0024-6115},
     journal={Proc. London Math. Soc. (3)},
      volume={85},
      number={1},
       pages={119\ndash 138},
         url={https://doi.org/10.1112/S002461150201359X},
      review={\MR{1901371}},
}

\bib{Hakim-Murnaghan}{article}{
      author={Hakim, Jeffrey},
      author={Murnaghan, Fiona},
       title={Distinguished tame supercuspidal representations},
        date={2008},
        ISSN={1687-3017},
     journal={Int. Math. Res. Pap. IMRP},
      number={2},
       pages={Art. ID rpn005, 166},
      review={\MR{2431732}},
}

\bib{Howe}{article}{
      author={Howe, Roger~E.},
       title={Tamely ramified supercuspidal representations of {${\rm
  Gl}_{n}$}},
        date={1977},
        ISSN={0030-8730},
     journal={Pacific J. Math.},
      volume={73},
      number={2},
       pages={437\ndash 460},
         url={http://projecteuclid.org/euclid.pjm/1102810618},
      review={\MR{0492087}},
}

\bib{Kaletha}{article}{
      author={Kaletha, Tasho},
       title={Regular supercuspidal representations},
        date={2019},
        ISSN={0894-0347},
     journal={J. Amer. Math. Soc.},
      volume={32},
      number={4},
       pages={1071\ndash 1170},
         url={https://doi.org/10.1090/jams/925},
      review={\MR{4013740}},
}

\bib{Kempf}{article}{
      author={Kempf, George~R.},
       title={Instability in invariant theory},
        date={1978},
        ISSN={0003-486X},
     journal={Ann. of Math. (2)},
      volume={108},
      number={2},
       pages={299\ndash 316},
         url={http://dx.doi.org.proxy.lib.umich.edu/10.2307/1971168},
      review={\MR{506989}},
}

\bib{Kim}{article}{
      author={Kim, Ju-Lee},
       title={Supercuspidal representations: an exhaustion theorem},
        date={2007},
        ISSN={0894-0347},
     journal={J. Amer. Math. Soc.},
      volume={20},
      number={2},
       pages={273\ndash 320 (electronic)},
}

\bib{Kim-Murnaghan}{article}{
      author={Kim, Ju-Lee},
      author={Murnaghan, Fiona},
       title={Character expansions and unrefined minimal {$K$}-types},
        date={2003},
        ISSN={0002-9327},
     journal={Amer. J. Math.},
      volume={125},
      number={6},
       pages={1199\ndash 1234},
  url={http://muse.jhu.edu/journals/american_journal_of_mathematics/v125/125.6kim.pdf},
      review={\MR{2018660}},
}

\bib{KimYu}{incollection}{
      author={Kim, Ju-Lee},
      author={Yu, Jiu-Kang},
       title={Construction of tame types},
        date={2017},
   booktitle={Representation theory, number theory, and invariant theory},
      series={Progr. Math.},
      volume={323},
   publisher={Birkh\"auser/Springer, Cham},
       pages={337\ndash 357},
      review={\MR{3753917}},
}

\bib{Morris}{article}{
      author={Morris, Lawrence},
       title={Tamely ramified intertwining algebras},
        date={1993},
        ISSN={0020-9910},
     journal={Invent. Math.},
      volume={114},
      number={1},
       pages={1\ndash 54},
         url={https://doi-org.proxy.lib.umich.edu/10.1007/BF01232662},
      review={\MR{1235019}},
}

\bib{Morris-depth-zero}{article}{
      author={Morris, Lawrence},
       title={Level zero {$\bf G$}-types},
        date={1999},
        ISSN={0010-437X},
     journal={Compositio Math.},
      volume={118},
      number={2},
       pages={135\ndash 157},
         url={https://doi-org.proxy.lib.umich.edu/10.1023/A:1001019027614},
      review={\MR{1713308}},
}

\bib{Moy-exhaustion}{article}{
      author={Moy, Allen},
       title={Local constants and the tame {L}anglands correspondence},
        date={1986},
        ISSN={0002-9327},
     journal={Amer. J. Math.},
      volume={108},
      number={4},
       pages={863\ndash 930},
         url={https://doi.org/10.2307/2374518},
      review={\MR{853218}},
}

\bib{MP1}{article}{
      author={Moy, Allen},
      author={Prasad, Gopal},
       title={Unrefined minimal {$K$}-types for {$p$}-adic groups},
        date={1994},
        ISSN={0020-9910},
     journal={Invent. Math.},
      volume={116},
      number={1-3},
       pages={393\ndash 408},
}

\bib{MP2}{article}{
      author={Moy, Allen},
      author={Prasad, Gopal},
       title={Jacquet functors and unrefined minimal {$K$}-types},
        date={1996},
        ISSN={0010-2571},
     journal={Comment. Math. Helv.},
      volume={71},
      number={1},
       pages={98\ndash 121},
}

\bib{Miyauchi-Stevens}{article}{
      author={Miyauchi, Michitaka},
      author={Stevens, Shaun},
       title={Semisimple types for {$p$}-adic classical groups},
        date={2014},
        ISSN={0025-5831},
     journal={Math. Ann.},
      volume={358},
      number={1-2},
       pages={257\ndash 288},
         url={https://doi.org/10.1007/s00208-013-0953-y},
      review={\MR{3157998}},
}

\bib{ReederYu}{article}{
      author={Reeder, Mark},
      author={Yu, Jiu-Kang},
       title={Epipelagic representations and invariant theory},
        date={2014},
        ISSN={0894-0347},
     journal={J. Amer. Math. Soc.},
      volume={27},
      number={2},
       pages={437\ndash 477},
         url={http://dx.doi.org/10.1090/S0894-0347-2013-00780-8},
      review={\MR{3164986}},
}

\bib{Secherre-Stevens}{article}{
      author={S\'echerre, Vincent},
      author={Stevens, Shaun},
       title={Repr\'esentations lisses de {${\rm GL}_m(D)$}. {IV}.
  {R}epr\'esentations supercuspidales},
        date={2008},
        ISSN={1474-7480},
     journal={J. Inst. Math. Jussieu},
      volume={7},
      number={3},
       pages={527\ndash 574},
         url={https://doi-org.proxy.lib.umich.edu/10.1017/S1474748008000078},
      review={\MR{2427423}},
}

\bib{Secherre-Stevens-types}{article}{
      author={S\'echerre, Vincent},
      author={Stevens, Shaun},
       title={Smooth representations of {$GL_m(D)$} {VI}: semisimple types},
        date={2012},
        ISSN={1073-7928},
     journal={Int. Math. Res. Not. IMRN},
      number={13},
       pages={2994\ndash 3039},
         url={https://doi.org/10.1093/imrn/rnr122},
      review={\MR{2946230}},
}

\bib{Springer-Steinberg}{incollection}{
      author={Springer, Tonny~Albert},
      author={Steinberg, Robert},
       title={Conjugacy classes},
        date={1970},
   booktitle={Seminar on {A}lgebraic {G}roups and {R}elated {F}inite {G}roups
  ({T}he {I}nstitute for {A}dvanced {S}tudy, {P}rinceton, {N}.{J}., 1968/69)},
      series={Lecture Notes in Mathematics, Vol. 131},
   publisher={Springer, Berlin},
       pages={167\ndash 266},
      review={\MR{0268192}},
}

\bib{Stevens}{article}{
      author={Stevens, Shaun},
       title={The supercuspidal representations of {$p$}-adic classical
  groups},
        date={2008},
        ISSN={0020-9910},
     journal={Invent. Math.},
      volume={172},
      number={2},
       pages={289\ndash 352},
         url={http://dx.doi.org.proxy.lib.umich.edu/10.1007/s00222-007-0099-1},
      review={\MR{2390287}},
}

\bib{Steinberg-torsion}{article}{
      author={Steinberg, Robert},
       title={Torsion in reductive groups},
        date={1975},
        ISSN={0001-8708},
     journal={Advances in Math.},
      volume={15},
       pages={63\ndash 92},
         url={https://doi.org/10.1016/0001-8708(75)90125-5},
      review={\MR{0354892}},
}

\bib{Yu}{article}{
      author={Yu, Jiu-Kang},
       title={Construction of tame supercuspidal representations},
        date={2001},
        ISSN={0894-0347},
     journal={J. Amer. Math. Soc.},
      volume={14},
      number={3},
       pages={579\ndash 622 (electronic)},
}

\bib{Zink}{article}{
      author={Zink, Ernst-Wilhelm},
       title={Representation theory of local division algebras},
        date={1992},
        ISSN={0075-4102},
     journal={J. Reine Angew. Math.},
      volume={428},
       pages={1\ndash 44},
         url={https://doi-org.proxy.lib.umich.edu/10.1515/crll.1992.428.1},
      review={\MR{1166506}},
}

\end{biblist}
\end{bibdiv}

\end{document}